\documentclass[aap,MSNbibl,seceqn,citesort,dvips]{arximspdf}
\usepackage{mathrsfs}

\doi{10.1214/11-AAP839} %kopijuoti is PTS
\volume{22}
\issue{6}
\pubyear{2012}
\firstpage{2460}
\lastpage{2504}

\makeatletter

\renewcommand{\epsilon}{\varepsilon}

\newcommand{\wwidetilde}[1]{\hspace*{0.8pt}\widetilde{\hspace*{-0.8pt}\widetilde{#1}}}
\newcommand{\widehattilde}[1]{\hspace*{0.8pt}\widehat{\hspace*{-0.8pt}\widetilde{#1}}}

\newcommand{\llangle}{\langle\!\langle}
\newcommand{\ggangle}{\rangle\!\rangle}

\newtheorem{theorem}{Theorem}[section]
\newtheorem{prop}{Proposition}[section]
\newproclaim{property}{Property}[section]

\newtheorem{lemma}{Lemma}[section]
\newtheorem{cor}{Corollary}[section]
\newproclaim{rmq}{Remark}[section]

\newcommand{\reps}[1]{#1^{R}_{\epsilon}}

\renewcommand{\Re}{\operatorname{\mathcal{R}\mathrm{e}}}
\newcommand{\Vvert}{|\!|\!|}

\newcommand{\eps}[1]{#1_{\epsilon}}
\newcommand{\mol}[1]{#1_{\eta}}

\newcommand{\sigb}{\bolds{\sigma}}

\makeatother

\begin{document}
\begin{frontmatter}

\title{A diffusion approximation theorem for a
nonlinear PDE with application to random birefringent optical
fibers\thanksref{T1}}
\runtitle{Diffusion approximation theorem}

\thankstext{T1}{Supported in part by the French ANR research Project
ANR-07-44 BLAN-0250 and by a grant from R\'{e}gion Ile-de-France.}

\begin{aug}
\author[A]{\fnms{A.} \snm{de Bouard}\ead[label=e1]{debouard@cmap.polytechnique.fr}}
\and
\author[A]{\fnms{M.} \snm{Gazeau}\corref{}\ead[label=e2]{gazeau@cmap.polytechnique.fr}}
\runauthor{A. de Bouard and M. Gazeau}
\affiliation{CMAP, CNRS and \'Ecole Polytechnique}
\address[A]{CMAP, CNRS UMR 7641\\
\'Ecole Polytechnique\\
91128 Palaiseau Cedex\\
France\\
\printead{e1}\\
\phantom{E-mail: }\printead*{e2}} %adresu isvedimo komanda gale!
\end{aug}

% HISTORY:
\received{\smonth{5} \syear{2011}}
\revised{\smonth{11} \syear{2011}}

% ABSTRACT
%
\begin{abstract}
In this article we propose a generalization of the theory of diffusion
approximation for random ODE to a nonlinear system of random
Schr\"{o}dinger equations. This system arises in the study of pulse
propagation in randomly birefringent optical fibers. We first show
existence and uniqueness of solutions for the random PDE and the
limiting equation. We follow the work of Garnier and Marty
[\textit{Wave Motion} \textbf{43} (2006) 544--560],
Marty [Probl\`emes d'\'evolution en milieux al\'eatoires: Th\'eor\`emes
limites, sch\'emas num\'eriques et applications en optique (2005)
Univ. Paul Sabatier], where a linear electric field is considered, and we
get an asymptotic dynamic for the nonlinear electric field.
\end{abstract}

% KEYWORDS
%

\pdfkeywords{35Q55, 60H15, Nonlinear Schrodinger equation,
stochastic partial differential equations,
white noise, diffusion limit}

\end{frontmatter}

%s1 #&#
\section{Introduction}
The Manakov PMD equation has been introduced by Wai and Menyuk in \cite
{3} to
study light propagation over long distance in random birefringent
optical fibers.
Due to the various length scales present in this problem, a small parameter
$\epsilon$ appears in the rescaled equation. Our aim in this paper is
to prove a
diffusion limit theorem for this equation for which we will have to
generalize the
perturbed test function method~\cite{papa,kushner,papanicolaou}
to the case
of infinite dimension. In~\cite{marty,2}, a limit theorem is
proved for the
linear part of the Manakov PMD equation using the Fourier transform and
the theory
of diffusion approximation for random ODE. Obviously the method in
\cite
{marty,2} does not work for a nonlinear PDE. In~\cite{bD,2}, a limit
theorem is
proved for a nonlinear scalar PDE driven by a one-dimensional noise.
The proof relies on the fact that the solution processes are continuous
functions
of the noise. These methods are no longer applicable to the limit
equation that we
will consider which is driven by a three-dimensional noise, because the solution
cannot be written as a continuous function of the noise. Indeed, in a general
setting a strong solution of a stochastic equation is only a measurable
function of
the initial data and the Brownian motion driving the equation. However,
in the case
of a one-dimensional noise, Doss~\cite{Doss} and Sussman \cite
{Sussman} proved
that the solution of such an equation can be written as a continuous
function of
the Brownian motion. This result has been extended by Yamato \cite
{yamato} to
multidimensional Brownian motions when the Lie algebra generated by the vector
fields of the equation is nilpopent of step~$p$. He actually proves the
equivalence
between the nilpotent hypothesis and the fact that the solution can be
written as
a continuous function of iterated Stratonovich integrals. In our case
the vector
fields driving the Manakov PMD equation are functions of the Pauli
matrices and the
nilpotent hypothesis of Yamato is not satisfied. This motivates the use
of the
perturbed test function method. Note that the method has been used for
a linear PDE
in~\cite{vovelle} and a PDE with bounded diffusion coefficients in
\cite{piatnitski}.

We are also interested in the mathematical analysis of both the Manakov PMD
and the limit equations. Using a unitary transformation, we are able to
establish
Strichartz estimates for the transformed equation, that are not
available for the
Manakov PMD equation. This result will then enable us to prove global existence
of solutions. The limiting equation is also studied. We use a
compactness method
to study the existence and uniqueness of solutions of this latter
equation, due to
the lack of nilpotent hypothesis and to the absence of unitary
transformation similar to the Manakov PMD case.
%We think that the generalized perturbed test function method would
%apply to more general PDEs
%(\nu(t), \eps{X}),
%where $G$ is a bounded function of the noise.

%s1.1 #&#
\subsection{Presentation of the model}

Optical fibers are thin, transparent and flexible fibers along which
the light propagates to transmit information over long distances and so
are of huge interest in modern communications. In a perfect fiber, the
two transverse components of the electric field are degenerate in the
sense that they propagate with the same characteristics: group
velocity, chromatic dispersion, refractive indices $(n_1=n_2)$, etc.
However, during the fabrication process the fiber may present defects
like an ellipticity of the core or suffer from mechanical distortions
like stress constraints or twisting~\cite{Agrawalappli,Agrawal}. These
phenomena induce modal birefringence $(n_1 \neq n_2)$ characterized by
an orientation angle $\theta$ and an amplitude $b$. If $n_1 >n_2$, we
then define a slow axe and a rapid axe corresponding, respectively, to
the mode indices $n_1$ and $n_2$. The orientation angle $\theta$
describes the rotation of the local polarization axes with respect to
the initial axes. The birefringence strength (or degree of modal
birefringence) is given by $b = |n_1-n_2| k_0 = k_1 -k_2$, where $k_1,
k_2 $ are the components of the wave vector and $k_0$ the wavenumber of
the incident light in vacuum. The beat\vspace*{1pt} length $L_B=\frac
{2\pi}{k_1 -k_2}$ indicates the length required for the polarization to
return to its initial state. There exist several types of birefringence
that do not have the same effect on the electric field. Usually linear
birefringence is studied (in the absence of Kerr effect, a linearly
polarized light remains linearly polarized), although it has been shown
that the birefringence could also be elliptic (occurring in case of
twisting, see Menyuk~\cite{5}). In case of a uniform anisotropy along
the fiber, the birefringence parameters $(\theta, b)$ are constant.
However, in realistic configurations, the anisotropy is not uniform
along the fiber. We assume, as in~\cite{31,4,1,3}, that the
birefringence is randomly varying, implying polarization mode
dispersion (PMD). The difference of velocity of the two modes, due to
random change of the birefringence (and so of the refractive indices),
induces coupling between the two polarized modes and pulse spreading:
PMD is one of the limiting factors of high bit rate transmission.

In~\cite{3}, Wai and Menyuk assumed that there is no
polarization-dependent loss and considered that communication fibers
are nearly linearly birefringent. We here use one of the models
introduced in~\cite{3} for which the local axes of birefringence are
bended with an angle $\theta$ randomly varying along the propagation
axe and that $b$ and $b'$ (the frequency derivative of $b$) are
constant along this axe.
Let us recall that the Pauli matrices are defined by
\[
\sigma_1=\pmatrix{
0 & 1\cr
1 & 0},\qquad
\sigma_2=\pmatrix{
0 & -i\cr
i & 0},\qquad
\sigma_3=\pmatrix{
1 & 0\cr
0 & -1},
\]
and let us consider the coupled nonlinear Schr\"{o}dinger equation
transformed into the frame of the local axes of birefringence \cite
{menyukappli,3}
%
%e1.1 #&#
\begin{eqnarray}\label{systemlocalaxe}
&&
i\,\frac{\partial\Psi}{\partial t}+ \widetilde{\Sigma}(t)\Psi+ ib'
\sigma_3\,\frac{\partial\Psi}{\partial x}+\frac{d_0}{2}\,\frac
{\partial^2
\Psi}{\partial x^2}\nonumber\\[-8pt]\\[-8pt]
&&\qquad{}+\frac{5}{6}|\Psi|^2\Psi+\frac
{1}{6}(\Psi^*\sigma_3\Psi)\sigma_3\Psi+\frac
{1}{3}N(
\Psi)=0,\nonumber
\end{eqnarray}
where $d_0$ is the group velocity dispersion parameter, $N( \Psi
)= (
\overline{\Psi}_1\Psi_2^2,
\overline{\Psi}_2\Psi_1^2
)^t$ and
\[
\widetilde{\Sigma}(t)= \pmatrix{
b & \displaystyle -\frac{i}{2}\,\frac{d\theta(t)}{dt}\vspace*{2pt}\cr
\displaystyle \frac{i}{2}\,\frac{d\theta(t)}{dt} & -b}.
\]
We recall that in the context of fiber optics, $x$ corresponds to the
retarded time while $t$ corresponds to the distance along the fiber. We
introduce\vspace*{2pt} a new vector field $\widetilde{\Psi}= \exp(-ibt
\sigma_3)\Psi$. The evolution\vspace*{1pt} of $\widetilde{\Psi}$ is
given by the previous equation (\ref{systemlocalaxe}) replacing
$\widetilde {\Sigma}$ and $N( \Psi)$, respectively, by
\[
\widetilde{\widetilde{\Sigma}}(t)= \pmatrix{
0 & \displaystyle -\frac{i}{2}\,\frac{d\theta(t)}{dt}e^{-2ibt}\vspace*{2pt}\cr
\displaystyle \frac{i}{2}\,\frac{d\theta(t)}{dt}e^{2ibt} & 0}
\quad\mbox{and}\quad
N( \widetilde{\Psi})=\pmatrix{
\displaystyle \overline{\widetilde{\Psi}}_1\widetilde{\Psi}_2^2
e^{-4ibt}\vspace*{2pt}\cr
\displaystyle \overline{\widetilde{\Psi}}_2\widetilde{\Psi}_1^2e^{4ibt}}.
\]
Following Wai and Menyuk~\cite{menyukappli,4,1,3}, we
denote by $l$ the fiber length. We also denote by $l_d$ the dispersion
length scale and $l_{nl}$ the nonlinear length scale related to Kerr
effect. The fiber autocorrelation length $l_c$ is the length over which
two polarization components remain correlated. We consider, as in \cite
{3}, a typical configuration where $l \sim l_d \sim l_{nl} \gg l_c
\gg L_B$, that is, we consider ``relatively small'' propagation
distances. Under these assumptions and the assumptions on $d\theta/dt$
below, the term $N( \widetilde{\Psi})$ is rapidly
oscillating and will be neglected~\cite{Agrawal,menyukappli,3}, its effect being averaged out to zero. As in \cite
{menyukappli,3}, we introduce a unitary matrix
%
%e1.2 #&#
\begin{equation}\label{defT}
T(t)=\pmatrix{
u_1(t) & \overline{u}_2(t)\cr
-u_2(t) & \overline{u}_1(t)},
\end{equation}
the solution of
%
%e1.3 #&#
\begin{equation}\label{EDOT}
i\,\frac{\partial T(t)}{\partial t}+ \widetilde{\widetilde{\Sigma}}(t)T(t)=0.
\end{equation}
We also consider, for $t \in\mathbb{R}_+$, the matrix
%
%e1.4 #&#
\begin{eqnarray}\label{matrixsigb}\qquad
\sigb( u(t)) &= &\pmatrix{
|u_1|^2-|u_2|^2 & 2\overline{u}_1\overline{u}_2\cr
2u_1u_2 & |u_2|^2-|u_1|^2 }
= \pmatrix{
m_3 & m_1-im_2\cr
m_1+im_2 & -m_3}
\nonumber\\[-8pt]\\[-8pt]
&=&\sigma_1m_1(t)+\sigma_2m_2(t)+\sigma_3m_3(t),\nonumber
\end{eqnarray}
which characterizes the linear birefringence and where $m_1, m_2, m_3$
are real-valued processes.
%and the coefficients
%a_1=\abs{u_1}^2-\abs{u_2}^2,
%a_2=2u_1u_2.
Then we can remove the rapid variation of the state of polarization in
the evolution of $\widetilde{\Psi}$ using the change of variable
$\widetilde{\Psi}(t)= T(t)X(t)$. We obtain
%
%e1.5 #&#
\begin{eqnarray}\label{system3u}
&&
i\,\frac{\partial X}{\partial t}+ib'\sigb( u(t)) \,\frac
{\partial X}{\partial x}+\frac{d_0}{2}\,\frac{\partial^2 X}{\partial
x^2}\nonumber\\[-8pt]\\[-8pt]
&&\qquad{}+ \frac{5}{6}|X|^2X +\frac{1}{6}(X^*\sigma_3X
)\sigma
_3X+\frac{1}{6}N_{u}(X )=0,\nonumber
\end{eqnarray}
where $N_{u}(X )=(N_{1, u}(X ),N_{2,
u}
(X ))^t$ satisfy
%
%e1.6 #&#
%e1.7 #&#
\begin{eqnarray}
\label{defN1}
N_{1, u}(X )&=&
(m_1^2+m_2^2)(2|X_2|^2-|X_1|^2)X_1 \nonumber\\
&&{}+
(m_1-im_2) m_3(2|X_1|^2-|X_2|^2)X_2 \\
&&{} +
(m_1-im_2)^2X_2^2\overline{X}_1+ (m_1+im_2
)m_3X_1^2\overline{X}_2,\nonumber\\
\label{defN2}
N_{2, u}(X )&=& (m_1^2+m_2^2
)(2|X_1|^2-|X_2|^2)X_2 \nonumber\\
&&{}-(m_1+im_2) m_3(2|X_2|^2-|X_1|^2)X_1
\\
&&{}- (m_1-im_2)m_3X_2^2\overline{X}_1+
(m_1+im_2)^2X_1^2\overline{X}_2.\nonumber
\end{eqnarray}
Assuming, as in~\cite{marty,2,1}, that the correlation length
of $d\theta/dt$ is much shorter than the birefringence beat length and
that $|d\theta/dt| \ll b$, we set $d\theta/dt=2\varepsilon_0
\alpha
(t)$, where $\varepsilon_0$ is a small dimensionless parameter and
$\alpha
$ a Markov process with good ergodic properties. Thus, we may replace
the process $u$ by $\nu$, with~\cite{marty,2}
%
%e1.8 #&#
\begin{eqnarray}\label{nu1}
d\nu(t)&=&i\sqrt{\gamma_c}\bigl(\sigma_1\nu(t)\circ
dW_1(t)+\sigma_2\nu
(t)\circ dW_2(t)\bigr)+i\gamma_s\sigma_3\nu(t)\,dt\nonumber\\
&=&i\sqrt
{\gamma_c}\bigl(\sigma_1\nu(t)\,dW_1(t)+\sigma_2\nu(t)\,dW_2(t)
\bigr)+i\gamma_s\sigma_3\nu(t)\,dt\\
&&{}-\gamma_c\nu(t)\,dt,\nonumber
\end{eqnarray}
where $|\nu_1(0)|^2+|\nu_2(0)|^2=1$, $W=
(W_1,W_2)$
is a $2d$ real-valued Brownian motion and $\circ$ denotes the
Stratonovich product. The second equation is the corresponding It\^{o}
equation. In addition, $\gamma_c, \gamma_s$ are two constants
determined by $\alpha$ and given by
\[
\gamma_c = \int_0^{\infty} \cos(2bt)\mathbb{E}(\alpha
(0)\alpha(t)
)\,dt  \quad\mbox{and}\quad  \gamma_s = \int_0^{\infty} \sin
(2bt)\mathbb
{E}(\alpha(0)\alpha(t) )\,dt.
\]
Then $\nu(t) \in\mathbb{S}^3$ a.s., the unit sphere in $\mathbb{C}^2
\sim\mathbb{R}^4$. We denote by $\Lambda$ the unique invariant
probability measure of $\nu$ (see Section~\ref{A1}) and by $\mathbb
{E}_{\Lambda}(\cdot)$ the expectation with respect to
$\Lambda$.
Thus, replacing $u$ by $\nu$ in (\ref{system3u}), we obtain a new
equation describing the evolution of the electric field envelope $X=
(X_1, X_2)^t$:
%
%e1.9 #&#
\begin{eqnarray}\label{system3}
&&
i\,\frac{\partial X}{\partial t}+\frac{d_0}{2}\,\frac{\partial^2
X}{\partial x^2}+ \frac{8}{9}|X|^2X\nonumber\\[-8pt]\\[-8pt]
&&\qquad =-ib'\sigb( \nu (t)) \,\frac{\partial X}{\partial
x}-\frac{1}{6}\bigl(N_{\nu}(X )- \mathbb{E}_{\Lambda}(N_{\nu}(X))\bigr);\nonumber
\end{eqnarray}
indeed, the process $m=(m_1,m_2,m_3)$ is now defined as $m=(g_1(\nu
),g_2(\nu),g_3(\nu))$ and it can be proved (see Section~\ref{A1}) that
\begin{eqnarray*}
\mathbb{E}_{\Lambda}(N_{1, \nu}(X ))&=&\tfrac
{2}{3}(2|X_2|^2-|X_1|^2
)X_1,\\
\mathbb{E}_{\Lambda}(N_{2, \nu}(X ))&=&\tfrac
{2}{3}(2|X_1|^2-|X_2|^2)X_2.
\end{eqnarray*}
We set
%
%e1.10 #&#
\begin{equation}\label{NLterms}
F_{\nu(t)}{( X(t)) } = \tfrac{8}{9}|X|^2X -\tfrac{1}{6}
\bigl(N_{\nu}
(X )- \mathbb{E}_{\Lambda}(N_{\nu}(X
))\bigr).
\end{equation}
Equation (\ref{system3}) is of great interest for the study of
dispersion because the main effects leading to signal distortions (Kerr
effect, chromatic dispersion, PMD) can be easily identified: on the
left-hand side, the first term describes the evolution of the pulse
along the fiber. The second one corresponds to the chromatic dispersion
and the last term to the Kerr effect averaged on the Poincar\'{e}
sphere. On the right-hand side of the equation, the first term
describes the linear PMD effect and the second term describes nonlinear PMD.

The Manakov PMD equation (\ref{system3}) is written in dimensionless
form. According to the length scales we consider, we set $\eps
{X}(t,x)=\frac{1}{\epsilon}X(\frac{t}{\epsilon^2},\frac
{x}{\epsilon
})$ and $ \eps{\nu}(t) = \nu(\frac{t}{\epsilon
^2})$,
where $\nu$ is\vadjust{\goodbreak} the solution of (\ref{nu1}); then the electric field
$\eps{X}$ has the following evolution:
%
%e1.11 #&#
\begin{equation}\label{manakovPMD}\qquad
i\,\frac{\partial\eps{X}(t)}{\partial t}+\frac{ib'}{\epsilon
}\sigb( \eps{\nu
}(t)) \,\frac{\partial\eps{X}(t)}{\partial x}+\frac
{d_0}{2}\,\frac{\partial^{2} \eps{X}(t)}{\partial x^{2}}+F_{\eps {\nu
}(t)}{( \eps{X}(t)) }=0,
\end{equation}
where the term $F_{\eps{\nu}(t)}{( \eps{X}(t)) }$ is given by (\ref
{NLterms}).

In various physical situations, the long time behavior of a phenomenon
subject to random perturbations requires to take care of the different
characteristic length scales of the problem. In this context
Papanicolaou, Stroock and Varadhan~\cite{papanicolaou} and
Blankenship and Papanicolaou~\cite{papa} introduced the approximation
diffusion theory for random ordinary differential equations. This
method has been used to study wave propagation in random media~\cite
{garnier} and, in particular, in randomly birefringent fibers \cite
{marty,2}, but only few results exist on limit theorems for random
PDEs. In the latter, the authors studied the evolution, in an optical
fiber, of the linear field envelope $\eps{X}$ given by\looseness=-1
\[
i\,\frac{\partial\eps{X}(t)}{\partial t}+\frac{ib'}{\epsilon
}\sigb( \eps{\nu
}(t))\,\frac{\partial\eps{X}(t)}{\partial x}+\frac
{d_0}{2}\,\frac{\partial^{2} \eps{X}(t)}{\partial x^{2}}=0
\]\looseness=0
and proved that the asymptotic dynamics, when $\epsilon$ goes to zero,
is given by
\[
i\,dX(t)+\biggl( \frac{d_0}{2}\,\frac{\partial^{2} X(t)}{\partial
x^{2}}\biggr)\,dt+i\sqrt{\gamma
}\sum
_{k=1}^3 \sigma_k\,\frac{\partial X(t)}{\partial x}\circ dW_k(t)=0,
\]
where $W=(W_1,W_2,W_3)$ is a $3$d Brownian motion, and $\gamma
=(b')^2/6\gamma_c$.
Note that the linear PMD effect reduces to one single parameter $\gamma
$ in front
of the three Brownian motions. Generalizing the perturbed test function method,
we will prove that the asymptotic dynamic of (\ref{manakovPMD}) is
given by the
stochastic nonlinear evolution:
%
%e1.12 #&#
\begin{eqnarray}\label{manakovlimite}
&&i\,dX(t)+\biggl( \frac{d_0}{2}\,\frac{\partial^{2} X(t)}{\partial x^{2}}
+F_{}{( X(t)) }
\biggr)\,dt\nonumber\\
&&\quad{}+i\sqrt
{\gamma}\sum_{k=1}^3 \sigma_k\,\frac{\partial X(t)}{\partial
x}\circ dW_k(t)\\
&&\qquad=0,\nonumber
\end{eqnarray}
where the nonlinear function $F$ reduces to $F_{}{( X(t)) }= \frac
{8}{9}|X(t)|^2X(t)$ that is simply the expectation, with respect to
the invariant measure $\Lambda$, of\break $F_{\eps{\nu}(t)}{( \eps{X}(t)) }$.
%Note that in the limit the nonlinear PMD effect averages to zero and
%it only remains the contribution of the linear PMD effect on the pulse
%evolution.
We will also make use of the following equivalent It\^{o} formulation:
%
%e1.13 #&#
\begin{eqnarray}\label{manakovlimiteito}
&&
i\,dX(t)+\biggl(\biggl(\frac{d_0}{2}-\frac{3i\gamma}{2}\biggr)\,\frac
{\partial^{2} X(t)}{\partial x^{2}} +F_{}{( X) }(t)\biggr)\,dt\nonumber\\
&&\quad{}+i\sqrt
{\gamma}\sum_{k=1}^3
\sigma
_k\,\frac{\partial X(t)}{\partial x} \,dW_k(t)\\
&&\qquad=0.\nonumber
\end{eqnarray}

Note that a different regime concerned with long propagation distances
and corresponding to $l \gg l_{nl}\sim l_{d}$ is of physical interest;
however, this regime would lead to another asymptotic analysis which is
beyond the scope of this paper.

This paper is organized as follows: in Section~\ref{section2} we give
notation that will be used along the paper and state the main results.
Section~\ref{section3} is devoted to the proof of well-posedness for
the Manakov PMD equation. In Section~\ref{section4} we study the local
well-posedness of the limiting equation (\ref{manakovlimite}). Finally,
in Section~\ref{section5} we prove the convergence in law of $\eps{X}$
to $X$ as $\epsilon$ goes to zero. This paper ends with Section \ref
{A1} where we recall some results obtained in~\cite{marty,2} about
the driving process $\nu$, and Section~\ref{A2} where proofs of
technical results used in Section~\ref{section5} are gathered. %In
%process given by the martingale problem.

%s1.2 #&#
\subsection{Notation and main results}\label{section2}

Before stating the main results of this article, let us give some
definitions and notation.

For all $p\geq1$, we define $\mathbb{L}^p(\mathbb{R})=
(L^p(\mathbb{R}; \mathbb{C}) )^2$ the Lebesgue spaces of
functions with values in $\mathbb{C}^2$. Identifying $\mathbb{C}$ with
$\mathbb{R}^2$, we define a scalar product on $\mathbb{L}^2(
\mathbb{R}) $ by
\[
(u,v)_{\mathbb{L}^2}=\sum_{i=1}^2\Re\biggl\{\int
_{\mathbb{R}}u_i\overline{v}_i\,dx\biggr\}.
\]
We denote by $\mathbb{W}^{m,p}, m\in\mathbb{N}^*, p\in\mathbb{N}^*$
the space of functions in $\mathbb{L}^p$ such that their $m$ first
derivatives are in $\mathbb{L}^p$. If $p=2$, then we denote $\mathbb
{H}^m( \mathbb{R}) = \mathbb{W}^{m,2}( \mathbb
{R}
)$, $m \in\mathbb{N}$. We will also use $\mathbb{H}^{-m}$ the
topological dual space of $\mathbb{H}^m$ and denote $\langle\cdot,\cdot
\rangle$ the
paring between $\mathbb{H}^m$ and $\mathbb{H}^{-m}$.
%Let $\Omega$ be an open subset of $\mathbb{R}^d$. We define the
%fractional Sobolev space $\mathbb{W}^{\alpha,p}(\Omega)$ as an
%intermediate family space between $\mathbb{L}^p(\Omega)$ and $
% \mathbb{W}^{\alpha,p}(\Omega)=\lbrace u \in\mathbb{L}^p(\Omega)
%; \frac{\abs{u(x)-u(y)}}{\abs{x-y}^{\alpha+d/p}} \in\mathbb{L}^p(
The Fourier transform of a tempered distribution $v \in\mathcal
{S}'(\mathbb{R})$
is either denoted by $\widehat{v}$ or $\mathcal{F}v$. If $s \in
\mathbb
{R}$, then
$\mathbb{H}^s$ is the fractional Sobolev space of tempered distributions
$v \in\mathcal{S}'(\mathbb{R})$ such that
$(1+|\xi|^2)^{s/2}\widehat{v} \in\mathbb{L}^2$.
Let $(E, \|\cdot\|_{E} )$ and $(F, \|\cdot\|_{F}
)$
be two
Banach spaces.
We denote by $\mathcal{L}( E,F)$ the space of linear continuous
functions from $E$ into $F$, endowed with its natural norm. If $I$ is
an interval
of $\mathbb{R}$ and $1\leq p \leq+\infty$, then $L^p
(I;E) $
is the space of strongly Lebesgue measurable functions $f$ from $I$
into $E$ such
that $t \mapsto\|f(t)\|_{E}$ is in $L^p(I)$. The space
$L^p(\Omega, E) $ is defined similarly where
$(\Omega, \mathcal{F}, \mathbb{P})$ is a probability space.
We denote by $L_w^{p}(I, E)$ the space $L^{p}(I,
E)$
endowed with the weak (or weak star) topology. For a real number $0<
\alpha<1$
and $p \geq1$, we denote by $W^{\alpha,p}([0,T],E
)$ the
fractional Sobolev space of functions $u$ in $L^p(0,T; E)$
satisfying
\[
\int_0^T\int_0^T\frac{\|u(t)-u(s)\|_{E}^p}{|t-s|^{\alpha p
+1}}\,ds\,dt<+\infty.
\]
The space $C^{\beta}([0,T]; E)$ is the space of
H\"{o}lder continuous functions of order \mbox{$\beta>0$} with values in $E$
and we denote by $\mathcal{M}(E)$ the set of probability measures on
$E$, endowed with the topology of the weak convergence $\sigma
(\mathcal{M}(E), C_b(E))$.

We will use the space
\[
\mathcal{K}= \bigl(C([0,T], \mathbb{H}_{\mathrm{loc}}^1
) \cap C_w
([0,T], \mathbb{H}^1)\cap L_w^{\infty}(0,T;
\mathbb{H}^2)\bigr) \times C([0,T], \mathbb{R}),
\]
%
% \[
% \mathcal{K}= (C([0,T], \mathbb{H}_{\mbox{\mathrm{loc}}}^1
%)
%([0,T], \mathbb{R}),
% \]
where $C_w([0,T], \mathbb{H}^m), m \in
\mathbb
{Z}$ is the space of functions $f$ in $L^{\infty}(0,T; \mathbb
{H}^m)$, weakly continuous from $[0, T]$ into $\mathbb{H}^m$. As
the solution of our limit equation will not necessary be global in
time, we need to introduce a space of exploding paths, as in \cite
{azencott}, by adding a point $\Delta$, which acts as a cemetery
point, at infinity in $\mathbb{H}^1$; then
\begin{eqnarray*}
&&
\mathcal{E}(\mathbb{H}^1) = \bigl\{ f \in C
([0,T], \mathbb{H}^{1} \cup\{\Delta\}),\\
&&\hspace*{46.5pt} f(t_0)= \Delta \mbox{ for } t_0 \in[0,T] \Rightarrow f(t) =\Delta
\mbox{ for } t \in[t_0,T] \bigr\}.
\end{eqnarray*}
We define a topology on $\mathbb{H}^1 \cup\{\Delta
\}$ such that the open sets of $\mathbb{H}^1 \cup\{
\Delta\}$ are the open sets of $\mathbb{H}^1$ and the
complementary in $\mathbb{H}^1 \cup\{\Delta
\}$
of the closed bounded sets in $\mathbb{H}^1$. %For any $T>0$, the space
%$C([0,T], \mathbb{H}^{1} \cup\Delta)$ is the
%space of continuous functions on $[0,T]$ with values in $
For any $f \in C([0,T], \mathbb{H}^{1} \cup
\{\Delta\})$ we denote the blowing-up time
$\tau
(f)$ by
\[
\tau(f) = \inf\{ t \in[0,T], f(t) = \Delta\}
\]
with the convention $\tau(f) = + \infty$ if $f(t) \neq\Delta$ for all
$t \in[0,T]$. We endow the space $\mathcal{E}(\mathbb
{H}^1
)$ with the topology induced by the uniform convergence in $\mathbb
{H}^1$ on every compact set of $[0, \tau(f))$.
%We say that a sequence $f_n \in\mathcal{E}(\mathbb{H}^1)$
%converges to $f \in\mathcal{E}(\mathbb{H}^1)$ if $f_n$
%uniformly converges to $f$ in $\mathbb{H}^1$ on every compact set of
%$[0, \tau(f))$.
% We set $\mathcal{Y}_m=C([0,T], \mathbb{H}^{m+1}
%)\cap C^{\beta}([0,T], \mathbb{H}^m)$, $
%([0,T], \mathbb{R})$. Moreover we denote respectively by $
%$\mathcal{Z}_m$ with radius $M$.

% To get tightness of the probability measure on a suitable space we
%need some classical compactness results that we state below:
% \begin{lemma}\label{compactimbedding}
% Let $\alpha> 0$ and $T>0$; then the embeddings
% \[
% C([0,T], \mathbb{H}^{2})\cap C^{\alpha}(
%[0,T], \mathbb{H}^{-1})\hookrightarrow C_w(
%[0,T], \mathbb{H}^1)
% \]
% and
% \[
% C([0,T], \mathbb{H}^{2})\cap C^{\alpha}(
%[0,T], \mathbb{H}^{-1}) \hookrightarrow C([0,T],
% \]
% are compact.
% \end{lemma}
% % \begin{pf}
% % We denote by $B_m(M)$ the closed ball of $\mathcal{Y}_m$ of radius
%M. To prove that $B_m(M)$ and $B_0(M)$ are compact in respectively $
%of $B_{\mathbb{H}^m}$ for the weak topology $\sigma(\mathbb{H}^m,
% % \end{proof}
% It follows from Lemma~\ref{compactimbedding}, Arzela Ascoli and
%Banach Alaoglu Theorems that
% \begin{cor}\label{compactimbeddingfork}
% For $\alpha, \gamma>0$ and $T>0$; the following embedding is compact
% \[
% ( C([0,T], \mathbb{H}^{2})\cap C^{\alpha}
%([0,T], \mathbb{H}^{-1}) ) \times C^{
% \]
% \end{cor}

Let $(\mathcal{A}, \mathcal{G},\mathbb{Q})$ be a probability
space endowed with the complete filtration $(\mathcal{G}_t
)_{t \geq0}$ generated by a two-dimensional Brownian motion
$W=(W_1,W_2)$ which is driving the diffusion process $\nu$
given by (\ref{nu1}). We first state an existence and uniqueness result
for (\ref{manakovPMD}).
%
%th1.1 #&#
\begin{theorem}\label{theorem1}
Let $\epsilon>0$ and suppose that $\eps{X}(0)=v \in\mathbb
{L}^2(\mathbb{R})$, then there exists a unique global solution $\eps
{X}$ to (\ref{manakovPMD}) such that, $\mathbb{Q}$-almost surely,
\[
\eps{X} \in C(\mathbb{R}_+, \mathbb{L}^2) \cap C^1
(\mathbb{R}_+, \mathbb{H}^{-2}) \cap\mathbb{L}^8_{\mathrm{loc}}
(\mathbb{R}_+, \mathbb{L}^4).
\]
Moreover, equation (\ref{manakovPMD}) preserves the $\mathbb{L}^2$
norm, that is, for all $t \in\mathbb{R}_+$
\[
\|\eps{X}(t)\|_{\mathbb{L}^2}=\|v\|_{\mathbb{L}^2}.
\]
If, in addition, $\eps{X}(0) =v \in\mathbb{H}^1$ $(\mbox
{resp., }
\mathbb{H}^2, \mbox{resp., } \mathbb{H}^3)$, then the
corresponding solution is in $C(\mathbb{R}_+, \mathbb
{H}^1)$
$[\mbox{resp., } C(\mathbb{R}_+, \mathbb{H}^2),
\mbox
{resp., } C(\mathbb{R}_+, \mathbb{H}^3)]$.
%whose $\mathbb{H}^1$ norm satisfy:
%)} &\leq& \norm{\pds{}{x}{\Psi(0)}}{\mathbb{L}^2}+CT\norm{
%CT^{1/2}(\norm{\Psi(0)}{\mathbb{L}^2}+ T^{1/2}\norm{\eps{
\end{theorem}

Let $(\Omega, \mathcal{F}, \mathbb{P})$ be a probability space on
which is defined a three-dimensional real-valued Brownian motion $W=
(W_1, W_2, W_3)$. We denote by $(\mathcal{F}_t)_{t
\in
\mathbb{R}_+}$ the complete filtration generated by $W$. The next
theorem gives existence and uniqueness of the local solution for (\ref
{manakovlimite})
%
%th1.2 #&#
\begin{theorem}\label{theorem2}
Let $X_0=v \in\mathbb{H}^1(\mathbb{R})$, then there exists a maximal
stopping time $\tau^*(v,\omega)$ and a unique strong solution $X$ (in
the probabilistic sense) to (\ref{manakovlimite}), such that $X \in
C([0, \tau^*), \mathbb{H}^1(\mathbb{R}) )$
$\mathbb{P}$-a.s. Furthermore, the $\mathbb{L}^2$ norm is almost surely\vadjust{\goodbreak}
preserved, that is, $\forall t \in[0, \tau^*), \|X(t)\|_{\mathbb
{L}^2}=\|v\|_{\mathbb{L}^2}$ and the following alternative holds for
the maximal existence time of the solution:
\[
\tau^*(v,\omega) = +\infty \quad\mbox{or}\quad  \limsup _{t
\nearrow
\tau^*(v,\omega)} \|X(t)\|_{\mathbb{H}^1}=+\infty.
\]
Moreover, if $v \in\mathbb{H}^2$, then $X \in C([0, \tau^*),
\mathbb{H}^2(\mathbb{R}) )$ and $\tau^*$ satisfies
%
%e1.14 #&#
\begin{equation}\label{limitnormH1}
\tau^*(v,\omega) = +\infty \quad\mbox{or}\quad  \lim _{t \nearrow
\tau
^*(v,\omega)} \|X(t)\|_{\mathbb{H}^1}=+\infty.
\end{equation}
\end{theorem}

Note that we do not obtain global existence for (\ref{manakovlimite}),
due to the lack of control of the evolution of the $\mathbb{H}^1$ norm
(see Remark~\ref{onlocalexistence}).

Using these existence theorems, we are able to prove a diffusion
approximation result for the nonlinear system of PDEs (\ref{manakovPMD}).
%Theorem~\ref{theorem2} gives only local existence of a solution for
%allow us to consider convergence in law of the processes in the case
%where blow up may occur for the limit equation. We add a point $
%open sets of $\mathbb{H}^1 \cup\{\Delta\}$ are
%the open sets of $\mathbb{H}^1$ and the complementary in $\mathbb{H}^1
%if we endow it with the distance $d_s$ induced by the ``stereographic
%projection''~\cite{bF}, which is defined as follows: let $\phi$ be
%the mapping from $\mathbb{H}^1 \cup\{\Delta\}$
%into $\mathbb{H}^1 \times[-1,1]$ defined by
%if $v \in\mathbb{H}^1$ and $\phi(\Delta)=(0,1)$. Then $\phi$ is a
%homeomorphism from $\mathbb{H}^1 \cup\{\Delta\}
%$ into its image. It suffices then to define the distance
% d_s(\overline{u},\overline{v})=\norm{\phi(\overline{u})-\phi(
%for $\overline{u}$ and $\overline{v}$ in $\mathbb{H}^1 \cup
%Theorem~\ref{theorem2} may be extended to a solution a.s in $C
%([0,T], \mathbb{H}^1 \cup\{\Delta\}
%)$ for any $T>0$, by setting $X(t)= \Delta$ for $t \in[
%
%th1.3 #&#
\begin{theorem}\label{theorem3}
Let $\eps{X}(0)= X_0 = v$ be in $\mathbb{H}^3( \mathbb
{R})$,
then the solution
$\eps{X}$ of (\ref{manakovPMD}) given by Theorem~\ref{theorem1}
converges in law to the
solution $X$ of (\ref{manakovlimite}) in
$\mathcal{E}(\mathbb{H}^1 )$, that is,
for all functions $f$ in $C_b(\mathcal{E}(\mathbb{H}^1
))$,
\[
\lim _{\epsilon\to0} \mathcal{L}( \eps{X} )(f)
=\mathcal{L}( X )(f).
\]
\end{theorem}

Note that we consider here the Manakov PMD equation (\ref{manakovPMD}),
but the method may be carried out to other nonlinear Schr\"{o}dinger
equations. Let us first emphasize the key points that allow us to prove
Theorem~\ref{theorem3}.

The first point is that the noise term is a linear function of the
unknown $\eps{X}$. This particular structure leads to a stochastic
partial differential equation for the limiting equation. The second
point is the fact that the Pauli matrices are Hermitian. This is
important to obtain the conservation of the $\mathbb{L}^2$ norm for
both equations. Finally, we use that the driving process $\nu$ is a
homogeneous Markov ergodic process defined on a compact state space
such that $\mathbb{E}_{\Lambda}( \sigb(y))=0$. The
hypothesis on the driving noise may be weakened as in the case of a
random ordinary differential equation assuming good mixing properties
(e.g., exponential decay of the covariance function). The
boundedness of $\sigb(\eps{\nu}(t))$ seems to be necessary.
It is used to prove uniform bounds in Lemma~\ref{lemmatension} for
tightness. %Other function of the noise may be considered. But in this
%case it may happen that only local existence results can be proved for
%performed. The proof of the convergence result can be slightly
%modified taking care of the maximal stopping time existence $\eps{
On the other hand, the lack of Strichartz estimates for the limiting
equation (\ref{manakovlimite}) is a negative aspect. Thus, we use that
$F(v)$ is locally Lipschitz in $\mathbb{H}^1( \mathbb{R})$
to prove existence and uniqueness of a local solution to
(\ref{manakovlimite}). But if $\sigb(\eps{\nu}(t))$
were a
one-dimensional process, larger dimension and larger power in the
nonlinear term could be considered.

Other types of nonlinear Schr\"{o}dinger equations may be considered
replacing, for example, $i\,\frac{\partial\eps{X}}{\partial x}$
by $\eps{X}$ and assuming
that the matrices $\sigma_k$ are real valued and symmetric. This latter
equation is simpler to handle using Strichartz estimates for the
fundamental solution and because $\sigb( \eps{\nu}(t))
\eps
{X}(t)$ can be treated as a perturbation as far as we are concerned with
existence of solutions.

%s2 #&#
\section{\texorpdfstring{The Manakov PMD equation: Proof of Theorem \protect\ref{theorem1}}
{The Manakov PMD equation: Proof of Theorem 1.1}}\label{section3}

The point here is that no Strichartz estimates are available for (\ref
{manakovPMD}) because of the lack of commutativity of the matrix $\sigb$
at a different time: $\sigb(\nu(t) ) \sigb(\nu
(s)
)\neq\sigb(\nu(s) )\sigb(\nu(t) )$.
Consequently, only local existence and uniqueness for initial data in
$\mathbb{H}^1$ can be easily proved directly on (\ref
{manakovPMD}). The idea of the proof is then to find a unitary
transformation such that Strichartz estimates are available for the
transformed equation. This change of unknown is given in the next result.
%
%le2.1 #&#
\begin{lemma}
Let us denote for $t \in\mathbb{R}_+$
\[
\eps{Z}(t)=\pmatrix{
\nu_{1,\epsilon}(t) & \overline{\nu}_{2,\epsilon}(t)\cr
-\nu_{2,\epsilon}(t) & \overline{\nu}_{1,\epsilon}(t)}
,
\]
where $\eps{\nu}= \nu(t / \epsilon^2) $, $\nu$ given by
(\ref{nu1}). Assuming that $\eps{X} \in C([0,T],
\mathbb
{L}^2)$, we set $\eps{\Psi}(t)=\eps{Z}(t)\eps{X}(t)$; then the
evolution of the electric field $\eps{\Psi}$ is given by the stochastic
It\^{o} equation
%
%e2.1 #&#
\begin{eqnarray}\label{system2pert}\qquad
&&i\,d\eps{\Psi}(t) + \biggl\{\frac{ib'}{\epsilon} \sigma_3\,\frac
{\partial\eps{\Psi}}{\partial x}+\frac{d_0}{2}\,\frac{\partial^2
\eps
{\Psi}}{\partial x^2}+\frac{5}{6}|\eps{\Psi}|^2\eps
{\Psi}
+\frac{1}{6}(\eps{\Psi}^*\sigma_3\eps{\Psi})\sigma
_3\eps{\Psi
}\biggr\}\, dt \nonumber\\[-8pt]\\[-8pt]
&&\qquad{}+\frac{\gamma_s}{\epsilon^2}\sigma_3\eps{\Psi}\,dt +\frac{i\gamma
_c}{\epsilon^2}\eps{\Psi}\,dt- \frac{\sqrt{\gamma_c}}{\epsilon
}\bigl(
\sigma_1\eps{\Psi}\,d\widetilde{W}_1(t)+\sigma_2\eps{\Psi
}\,d\widetilde
{W}_2(t)\bigr) = 0,\nonumber
\end{eqnarray}
where $\widetilde{W}_j(t) = \epsilon W_j(t/\epsilon^2 ),
j=1,2$, and with initial conditions
\[
\eps{\Psi}(0)=\pmatrix{
\nu_{1,\epsilon}(0)v_1+\overline{\nu}_{2,\epsilon}(0)v_2\cr
-\nu_{2,\epsilon}(0)v_1+
\overline{\nu}_{1,\epsilon}(0)v_2}
= \psi_0.
\]
%
% Moreover the matrices $\eps{\widetilde{\Sigma}}$ are given by
% \[
% \compoeps{\widetilde{\Sigma}}{0}(t)=\pmatrix{
% \gamma_s(\abs{\compoeps{\nu}{1}}^2- \abs{\compoeps{\nu}{2}}^2
%) & -2\gamma_s\compoeps{\nu}{1}\compoeps{\overline{\nu}}{2} \\
% -2\gamma_s\compoeps{\overline{\nu}}{1}\compoeps{\nu}{2} & \gamma_s
%(\abs{\compoeps{\nu}{2}}^2- \abs{\compoeps{\nu}{1}}^2)
% \end{pmatrix},
% \]
% \[
% \compoeps{\widetilde{\Sigma}}{1}(t)=\sqrt{\gamma_c}\pmatrix{
% 2\Re(\compoeps{\nu}{1}\compoeps{\nu}{2}) & (
% \compoeps{\overline{\nu}}{1}^2-\compoeps{\nu}{2}^2 & -2\Re(
% \end{pmatrix},\]\[\compoeps{\widetilde{\Sigma}}{2}(t)=\sqrt{\gamma_c}
% 2\Im(\compoeps{\overline{\nu}}{2}\compoeps{\overline{\nu}}{1}
%) & i(\compoeps{\nu}{1}^2+\compoeps{\overline{\nu}}{2}^2
%) \\
% -i( \compoeps{\overline{\nu}}{1}^2+\compoeps{ \nu}{2}^2) &
%-2\Im(\compoeps{\overline{\nu}}{2}\compoeps{\overline{\nu}}{1}
%)
% \end{pmatrix}.
% \]
\end{lemma}
\begin{pf}
% By integration by parts, we get
% \[
% d\eps{\Psi}(t)=d\eps{Z}(t) \eps{X}(t) + \eps{Z}(t) \pds{}{t}{
% \]
% The first term of the above expression is given by
% \[
% d\eps{Z}(t) \eps{X}(t)&=&d\eps{Z}(t)\eps{Z}^{-1}(t) \eps{\Psi}(t)\\&=&
%)+\compoeps{\Psi}{2} (\compoeps{\nu}{1}\,d\compoeps{\overline{
% -\compoeps{\Psi}{1}(\compoeps{\overline{\nu}}{1}\,d\compoeps{
% \end{pmatrix}.
% \]
Using the equation satisfied by $\eps{\nu}$ and because $|\nu
_{1,\epsilon}(t)|^2 +\break |\nu_{2,\epsilon}(t)|^2=1$ for any $t\geq0$,
we obtain
\begin{eqnarray*}
i\,d\eps{Z}(t) \eps{Z}^{-1}\eps{\Psi}(t)&=& -\frac{\gamma_s}{\epsilon
^2}\sigma_3\eps{\Psi}\,dt-\frac{i\gamma_c}{\epsilon^2}\eps{\Psi
}\,dt + \frac
{\sqrt{\gamma_c}}{\epsilon}\sigma_1\eps{\Psi}\,d\widetilde{W}_1(t)\\
&&{}+ \frac
{\sqrt{\gamma_c}}{\epsilon}\sigma_2\eps{\Psi}\,d\widetilde{W}_2(t).
\end{eqnarray*}
The nonlinear part of (\ref{system2pert}) is obtained as in
the derivation of (\ref{system3u}).
\end{pf}

% \[
% \eps{\sigb}(t)=\sigma_1\compoeps{m}{1}(t)+\sigma_2\compoeps{m}{2}(t)+
% \]
%Note that the above transformation splits the time dependent unbounded
%operator $\frac{ib'}{\epsilon}\sigb(\nu(t) )\pds{}{x}{}$
%into a time dependent bounded operator and a time independent
%unbounded operator. This splitting allows us to obtain global
%wellposedness for Equation \eqref{system2pert}.
We first investigate the behavior of the linear equation
%
%e2.2 #&#
\begin{equation}\label{linearsystem2pert}
i\,\frac{\partial\eps{\Psi}}{\partial t} + \frac{1}{\epsilon
}ib' \sigma_3\,\frac
{\partial
\eps{\Psi}}{\partial x}+\frac{d_0}{2}\,\frac{\partial^2 \eps{\Psi
}}{\partial x^2} =0
\end{equation}
with initial condition $\eps{\Psi}(0)=\psi_0\in\mathbb{L}^2$.
%
%pr2.1 #&#
\begin{prop}\label{semigroup}
The unbounded matrix operator $\eps{H}= \frac{id_0}{2}I_2\,\frac
{\partial^{2} }{\partial x^{2}} -
\frac{b'}{\epsilon}\sigma_3\,\frac{\partial}{\partial x}$
defined on $\mathscr
{D}(
\eps{H}) = \mathbb{H}^2$ is the infinitesimal generator of a
unique strongly continuous unitary group $\eps{U}(t)$ on $\mathbb
{L}^2$. Moreover,\vadjust{\goodbreak} $\eps{U}(t)$ may be expressed as a convolution
kernel, that is, for $\psi_0 \in\mathcal{S}(\mathbb{R})$
\begin{eqnarray*}
\eps{U}(t) \psi_0 &=& \eps{A}(t) \star\psi_0 \\
&=& \frac{1}{\sqrt{2\pi i
d_0t}}\pmatrix{
\displaystyle \exp\biggl\{\frac{i}{2}\frac{(x-b't/\epsilon
)^2}{d_0t}
\biggr\} & 0
\vspace*{2pt}\cr
0 &
\displaystyle \exp\biggl\{\frac{i}{2}\frac{(x+b't/\epsilon
)^2}{d_0t}\biggr\}}
\star\psi_0.
\end{eqnarray*}
\end{prop}
\begin{pf}%[Proof of Proposition~\ref{semigroup}]
Assuming $\psi_0 \in\mathcal{S}(\mathbb{R})$ and
taking the
Fourier transform, in the space variable, of (\ref{linearsystem2pert}),
we obtain readily
\[
\frac{\partial\eps{\widehat{\Psi}}}{\partial t}= -\frac
{1}{\epsilon}ib' \sigma
_3\xi\eps
{\widehat{\Psi}}-i\,\frac{d_0\xi^2}{2}\eps{\widehat{\Psi}}.
\]
Since $\sigma_3$ does not depend on time, we obtain
\[
\eps{\widehat{\Psi}}(t)=\eps{R}(t)\widehat{\psi}_0= \pmatrix{
\displaystyle \exp\biggl\{-\frac{id_0}{2}\xi^2t - i\frac{b'}{\epsilon}\xi
t\biggr\} & 0
\vspace*{2pt}\cr
0 & \displaystyle \exp\biggl\{-\frac{id_0}{2}\xi^2t + i\frac{b'}{\epsilon}\xi
t\biggr\}}
\widehat{\psi}_0.
\]
The statement of Proposition~\ref{semigroup} follows then in a
classical way, setting $\eps{A}(t) = \mathcal{F}^{-1}(\eps
{R}(t))$.
\end{pf}

The explicit formulation of the kernel given in Proposition \ref
{semigroup} allows immediately to get the following dispersive
estimates: if $p \geq2$, $t \neq0$, then $\eps{U} \in\mathcal
{L}( \mathbb{L}^{p'},\mathbb{L}^{p} )$ where $p'$ is such
that $\frac{1}{p}+\frac{1}{p'}=1$ and for all $\psi_0 \in\mathbb{L}^{p'}$,
%
%e2.3 #&#
\begin{equation}\label{decay}
\|\eps{U}(t)\psi_0\|_{\mathbb{L}^{p}} \leq(2\pi
|d_0||t|)^{-1/2 +1/p}\|\psi_0\|_{\mathbb{L}^{p'}}.
\end{equation}
Using then classical arguments (see~\cite{caz,ginibre}), one may
prove Strichartz inequalities for $\eps{U}(t)$.
%
%pr2.2 #&#
\begin{prop}\label{strichartz}
The following properties hold:
\begin{longlist}[(2)]
\item[(1)] For every $\psi_0 \in\mathbb{L}^2(\mathbb{R})$,
$\eps
{U}(\cdot)\psi_0 \in L^8(\mathbb{R}; \mathbb{L}^4) \cap
C
(\mathbb{R}; \mathbb{L}^2)$. Furthermore, there exists a constant
$C$ such that
\[
\|\eps{U}(\cdot)\psi_0\|_{L^8(\mathbb{R}; \mathbb{L}^4)}
\leq C \|\psi_0\|_{\mathbb{L}^2}  \qquad\mbox{for every }
\psi_0
\in\mathbb{L}^2.
\]
\item[(2)] Let $I$ be an interval of $\mathbb{R}$ and $t_0 \in I$. Let $f
\in L^{8/7}(I, \mathbb{L}^{4/3})$, then the function
\[
t \mapsto\int^{t}_{t_0}\eps{U}(t-s)f(s)\,ds
\]
belongs to $L^{8}(I, \mathbb{L}^{4}) \cap C(I,
\mathbb
{L}^2)$. Furthermore, there exists a constant $C$ independent of
$I$ such that for every $f \in L^{8/7}(I, \mathbb{L}^{4/3})$
\[
\biggl\|\int^{\cdot}_{t_0}\eps{U}(\cdot-s)f(s)\,ds\biggr\|_{L^{8}(I, \mathbb
{L}^{4})\cap L^{\infty}(I, \mathbb{L}^{2})}
\leq
C \|f\|_{L^{8/7}(I, \mathbb{L}^{4/3})}.\vadjust{\goodbreak}
\]
\end{longlist}
\end{prop}

We now turn to the study of the nonlinear problem. We will use, as is
classical, a cutoff argument on the nonlinear term which is not
Lipschitz. The cutoff we consider here is of the same form as the one
considered in~\cite{bouard}. We first prove an existence and
uniqueness result for this truncated equation, then deduce from this
result the existence of a unique solution for (\ref
{system2pert}). We denote:
\[
f(\eps{\Psi})=\tfrac{5}{6}|\eps{\Psi}
|^2\eps{\Psi}
+\tfrac{1}{6}(\eps{\Psi}^*\sigma_3\eps{\Psi})\sigma
_3\eps{\Psi}.
\]
Let $\Theta\in C_c^{\infty}(\mathbb{R})$ with
supp$\Theta
\subset[-2;2]$ such that $\Theta(x)= 1$ for $|x|\leq1$ and
$0 \leq\Theta(x) \leq1$ for $x \in\mathbb{R}$. Let $R>0$
and $\Theta_R(x)=\Theta(x/R )$. We then consider the
following equation:
%
%e2.4 #&#
\begin{eqnarray}\label{system2perttronquee}
\reps{\Psi}(t)&=& \eps{U}(t)\psi_0 + \frac{i\gamma_s}{\epsilon
^2}\int
_0^t\eps{U}(t-s)\sigma_3\reps{\Psi}(s)\,ds \nonumber\\
&&{}- \frac{\gamma
_c}{\epsilon
^2}\int_0^t\eps{U}(t-s)\reps{\Psi}(s)\,ds\nonumber\\
&&{} + i\int_0^t\eps{U}(t-s) \Theta_R\bigl(\|\reps{\Psi}\|
_{L^8(0,s;\mathbb{L}^4) } \bigr) f(\reps{\Psi
}(s))
\,ds \\
&&{} - \frac{i\sqrt{\gamma_c}}{\epsilon}\int_0^t\eps{U}(t-s)\sigma
_1 \reps
{\Psi}(s)\,d\widetilde{W}_1(s)\nonumber\\
&&{}- \frac{i\sqrt{\gamma_c}}{\epsilon
}\int
_0^t\eps{U}(t-s)\sigma_2 \reps{\Psi}(s)\,d\widetilde
{W}_2(s),\nonumber
\end{eqnarray}
which is the mild form of the It\^{o} equation,
%
%e2.5 #&#
\begin{eqnarray}\label{system2perttronque}\quad
&&
id\reps{\Psi}(t) + \biggl\{\frac
{ib'}{\epsilon}
\sigma_3\,\frac{\partial\reps{\Psi}(t)}{\partial x}+\frac
{d_0}{2}\,\frac
{\partial^2 \reps{\Psi}(t)}{\partial x^2}+\frac{\gamma_s}{\epsilon
^2}\sigma_3\reps{\Psi}(t)+\frac{i}{\epsilon^2}\gamma_c\reps{\Psi
}(t)\biggr\}\, dt
\nonumber\\
&&\qquad{}
-\frac{\sqrt{\gamma_c}}{\epsilon}\sigma_1\reps{\Psi
}\,d\widetilde{W}_1(t)-
\frac{\sqrt{\gamma_c}}{\epsilon}\sigma_2\reps{\Psi}\,d\widetilde{W}_2(t)
\\
&&\qquad{}+\Theta_R\bigl(\|\reps{\Psi}\|_{L^8(0,t;\mathbb{L}^4
) }
\bigr)f(\reps{\Psi}(t))\,dt = 0\nonumber
\end{eqnarray}
with initial condition $\reps{\Psi}(0)=\psi_0$.
%
%pr2.3 #&#
\begin{prop}\label{existencesystem2}
Let $\reps{\Psi}(0)=\psi_0 \in\mathbb{L}^2(\mathbb{R})$.
Let $T>0$ and $\mathcal{U}^T_c= C([0,T]$; $\mathbb{L}^2)
\cap
L^8(0, T; \mathbb{L}^4)$; then (\ref
{system2perttronquee}) has a unique strong adapted solution $\reps
{\Psi}
\in L^{8}(\mathcal{A}; \mathcal{U}^T_c)$, for any $T>0$.
\end{prop}
\begin{pf}%[Proof of Proposition~\ref{existencesystem2}]
We use a fixed point argument in the Banach space $L^{8}(\mathcal
{A}; \mathcal{U}^T_c )$ for sufficiently small time $T$ depending
on $R$. %where $\mathcal{U}=L^{\infty}(0,T; \mathbb{L}^2)
%the solution afterwards.
We first need to establish estimates on the stochastic integrals
\[
J_{j,\epsilon}\eps{\Psi}(t)=\int_0^t\eps{U}(t-s)\sigma_j\eps
{\Psi
}(s)\,d\widetilde{W}_j(s), \qquad  j=1,2.
\]

%le2.2 #&#
\begin{lemma}\label{majintsto}
Let $T>0$; then for each adapted process $\eps{\Psi} \in L^{8}
(\mathcal{A}; \mathcal{U}_c^T )$ and for $j=1,2$ the stochastic
integral $J_{j,\epsilon}\eps{\Psi}$ belongs to $L^{8}
(\mathcal
{A}; \mathcal{U}_c^T)$. Moreover, for any $T>0$ and $t$ in $[0,
T]$ we have the estimates
\[
\mathbb{E}\bigl(\|J_{j,\epsilon}\eps{\Psi}\|_{L^8(0, T;
\mathbb{L}^4)\cap L^{\infty}(0, T; \mathbb{L}^2
)}^{8}\bigr)
\leq CT^{4}\mathbb{E}\bigl( \|\eps{\Psi}\|_{L^{\infty}
(0,T;\mathbb{L}^{2})}^{8}\bigr).
\]
\end{lemma}
\begin{pf}%[Proof of Lemma~\ref{majintsto}]
Since $\eps{\Psi} \in L^{8}(\mathcal{A}; \mathcal{U}_c^T
)$
and is adapted, we may apply the Burk\-holder--Davis--Gundy inequality in
the Banach space $\mathbb{L}^4(\mathbb{R})$ (which is UMD
space~\cite{brzezniak}):
\begin{eqnarray*}
\mathbb{E}\bigl(\|J_{j,\epsilon}\eps{\Psi}\|_{L^8(0, T;
\mathbb{L}^4)}^{8} \bigr)
&=&\mathbb{E}\biggl( \int_0^T\biggl\|\int_0^t\eps{U}(t-s)\sigma_j\eps
{\Psi}(s)\,d\widetilde{W}_j(s)\biggr\|_{\mathbb{L}^4}^8\,dt\biggr)\\
&\leq&\int_0^T\mathbb{E}\biggl(\sup_{0 \leq u \leq
t}\biggl\|\int_0^u\eps{U}(t-s)\sigma_j\eps{\Psi}(s)\,d\widetilde
{W}_j(s)\biggr\|_{\mathbb{L}^4}^8\biggr)\,dt\\
&\leq& C\mathbb{E}\biggl(\int_0^T \biggl( \int_0^t \|\eps
{U}(t-s)\sigma_j\eps{\Psi}(s)\|_{\mathbb{L}^4}^{2} \,ds\biggr)^4\,dt
\biggr).
\end{eqnarray*}
Using the H\"{o}lder inequality in time, Fubini and a change of variable,
\begin{eqnarray*}
&&
\mathbb{E}\biggl(\int_0^T \biggl( \int_0^t \|\eps{U}(t-s)\sigma
_j\eps{\Psi}(s)\|_{\mathbb{L}^4}^{2} \,ds\biggr)^4\,dt
\biggr)\\
&&\qquad\leq
T^3\mathbb
{E}\biggl(\int_0^T\|\eps{U}(\cdot)\sigma_j\eps{\Psi}(s)\|_{L^8
(0,T;\mathbb{L}^4) }^{8}\,ds \biggr).
\end{eqnarray*}
On the other hand, by Proposition~\ref{strichartz},
\begin{eqnarray*}
\mathbb{E}\biggl(\int_0^T\|\eps{U}(\cdot)\sigma_j\eps{\Psi}(s)\|
_{L^8(0,T;\mathbb{L}^4) }^{8}\,ds \biggr)&\leq&
C\mathbb
{E}
\biggl(\int_0^T\|\eps{\Psi}(s)\|_{\mathbb{L}^2}^{8}\,ds \biggr)\\
&\leq&CT\mathbb{E}\bigl(\|\eps{\Psi}\|_{L^{\infty}
(0,T;\mathbb{L}^2) }^{8}\bigr).
\end{eqnarray*}
Combining these inequalities leads to the estimate in $L^8
(0,T;\mathbb{L}^4) $. The other estimate is proved using the
Burkholder inequality in Hilbert space and the unitary property of
the group $\eps{U}$. Finally, $\eps{U}(t)$ being a unitary semigroup in
$\mathbb{L}^2$, Theorem 6.10 in~\cite{DPZ} tells us that, provided
$ \eps{\Psi} \in L^{8}(\mathcal{A}, L^2(0,T; \mathbb
{L}^2))$, then $J_{j,\epsilon} \eps{\Psi}(\cdot)$ has
continuous modification with values in $\mathbb{L}^2(\mathbb
{R})$.
\end{pf}

Given $\reps{\Psi} \in L^{8}(\mathcal{A}; \mathcal
{U}^T_c)$,
we denote by $\mathcal{T}\reps{\Psi}(t)$ the right-hand side of
(\ref{system2perttronquee}). Since the group $\eps{U}(\cdot)$ maps
$\mathbb
{L}^2(\mathbb{R}) $ into $C(\mathbb{R}, \mathbb
{L}^2(\mathbb{R}) )$, Proposition \ref
{strichartz} and
Lem\-ma~\ref{majintsto} easily imply that the mapping $\mathcal{T}$ maps
$L^{8}(\mathcal{A}; \mathcal{U}^T_c)$ into itself.
%Let us now prove that $\mathcal{T}$ is a contraction mapping in $L^{8}
%(\mathcal{A}; \mathcal{U}^T_c)$, if $T$ is chosen
%sufficiently small (depending on $R$ and $\epsilon$).
Let now $\reps{\Psi}$ and $\reps{\Phi}$ being adapted processes with
values in $L^{8}(\mathcal{A}; \mathcal{U}^T_c)$, then using
Proposition~\ref{strichartz}, the same arguments as in~\cite{bouard}
for the cutoff and Lemma~\ref{majintsto} applied to $J_{j,\epsilon
}(\reps{\Phi}(t)-\reps{\Psi}(t))$, we get
\[
\mathbb{E}(\|\mathcal{T}\reps{\Psi}-\mathcal{T}\reps{\Phi
}\|_{\mathcal{U}^T_c}^{8})^{1/8}\leq\biggl( \frac
{CT}{\epsilon
^2}+ \frac{CT^{1/2}}{\epsilon} + C(R)T^{1/2} \biggr)\mathbb{E}
(\|\reps{\Psi}-\reps{\Phi}\|_{\mathcal{U}^T_c}^{8})^{1/8}.
\]
We conclude that $\mathcal{T}$ is a contraction mapping if $T$ is
chosen such that $CT/\epsilon^2+CT^{1/2}/\epsilon+ C(R)T^{1/2}<1$. As
usual, iterating the procedure, we deduce the existence of a unique
solution of (\ref{system2perttronquee}) in $L^{8}
(\mathcal
{A}; \mathcal{U}^T_c)$ for all $T>0$.~%
\end{pf}

Our aim is now to get global existence for the process $\eps{\Psi}$,
the solution of (\ref{system2pert}) which may be constructed
from the above results. Let us set
\[
\reps{\kappa}(\psi_0,\omega)=\inf\bigl\{ t \geq0,
\|\reps{\Psi}\|_{L^8(0,t;\mathbb{L}^4) }\geq R
\bigr\},
\]
which is a $\eps{\mathcal{G}}(t)$ stopping time. It can be proved using
Strichartz estimates and the integral formulation (\ref
{system2perttronquee}) (see~\cite{bouard,bouardDebussche}) that
$\reps{\kappa}$ is nondecreasing with $R$ and that $\reps{\Psi
}=\eps
{\Psi}^{R'}$ on $[0, \reps{\kappa}]$ for $R<R'$. Thus, we are able to
define a local solution $\eps{\Psi}$ to (\ref{system2pert}) on
the random interval $[0, \eps{\kappa}^*(\psi_0 ))$,
where $
\eps{\kappa}^*(\psi_0 )= \lim_{R \to+\infty} \reps
{\kappa
}$, by setting $\eps{\Psi}(t)=\reps{\Psi}(t)$ on $[0, \reps{\kappa}]$.
It remains to prove that $\eps{\kappa}^*=+\infty$ almost surely. From
the construction of the stopping time $\eps{\kappa}^*$ it is clear
that a.s.,
%
%e2.6 #&#
\begin{equation}\label{alternative1}
\mbox{if } \eps{\kappa}^*(\psi_0 )< +\infty \qquad
\mbox{then }
\lim_{t \nearrow\eps{\kappa}^*(\psi_0 )} \|\reps
{\Psi}\|_{L^8(0,t; \mathbb{L}^4) }= +\infty.
\end{equation}
The arguments are adapted from~\cite{bouard}. We first prove the
following lemma:
%
%le2.3 #&#
\begin{lemma}\label{L2norm}
Let $\eps{\Psi}(0)=\psi_0$ be as in Proposition~\ref{existencesystem2}
and $\reps{\Psi}$ be the corresponding solution of (\ref
{system2perttronque}); then for any $t < T $
\[
\|\reps{\Psi}(t)\|_{\mathbb{L}^2}=\|\psi_0\|_{\mathbb{L}^2}\qquad
\mbox{a.s.},
\]
and there is a constant $\eps{M}>0$, depending on $T$ and $\|\psi_0\|
_{\mathbb{L}^2}$, but independent of~$R$, such that
%
%e2.7 #&#
\begin{equation}\label{boundcont}
\mathbb{E}\bigl(\|\reps{\Psi}\|_{L^8(0,T;\mathbb
{L}^4)}\bigr) \leq\eps{M}(T).
\end{equation}
%
% From this inequality it can be deduced that
% \[
% \eps{\kappa}^*(\psi_0 )= +\infty\mbox{ or } \lim_{t \to
% \]
\end{lemma}
\begin{pf}
To prove that the $\mathbb{L}^2$ norm of the solution $\reps{\Psi}$ of
(\ref{system2perttronque}) is constant in time, we apply formally the
It\^{o} formula to $\frac{1}{2}\|\reps{\Psi}(t)\|_{\mathbb{L}^2}^2$ and
notice that by integration by parts
\[
\biggl( b'\sigma_3 \,\frac{\partial\reps{\Psi}}{\partial x},\reps
{\Psi} \biggr)_{\mathbb
{L}^2}=-\biggl( \reps{\Psi},b'\sigma_3\,\frac{\partial\reps{\Psi
}}{\partial x}
\biggr)_{\mathbb{L}^2}=0.
\]
Since $\sigma_j^*=\sigma_j,  j=1,2,3$, where $*$ stands for the
conjuguate transpose, we get
\[
( \reps{\Psi}(t),i\sigma_j\reps{\Psi}(t) )_{\mathbb{L}^2}=0
\qquad\mbox{for
} j=1,2,3.
\]
Moreover, because the It\^{o} corrections cancel with the damping term
$-\frac{\gamma_c}{\epsilon^2}\reps{\Psi}$ of~(\ref
{system2perttronque}), we get $\|\reps{\Psi}(t)\|_{\mathbb{L}^2}=\|
\psi_0\|_{\mathbb{L}^2}, \forall t \leq T$. The\vspace*{1pt}
computations can be made rigorous by a regularization procedure.

In order to prove (\ref{boundcont}), we follow the procedure in \cite
{bouard,bouardDebussche}. %Let us define the stopping time $\reps{
% \[
% \reps{\widetilde{\kappa}}=\inf\{ t \in[0, \eps{\kappa}^*
%(\psi_0 )), \norm{\eps{\Psi}(t)}{\mathbb{L}^2}^2\geq R
% \]
%We will proceed by contradiction and to this purpose we assume that $
Using the integral formulation (\ref{system2perttronquee}), the
conservation of the $\mathbb{L}^2$-norm\vadjust{\goodbreak} and Proposition~\ref
{strichartz}, we obtain for a.e. $\omega\in\Omega$ and for all time
$T_1$ such that $T \geq T_1>0$
%
%e2.8 #&#
\begin{equation}\label{majalternative}
\|\reps{\Psi}\|_{L^8(0,T_1; \mathbb{L}^4)}
\leq\eps{K}(\omega)+ CT_1^{1/2} \|\reps{\Psi}\|_{L^8
(0,T_1; \mathbb{L}^4) }^{3},
\end{equation}
where
\[
\eps{K}(\omega)=C\biggl(1+\frac{T}{\epsilon^2}\biggr)\|\psi_0\|
_{\mathbb{L}^2}+\frac{1}{\epsilon}\sum_{j=1}^{2}\|J_{j,\epsilon
}\reps{\Psi}\|_{L^8( 0,T; \mathbb{L}^4) }.
\]
%
% Let us set $y= \norm{\eps{\Psi}}{L^8(0,T; \mathbb{L}^4)} $
%and we study the function $h$ defined by
% \[
% h: y \mapsto\eps{K}(\omega)+CT^{1/2}y^3-y,   y \in\mathbb{R}_+
% \]
From inequality (\ref{majalternative}) it follows that $\|\eps{\Psi
}\|_{L^8(0,T_1; \mathbb{L}^4)} \leq2\eps{K}(\omega
)$ if
$T_1$ is chosen, for example, such that $T_1(\omega)=\inf(T,
2^{-6}(C^{1/2}\eps{K})^{-4} )$.
% Indeed $\norm{\eps{\Psi}}{L^8(0,T; \mathbb{L}^4)}
% \[
% CT^{1/2}2^3\eps{K}(\omega)^3 \leq\eps{K}(\omega)
%)^2.
% \]
If $T_1 < T$ we can reiterate the process on small time intervals
$[lT_1, (l+1)T_1] \subset[0, T]$ (keeping $R$ fixed and varying $l$)
to get $\|\eps{\Psi}\|_{L^8(lT_1,(l+1)T_1; \mathbb{L}^4)}
\leq2\eps{K}(\omega)$. Summing these estimates, using $T_1=
2^{-6}C^{-2}( \eps{K})^{-4}$ and the Young inequality, we obtain
\[
\|\reps{\Psi}\|_{L^8(0,T; \mathbb{L}^4)} \leq
C(T)( \eps{K}(\omega)) ^5.
\]
Taking the expectation in the above inequality, using the H\"{o}lder
inequality and Lemma~\ref{majintsto}, we get the following estimate:
%
%e2.9 #&#
\begin{equation}\label{bornecontradiction}\quad
\mathbb{E}\bigl( \|\reps{\Psi}\|_{L^8(0,T; \mathbb
{L}^4)}\bigr)
\leq C(T)\biggl( \biggl(1+\frac{T}{\epsilon^2}\biggr)^5\|\psi
_0\|_{\mathbb{L}^{2}}^5+ \frac{CT^{5/2}}{\epsilon^5}\|\psi_0\|
_{\mathbb{L}^{2}}^5\biggr),
\end{equation}
from which (\ref{boundcont}) follows.
% Now if we assume that
% \[
% \mathbb{Q}( \eps{\kappa}^*(\psi_0 ) < +\infty\
% \]
%Then the probability that there exists a time $T \geq\eps{
%for $R$ sufficiently large is positive. Consequently we deduce from
%inequality \eqref{bornecontradiction} that for all $t$ in $[0, \eps{
\end{pf}

We easily deduce from Lemma~\ref{L2norm} and (\ref{alternative1}) that
$\eps{\kappa}^* =+\infty$ a.s. and as in~\cite{bouard} the existence
and uniqueness of a solution $\eps{\Psi}$ of (\ref{system2pert}), a.s.
in $\mathcal{U}^T_c$ for any $T>0$.

To end the proof of Theorem~\ref{theorem1}, we have to extend those
results to the process~$\eps{X}$. For a.e. $\omega$ in $\mathcal{A}$ and
for each $t \geq0$ we set $\eps{X}(t)=\eps{Z}^{-1}(t)\eps{\Psi
}(t)$. By definition of the process $\eps{Z}^{-1}(t)$ [which, in
particular, is measurable with respect to $\eps{\mathcal{G}}(t)$] and
properties of $\eps{\Psi}$, we easily deduce that $\eps{X}(t)$ is
adapted and continuous with values in $\mathbb{L}^2$, and satisfy
(\ref{manakovPMD}), hence is $C^1$ with values in $\mathbb{H}^{-2}$. By
unitarity of $\eps{Z}$ we also deduce that for all $t \geq0$
\[
\|\eps{\Psi}(t)\|_{\mathbb{L}^2}^2=( \eps{X}(t),\eps
{Z}^{-1}(t)\eps{Z}(t)\eps{X}(t) )_{\mathbb{L}^2}=\|\eps{X}(t)\|
_{\mathbb{L}^2}^2,
\]
and since the coefficients of $\eps{Z}^{-1}(t)$ are a.s. uniformly
bounded, $\eps{X} \in L^8_{\mathrm{loc}}(\mathbb{R}_+,\break \mathbb
{L}^4)$
a.s.; Theorem~\ref{theorem1} is proved.

We now extend the previous global existence results to more regular
initial data. $T$ being fixed, we denote
\[
\mathcal{V}^T=L^{\infty}(0,T; \mathbb{H}^1)\cap
L^8
(0,T; \mathbb{W}^{1,4})
\]
and
\[
\mathcal
{V}^T_c=C(0,T; \mathbb{H}^1)\cap L^8(0,T; \mathbb
{W}^{1,4}).
\]

%pr2.4 #&#
\begin{prop}\label{existenceH1H2}
Let $\eps{\Psi}(0)=\psi_0 \in\mathbb{H}^1$ and let $T>0$; then
equation (\ref{system2pert}) has a unique strong solution $\eps{\Psi}$
with trajectories in $C(0,T; \mathbb{H}^1)$.\vadjust{\goodbreak}
\end{prop}
\begin{pf}%[Proof of Proposition~\ref{existenceH1H2}]
Let $\psi_0$ be in $\mathbb{H}^1$. Given $\reps{\Psi} \in L^8
(\mathcal{A};\mathcal{V}^T )$, we denote by $\mathcal{T}\reps
{\Psi
}(t)$ the right-hand side of (\ref{system2perttronquee}) and $\mathcal
{U}^T=L^{\infty}(0,T; \mathbb{L}^2) \cap L^8(0,T;
\mathbb{L}^4)$. By Proposition~\ref{strichartz}, Lemma \ref
{majintsto} applied to $\partial_x\reps{\Psi}$ and the H\"{o}lder
inequality, we deduce that
\[
\|\mathcal{T}\partial_x\reps{\Psi}\|_{L^8(\mathcal{U}^T)}
\leq C \|\partial_x\psi_0\|_{\mathbb{L}^{2}} +\biggl( \frac
{CT}{\epsilon^2} + \frac{CT^{1/2}}{\epsilon} + CT^{1/2}4R^2\biggr)
\|\partial_x\reps{\Psi}\|_{L^8(\mathcal{U}^T)}.
\]
Therefore, we conclude that choosing $R_0= 2C \|\Psi_0\|_{\mathbb
{H}^1}$, $\mathcal{T}$ maps the closed ball of $L^8(\mathcal
{A};\mathcal{V}^T)$ with radius $R_0$ into itself, provided $T$ is
small enough depending only on $R$ and $\epsilon$, but not on $R_0$.
Combining with the fact that $\mathcal{T}$ is a contraction in $L^8
(\mathcal{A};\mathcal{U}^T)$ and that the balls of $L^8
(\mathcal{A};\mathcal{V}^T)$ are closed for the norm in $L^8
(\mathcal{A};\mathcal{U}^T)$, we conclude to the existence of a unique
fixed point $\reps{\Psi} \in L^8(\mathcal{A};\mathcal {V}^T)$. Using
Proposition~\ref{strichartz} and Lemma \ref {majintsto}, we get
continuity of the solution in $\mathbb{H}^1$. Since the cutoff only
depends on the $L^8(0,T, \mathbb{L}^4 (\mathbb {R} ) )$ norm, we deduce
that there is a unique global solution $\eps{\Psi}$ to
(\ref{system2pert}) with paths in $C ([0,T]; \mathbb{H}^1)$. Since the
transformation $\eps{Z}$ does not depend on $x$, we conclude that these
results still hold true for~$\eps{X}$.
\end{pf}
%
%pr2.5 #&#
\begin{prop}
Let $\eps{\Psi}(0)=\psi_0 \in\mathbb{H}^m$, $m=2,3$. Let $T>0$; then
equation (\ref{system2pert}) has a unique strong solution $\eps{\Psi}$
with paths in $C([0,T]; \mathbb{H}^m)$, $m=2,3$.
\end{prop}
\begin{pf}
We consider equation (\ref{system2perttronque}) but with $\Theta_R(\|
\reps{\Psi}\|_{L^8(0,t; \mathbb{L}^4)})$ replaced by
$\Theta_R(\|\reps{\Psi}(t)\|_{ \mathbb{H}^1}^2)$. Given
$\reps{\Psi}$ in $L^8(\mathcal{A}; L^{\infty}(0,T;
\mathbb
{H}^2( \mathbb{R}) ) )$, we denote by\break
$\mathcal
{T}\reps{\Psi}(t)$ the right-hand side of the integral formulation of
this equation. We easily prove that $\mathcal{T}$ maps the closed ball
of $L^8(\mathcal{A}; L^{\infty}(0,T; \mathbb{H}^2(
\mathbb{R}) ) )$ with radius $R_0$ into itself, for
$R_0= 2C \|\Psi_0\|_{\mathbb{H}^2}$, provided that $T$ is small
enough, depending only on $R$ and $\epsilon$, but not on $R_0$. Using
that this ball is closed for the norm in $L^8(\mathcal{A};
L^{\infty}(0,T; \mathbb{H}^1( \mathbb{R})
)
)$ and that $\mathcal{T}$ is a contraction for the norm in $L^8
(\mathcal{A}; L^{\infty}(0,T; \mathbb{H}^1( \mathbb
{R})
) )$, we deduce that there exists a unique solution $\eps
{\Psi}$ with paths in $C(0,T; \mathbb{H}^2( \mathbb
{R})
)$ a.s., which is global since the solution is global in $\mathbb
{H}^1$. Existence and uniqueness in $\mathbb{H}^3$ can be proved by the
same arguments. Again those results are easily extended to $\eps{X}$
and this concludes the proof of Theorem~\ref{theorem1}.
\end{pf}

%s3 #&#
\section{\texorpdfstring{The limiting equation: Proof of Theorem \protect\ref{theorem2}}
{The limiting equation: Proof of Theorem 1.2}}\label{section4}

In order to prove a local existence and uniqueness result for the
system (\ref{manakovlimite}), we use a compactness approach (see, e.g.,
\cite{flandoli}) motivated by the fact that we do not know
if Strichartz estimates are available for (\ref{manakovlimite}).
Indeed, no transformation similar to the Manakov PMD case seems
to be available, as the equation $dX(t)=-\sqrt{\gamma}\sum_{k=1}^3
\sigma_k\,\frac{\partial X(t)}{\partial x}\,dW_k(t)$ cannot be
solved in a simple way.
% proceed as follows: we construct a random linear semigroup of the
%linear part of Equation \eqref{manakovlimite}
%
% for which no Strichartz estimates are available. Then we use the
%Duhamel formulation and the Lipschitz property in $\mathbb{H}^1$ of
%$F$ to get local existence.
%
% The semigroup is constructed through classical regularization
%procedure: we show the wellposedness of a regularization of Equation
%of martingale to prove that we can construct a weak solution $(
%new filtered probability space $(\widetilde{\Omega}, \widetilde{
%)$. We get back the existence of a unique strong solution on $
%(\Omega, \mathcal{F}, \mathcal{F}_t,\mathbb{P})$ by strong
%unicity of the solution (Yamada-Watanabe Theorem).
We first prove existence of a unique solution in $\mathbb{H}^1$ for the
linear part of the equation, defining then a random propagator, and
then consider the nonlinear part as a perturbation. We will strongly
use the fact that the nonlinearity is locally Lipschitz in $\mathbb
{H}^1$. The regularity in $\mathbb{H}^2$ will follow with the same
arguments as for (\ref{system2pert}).
Let us consider the linear part of (\ref{manakovlimite}),
%
%e3.1 #&#
\begin{eqnarray}\label{lineareq}
dX(t)&=&\biggl( i\,\frac{d_0}{2}\,\frac{\partial^{2} X}{\partial
x^{2}}\biggr)\,dt-\sqrt{\gamma
}\sum
_{k=1}^3 \sigma_k\,\frac{\partial X(t)}{\partial x}\circ
dW_k(t)\nonumber\\[-8pt]\\[-8pt]
&=&\biggl( i\,\frac{d_0}{2}+\frac{3\gamma}{2}\biggr)\,\frac{\partial
^{2} X}{\partial x^{2}}\,dt-\sqrt
{\gamma}\sum_{k=1}^3 \sigma_k\,\frac{\partial X(t)}{\partial x}\,dW_k(t)\nonumber
\end{eqnarray}
with initial data $X(0)=v \in\mathbb{H}^2$. We introduce, for $\eta
>0$, the mollifier $\mol{J} =(I-\eta\,\frac{\partial^{2}
}{\partial x^{2}}
)^{-1}$. We
denote by $\mol{X}$ the solution of the regularized It\^{o} equation
%
%e3.2 #&#
\begin{equation}\label{lineareqreg}\quad
d\mol{X}(t)=\biggl( i\,\frac{d_0}{2}+\frac{3\gamma}{2}\biggr)\,\frac
{\partial^{2} \mol{J}^ 2\mol{X}}{\partial x^{2}}\,dt-\sqrt{\gamma
}\sum_{k=1}^3 \sigma
_k\,\frac{\partial\mol{J}\mol{X}(t)}{\partial x}\,dW_k(t)
\end{equation}
and $\mol{X}(0)=v \in\mathbb{H}^2$. Since the operators $\partial
^2_x\mol{J}^2$ and $\partial_x\mol{J}$ are bounded from $\mathbb{H}^1$
into $\mathbb{H}^1$ (with constants depending on $\eta$), we easily
get, thanks to the Doob inequality, the Fubini theorem, the It\^{o}
isometry and the independence of $(W_k)_{k=1,2,3}$, the
existence and uniqueness of a solution $\mol{X}$ to (\ref{lineareqreg})
with paths in $C([0,T], \mathbb{H}^2)$ for any $T>0$.
Moreover, it is easy to see that the $\mathbb{H}^2$ norm of $\mol{X}$
is conserved since the Pauli matrices are Hermitian. Consequently, the process
\[
\mol{M}(t)=-\mol{X}(t)+\mol{X}(0) + \int_0^t\biggl(\frac{id_0}{2}
+ \frac
{3\gamma}{2}\biggr)\,\frac{\partial^{2} \mol{J}^2\mol{X}}{\partial x^{2}}\,ds
\]
is a $\mathcal{F}_t$ martingale with paths in $C([0,T], \mathbb
{L}^2)$. Let us compute the quadratic variation. Let
$a=(a_1,a_2)^t$ and $b=(b_1,b_2)^t$ be in $\mathbb{L}^{2}$ and $T
\geq t \geq s \geq0$; then
\begin{eqnarray*}
&&\mathbb{E}\bigl(( a,\mol{M}(t) )_{\mathbb{L}^2}(
b,\mol{M}(t) )_{\mathbb{L}^2}-( a,\mol{M}(s) )_{\mathbb{L}^2}(
b,\mol{M}(s) )_{\mathbb{L}^2}|\mathcal{F}_s\bigr)\\
&&\qquad=\gamma\sum_{k=1}^3\mathbb{E}\biggl(\int_s^t\biggl( a,\sigma
_k\,\frac{\partial\mol{J}\mol{X}}{\partial x} \biggr)_{\mathbb
{L}^2}\biggl( b,\sigma_k\,\frac{\partial\mol{J}\mol{X}}{\partial x}
\biggr)_{\mathbb{L}^2}\,du
\Big|\mathcal{F}_s\biggr).
\end{eqnarray*}
We deduce that the quadratic variation of $\mol{M}(t)$ is given by
%
%e3.3 #&#
\begin{equation}\label{qv}\qquad
( b,\llangle\mol{M}(t)\ggangle a )_{\mathbb{L}^2}
=\gamma\sum_{k=1}^3\int_0^t\biggl( a,\sigma_k\,\frac{\partial\mol
{J}\mol{X}}{\partial x} \biggr)_{\mathbb{L}^2}\biggl( b,\sigma_k\,\frac
{\partial\mol{J}\mol{X}}{\partial x}
\biggr)_{\mathbb{L}^2}\,du.
\end{equation}
Using the conservation of the $\mathbb{H}^2$ norm and equation (\ref
{lineareqreg}), we get for all $0\leq\alpha<\frac{1}{2}$
%
%e3.4 #&#
\begin{equation}\label{est}
\mathbb{E}\bigl(\|\mol{X}\|_{C^{\alpha}([0,T]; \mathbb
{L}^{2}) }\bigr)\leq C_{\alpha}(T),
\end{equation}
where $C_{\alpha}(T)$ is a constant independent of $\eta$. Using the
Ascoli--Arzela and Banach--Alaoglu theorems, the Markov inequality and
inequality (\ref{est}), we get that the sequence $(\mathcal {L} (
\mol{X}))_{\eta>0}$ is tight on $ C_w ([0,T],\mathbb {H}^1(\mathbb{R})
) \cap L_w^{\infty}(0,T, \mathbb {H}^2)$. The\vspace*{1pt} Skorokhod
theorem~\cite{billingsley,ethier} implies that on some probability
space $(\widetilde {\Omega}, \widetilde{\mathcal{F}},
\widetilde{\mathcal {F}_t},\widetilde {\mathbb{P}})$, there exist a
sequence of stochastic processes $(\mol{\widetilde{X}})_{\eta>0}$, and
a process $\widetilde {X}$, such that
\[
\mathcal{L}( \mol{\widetilde{X}})=\mathcal{L}(
\mol
{X}),\qquad
\mathcal{L}( \widetilde{X})=\mathcal{L}( X)
\]
and $\lim _{\eta\to0} \mol{\widetilde{X}}=\widetilde{X}$,
$\widetilde{\mathbb{P}}$-a.s. in $C_w([0,T],\mathbb{H}^1)
\cap L_w^{\infty}(0,T, \mathbb{H}^2)$. For all $\eta>0$ and
$t \in[0,T]$ we define the process
%The random variable $\mol{\widetilde{X}}$ is chosen to satisfy the
%martingale problem
%
\[
\widetilde{{M}}_{\eta}(t)=-\mol{\widetilde{X}}(t)+\mol{\widetilde
{X}}(0) +
\int_0^t\biggl(\frac{id_0}{2} + \frac{3\gamma}{2}\biggr)\,\frac
{\partial^{2} \mol{J}^2\mol{\widetilde{X}}}{\partial x^{2}}(s)\,ds.
\]
We deduce from the above laws equality that $\widetilde{{M}}_{\eta}(t)$ is
a square integrable continuous martingale with values in $\mathbb
{L}^{2}$ with respect to the filtration $\widetilde{\mathcal{F}}_t$
%Moreover we proved that there is a unique global strong solution $
%is a continuous Hilbertian martingale, square integrable, with values
%in $\mathbb{H}^{-1}$.
and that the quadratic variation $\llangle\mol{\widetilde{M}}(t) \ggangle$ is
given by formula (\ref{qv}) replacing $\mol{X}$ by $\mol{\widetilde
{X}}$. Let $a \in\mathbb{H}^{1}$, then by the above martingale
property we get for all $s \leq t$
\[
\mathbb{E}\bigl(\bigl( a,\mol{\widetilde{M}}(t)-\mol
{\widetilde{M}}(s) \bigr)_{\mathbb{L}^2}|\widetilde{\mathcal
{F}_s}\bigr)=0.
\]
Using the almost sure convergence in $C_w([0,T],\mathbb
{H}^{1}
(\mathbb{R}) )$ of $\mol{X}$, the boundedness in
$\mathbb
{H}^{-1}$ of the operator $\mol{J}$ and the conservation of the
$\mathbb
{H}^1$ norm, we get the almost sure convergence in $C_w
([0,T],\mathbb{H}^{-1}(\mathbb{R}) )$ of $\mol
{\widetilde{M}}$ to $\widetilde{M}$, where
\[
\widetilde{M}(t)=\widetilde{X}(t)-\widetilde{X}(0)-\int_0^t
\biggl(\frac
{id_0}{2} + \frac{3\gamma}{2}\biggr)\,\frac{\partial^{2} \widetilde
{X}}{\partial x^{2}}(s)\,ds.
\]
Hence, $\widetilde{M}$ is a weakly continuous martingale with values in
$\mathbb{H}^{-1}$. Moreover, using the a.s. convergence in $C_w
([0,T],\mathbb{H}^{1}(\mathbb{R}) )$ and the dominated
convergence theorem, we get for all $ t,s \in[0, T], t \geq s$
and for any $a,b \in\mathbb{H}^1$,
\begin{eqnarray*}
&&
\lim _{\eta\to0} \mathbb{E}( \langle b,\llangle
\mol{\widetilde{M}}(t) \ggangle a \rangle| \widetilde{\mathcal
{F}}_s)\\
&&\qquad=
\gamma\sum_{k=1}^3\mathbb{E}\biggl( \int_0^t\biggl\langle
a,\sigma_k\,\frac{\partial\widetilde{X}}{\partial x}(u)
\biggr\rangle\biggl\langle b,\sigma
_k\,\frac{\partial\widetilde{X}}{\partial x}(u) \biggr\rangle
\,du\Big| \widetilde
{\mathcal{F}}_s\biggr).
\end{eqnarray*}
Thus, the quadratic variation $\langle b,\llangle\widetilde{M}(t) \ggangle a
\rangle$
is given, for all $t \in[0, T]$, by
%
%e3.5 #&#
\begin{equation}\label{quadraticvariationlimit}
\langle b,\llangle\widetilde{M}(t)\ggangle a \rangle
=\gamma\sum_{k=1}^3\int_0^t\biggl\langle a,\sigma_k\,\frac{\partial
\widetilde{X}}{\partial x}(u) \biggr\rangle\biggl\langle b,\sigma_k\,\frac
{\partial\widetilde{X}}{\partial x}(u)
\biggr\rangle \,du.
\end{equation}
Noticing that $\widetilde{M}(0)=0 $ and using the representation
theorem for continuous square integrable martingales, we\vspace*{2pt}
obtain that, on a\vspace*{1pt} possibly enlarged space $(\widetilde{\Omega},
\widetilde {\mathcal{F}}, \widetilde{\mathcal{F}}_t,\widetilde{\mathbb
{P}} )$, one can find a Brownian motion $\widetilde{W} = (\widetilde
{W}_1,\widetilde{W}_2,\widetilde{W}_3 )$ such that
\[
\langle a,\widetilde{M}(t) \rangle=\sqrt{\gamma}\int_0^t \sum_{k=1}^3
\biggl\langle a,\sigma_k\,\frac{\partial\widetilde{X}}{\partial
x}(s) \biggr\rangle \,d\widetilde{W}_k(s).
\]
Thus, we deduce that $(\widetilde{X}, \widetilde{W})$ is a weak
solution of (\ref{lineareq}) on $(\widetilde {\Omega },
\widetilde{\mathcal{F}}, \widetilde{\mathcal{F}}_t,\widetilde {\mathbb
{P}})$ with values in $C_w([0,T],\mathbb{H}^1 (\mathbb {R}) )\cap
L^{\infty}(0,T, \mathbb{H}^2)$. To conclude the proof, we have to prove
pathwise uniqueness of the solution and strong continuity\vspace*{2pt}
in $\mathbb{H}^1$. Since $\widetilde{X} \in L^{\infty}(0,T,
\mathbb{H}^2)$ is the solution of (\ref{lineareq}), we easily deduce
that $\widetilde{X} \in C^{\alpha } ([0,T], \mathbb{L}^2)$ for any
$\alpha\in[0, 1/2)$. By interpolation we obtain that $\widetilde{X} \in
C([0,T], \mathbb {H}^1)$. It follows, using the It\^{o} formula, that
pathwise uniqueness holds for (\ref{lineareq}) in $C([0,T], \mathbb
{H}^1)$. This implies, by the Yamada--Watanabe theorem, that the
solution exists in the strong sense. Thus, we can define a random
unitary propagator $U(t,s)$ which is strongly continuous from $\mathbb
{H}^2$ into $\mathbb{H}^1$. This random propagator can be extended to a
random propagator from $\mathbb{H}^1$ into $\mathbb{H}^1$ using the
continuity of $X$ in $\mathbb{H}^1$, the density of $\mathbb{H}^2$ into
$\mathbb{H}^1$ and the isometry property of $U(t,s)$ in $\mathbb{H}^1$.
%%This conclude that the random propagator $U$ is strongly continuous
%from $\mathbb{H}^1$ to $\mathbb{H}^1$.

The local existence of the nonlinear problem (\ref{manakovlimite}) in
$\mathbb{H}^1$ follows from the construction of the random propagator
$U$: we consider a cutoff function $\Theta\in C_c^{\infty}
(\mathbb
{R})$, $\Theta\geq0$ satisfying
\[
\Theta_R(\|X(t)\|_{\mathbb{H}^1}^2)=
\cases{
1, &\quad if $\|X(t)\|_{\mathbb{H}^1}^2\leq R$,\vspace*{2pt}\cr
0, &\quad if $\|X(t)\|_{\mathbb{H}^1}^2 \geq2R$,}
\]
and first construct a solution $X^R$ of the cutoff equation,
%
%e3.6 #&#
\begin{eqnarray}\label{manakovlimitetronquee}
&&
i\,dX^R(t)+\Biggl( \frac{d_0}{2}\,\frac{\partial^{2} X^R}{\partial
x^{2}} +\Theta_R(\|
X^R(t)\|_{\mathbb{H}^1}^2)F_{}{( X^R) }(t)\Biggr)\,dt\nonumber\\
&&\quad{} +i\sqrt
{\gamma
}\sum_{k=1}^3 \sigma_k\,\frac{\partial X^R(t)}{\partial x}\circ
dW_k(t)\\
&&\qquad=0\nonumber
\end{eqnarray}
with initial data $X^R(0)=v \in\mathbb{H}^1$ and whose integral
formulation is given a.e. by
%
%e3.7 #&#
\begin{equation}\label{duhamellimite}
X^R(t)=U(t,0)v+i \int_0^t \Theta_R( \|X^R(s)\|_{\mathbb
{H}^1}^2)U(t,s)F_{}{( X^R(s)) }\,ds.
\end{equation}
The existence and uniqueness of $X^R\in L^{\rho}(\Omega; C
(0,T;\mathbb{H}^1 ) )$, the solution of (\ref
{duhamellimite}), is easily obtained by a fixed point argument since the
nonlinear term is globally Lipschitz. Introducing the nondecreasing
stopping time
\[
\tau^R = \inf\{ t \geq0, \|X^R(t)\|_{\mathbb{H}^1}^2
\geq R\},
\]
we may then define a local solution $X$ to (\ref
{manakovlimite}) on a random interval $[0, \tau^*(v))$, where $\tau
^*(v)= \lim_{R \to+\infty} \tau^R$ almost surely, by setting
$X(t)=X^R(t)$ on $[0, \tau^R]$. Then for any stopping time $\tau<
\tau
^*$ we have constructed a unique local solution with paths a.s. in
$C([0, \tau], \mathbb{H}^1 ) $.
It follows from the construction of the stopping time $\tau^*$ that if
$\tau^*< +\infty$, then $\limsup_{t \to\tau^*} \|X(t)\|_{\mathbb
{H}^1}= +\infty$. Let us now prove that if $v \in\mathbb{H}^2$, then
the maximal stopping time satisfies the following alternative:
%
%e3.8 #&#
\begin{equation}\label{limitnormh1}
\tau^* = +\infty
\quad\mbox{or}\quad \lim_{t \to\tau^*} \|X(t)\|_{\mathbb
{H}^1}= +\infty.
\end{equation}
We note that the random propagator commutes with derivation. Hence, if
$v \in\mathbb{H}^2$, then $U(\cdot,0)v \in C([0,T], \mathbb
{H}^2
)$. We easily deduce, using (\ref{lineareq}) and interpolating
$\mathbb{H}^1$ between $\mathbb{H}^2$ and $\mathbb{L}^2$, that $U(\cdot,0)v
\in C^{\beta}([0,T], \mathbb{H}^1)$ for $\beta\in[0,1/4)$.
By a fixed point argument in $\mathbb{H}^2$ and equation (\ref
{duhamellimite}), we conclude that $X \in C^{\beta}([0,\tau],
\mathbb{H}^1)$ for any stopping time $\tau< \tau^*$ and for the
same maximal time existence $\tau^*$. Hence, using the condition on
$\tau^*$ and uniform continuity of $X$ in $\mathbb{H}^1$, we get that
(\ref{limitnormh1}) holds.
%
%re3.1 #&#
\begin{rmq}\label{onlocalexistence}
We were not able to prove the global well-posedness for (\ref
{manakovlimite}). Due to the lack of Strichartz estimates, we cannot
control the evolution of the $\mathbb{H}^1$ norm. Even though the
deterministic energy provides a control on the $\mathbb{H}^1$ norm
because we are in the subcritical case, its evolution for a solution of
(\ref{manakovlimite}), which is given in the next lemma,
involves terms which are not well controlled. However, we cannot really
conclude to the real occurrence of blow up or not in this model. It is
clear that on a physical point of view such a phenomenon should not occur.
\end{rmq}
%
%le3.1 #&#
\begin{lemma}
Let the functional $H$ be defined for $u \in\mathbb{H}^1
(\mathbb
{R})$ by
\[
H(u)=\frac{d_0}{4}\int_{\mathbb{R}}\biggl|\frac{\partial u}{\partial
x}\biggr|^2\,dx -\frac
{2}{9}\int_{\mathbb{R}}|u|^4\,dx.
\]
Then for any stopping time $\tau$ such that $\tau< \tau^*$, we have
\begin{eqnarray*}
H(X(\tau)) &=& H( X_0) +\sqrt{\gamma}\frac
{8}{9}\sum_{k=1}^3\int_0^{\tau}\biggl\langle|X|^2X,\sigma_k \,\frac
{\partial X}{\partial x} \biggr\rangle
\circ dW_k(s)\\
&=&H( X_0) +\sqrt{\gamma}\frac{8}{9}\sum_{k=1}^3\int
_0^{\tau
}\biggl\langle|X|^2X,\sigma_k \,\frac{\partial X}{\partial x}
\biggr\rangle
\,dW_k(s)\\
&&{}+\frac{2\gamma}{9} \int_0^{\tau}\int_{\mathbb{R}}(
\partial_x |X_1|^2 +\partial_x |X_2|^2)^2\,dx\,ds
\\
&&{}-\frac{4}{9}\gamma\int_0^{\tau}\int_{\mathbb{R}}\biggl|X_1\,\frac
{\partial X_2}{\partial x}-\frac{\partial X_1}{\partial
x}X_2\biggr|^2\,dx\,ds\\
&&{}+\frac{12}{9}\gamma\int_0^{\tau}\int_{\mathbb{R}}\partial
_x|X_1|^2\,\partial_x|X_2|^2\,dx\,ds.
\end{eqnarray*}
\end{lemma}
\begin{pf}
The first equality follows by Stratonovich differential calculus
applied to the functional $H$ and because the process $X$ is the
solution of (\ref{manakovlimite}). The calculation can be made rigorous
by localization ($H$ is $C^2$ but not bounded) and regularization
through convolution. The second equality is obtained writing the
evolution of $H$ in its It\^{o} formulation, that is,
\begin{eqnarray*}
H(X(\tau)) &=& H( X_0) +\sqrt{\gamma}\frac
{8}{9}\sum_{k=1}^3\int_0^{\tau}\biggl\langle|X|^2X,\sigma_k \,\frac
{\partial X}{\partial x} \biggr\rangle \,dW_k(s)\\
&&{} + \frac{24}{9}\gamma\int_0^{ \tau}\langle X,\partial_xX\Re
(X.\partial_x\overline{X}) \rangle \,ds\\
&&{}- \frac{8}{9}\gamma\sum_{k=1}^3 \int_0^{ \tau}\langle X,\sigma
_k\partial_xX\Re(X.\overline{\sigma}_k\partial_x\overline
{X}) \rangle \,ds,
\end{eqnarray*}
where we used the unitary of the Pauli matrices and $\sigma_k = \sigma
_k^*$, for $k=1,2,3$. Easy calculations lead to the expression given above.
\end{pf}

%s4 #&#
\section{\texorpdfstring{Diffusion limit of the Manakov PMD equation: Proof of Theorem \protect\ref{theorem3}}
{Diffusion limit of the Manakov PMD equation: Proof of Theorem 1.3}}\label{section5}

The aim of this part is the proof of the convergence result given in
Theorem~\ref{theorem3}. For this purpose we have to cutoff equation
(\ref{manakovPMD}) in order to get uniform bounds, with respect to
$\epsilon$, of high order moments of the $\mathbb{H}^2$ norm of the
solution. Let us denote by $\reps{X}$ the solution of the cutoff equation
%
%e4.1 #&#
\begin{equation}\label{manakovPMDtronquee}
\cases{
\displaystyle  i\,\frac{\partial\reps{X}(t)}{\partial t}
+\frac{ib'}{\epsilon}\sigb
(\eps{\nu}(t)) \,\frac{\partial\reps{X}}{\partial
x}+\frac{d_0}{2}\,\frac{\partial^{2} \reps{X}}{\partial
x^{2}}\vspace*{2pt}\cr
\qquad{}+\Theta_R(\|\reps{X}(t)\|_{\mathbb{H}^1}^2
)F_{\eps{\nu }(t)}{( \reps{X}) }=0,\vspace*{2pt}\cr
\displaystyle X_0=v \in\mathbb{H}^3(\mathbb{R}).}
\end{equation}
The proof will consist of the following steps:
\begin{longlist}[(5)]%[\ref{section5}.1]
% \item We first prove that the sequence $(\reps{X})_{
%
\item[(1)] We prove uniform bounds on the solution $\reps{X}$ of (\ref
{manakovPMDtronquee}). These bounds will enable us to prove tightness on
$\mathcal{K}$.
\item[(2)] We use the perturbed test function method to get convergence of
the generators in some sense~\cite{garnier,kushner,papanicolaou}. This method formally gives a candidate for the limit process.
\item[(3)] Setting $\reps{Z}=( \reps{X}, \|\reps{X}(\cdot)\|_{\mathbb
{H}^1}^2)$, we then prove that the family of laws $\mathcal
{L}( \reps{Z})=\mathbb{P}\circ(\reps{Z})^{-1}$ is
tight on $\mathcal{K}$
and we deduce that the process $\reps{Z}$ converges in law, up to a
subsequence.
\item[(4)] Combining the previous steps and using the martingale problem
formulation, we identify the limit and conclude to the weak convergence
of the whole sequence $\reps{X}$.
\item[(5)] Finally, we get rid of the cutoff and we conclude that the
sequence $(\eps{X})_{\epsilon>0}$ converges in law to
$X$ in
$\mathcal{E}(\mathbb{H}^1)$ using the Skorokhod theorem.
\end{longlist}

%
%s4.1 #&#
\subsection{\texorpdfstring{Uniform bounds on $\reps{X}$}{Uniform bounds on X R epsilon}}

Recall that a unique solution $\reps{\Psi} \in C(\mathbb{R}_+,\break
\mathbb{H}^3)$ of the following equation exists (see Section
\ref
{section3}):
%
%e4.2 #&#
\begin{eqnarray}\label{system2perttronqueeH1}\qquad
&&
id\reps{\Psi}(t) + \biggl\{\frac{ib'}{\epsilon} \sigma
_3\,\frac
{\partial\reps{\Psi}(t)}{\partial x}+\frac{d_0}{2}\,\frac{\partial^2
\reps{\Psi}(t)}{\partial x^2}+\frac{\gamma_s}{\epsilon^2}\sigma
_3\reps
{\Psi}(t)+\frac{i}{\epsilon^2}\gamma_c\reps{\Psi}(t)
\biggr\}\, dt
\nonumber\\
&&\qquad{}-\frac{\sqrt{\gamma_c}}{\epsilon}\sigma_1\reps{\Psi
}\,d\widetilde{W}_1(t)-
\frac{\sqrt{\gamma_c}}{\epsilon}\sigma_2\reps{\Psi}\,d\widetilde{W}_2(t)\\
&&\qquad{}+\Theta_R(\|\reps{\Psi}(t)\|_{\mathbb{H}^1}^2
)f(\reps
{\Psi}(t))\,dt = 0.\nonumber
\end{eqnarray}
A solution $\reps{X}$ to (\ref{manakovPMDtronquee}) is then easily
deduced from $\reps{X}(t) = \eps{Z}^{-1}(t)\reps{\Psi}(t)$.

%le4.1 #&#
\begin{lemma}\label{boundinH2}
Let $\psi_0 \in\mathbb{H}^3 $ and $\reps{\Psi}$ be the solution of
(\ref{system2perttronqueeH1}); then for all $T>0$ there exists a
positive constant $C(R,T ) $ independent of $\epsilon$, such
that, a.s. for every $t$ in $[0,T]$,
\[
\|\reps{\Psi}(t)\|_{\mathbb{H}^3} \leq C(R,T ).
\]
Similar bounds hold for $\reps{X}(t)= \eps{Z}^{-1}(t)\reps{\Psi}(t)$
for any $t \in[0,T]$ since $\eps{Z}^{-1}$ is almost surely bounded.
\end{lemma}
\begin{pf}%[Proof of Lemma~\ref{boundinH2}]
The bounds on the $\mathbb{H}^3$ norm are obtained using an energy
method. Using a regularization procedure, the It\^{o} formula applied
to\break
$\|\partial_x\reps{\Psi}(t)\|_{\mathbb{L}^2}^2$ and equation
(\ref{system2perttronqueeH1}), we obtain for all $t \in[0,T]$
\begin{eqnarray*}
\|\partial_x\reps{\Psi}(t)\|_{\mathbb{L}^2}^2 &=& \|\partial_x\psi
_0\|_{\mathbb{L}^2}^2 + 2\int_0^t \langle\partial_x\reps{\Psi
}(s),d\,\partial_x\reps{\Psi}(s) \rangle\\
&&{}+\frac{2\gamma_c}{\epsilon
^2} \int
_0^t \|\partial_x\reps{\Psi}(s)\|_{\mathbb{L}^2}^2\,ds,
\end{eqnarray*}
hence,
\begin{eqnarray*}
\|\partial_x\reps{\Psi}(t)\|_{\mathbb{L}^2}^2 &\leq& \|
\partial_x\psi_0\|_{\mathbb{L}^2}^2\\
&&{}+2\int_0^t\Theta_R(\|
\reps{\Psi}(s)\|_{\mathbb{H}^1}^2)\|\partial_xf
(\reps{\Psi}(s))\|_{\mathbb{L}^2}\|\partial_x\reps{\Psi
}(s)\|_{\mathbb{L}^2}\,ds \\
&\leq& \|\partial_x\psi_0\|_{\mathbb{L}^2}^2+C(R) \int
_0^t\|\partial_x\reps{\Psi}(s)\|_{\mathbb{L}^2}^2\,ds.
\end{eqnarray*}
By the Gronwall lemma we deduce that
\[
\|\partial_x\reps{\Psi}(t)\|_{\mathbb{L}^2}^2 \leq\|\partial
_x\psi_0\|_{\mathbb{L}^2}^2\exp(C(R)T).
\]
Using the same procedure for $\|\partial^2_x\reps{X}\|_{\mathbb
{L}^2}^2$, the Gagliardo--Nirenberg and Young inequalities,
\begin{eqnarray*}
&&\|\partial^2_x\reps{\Psi}(t)\|_{\mathbb{L}^2}^2 -
\|\partial^2_x\psi_0\|_{\mathbb{L}^2}^2\\
&&\qquad\leq C\int_0^t\Theta_R(\|\reps{\Psi}(s)\|_{\mathbb
{H}^1}^2)\bigl(\bigl( \|\reps{\Psi}(s)\|_{\mathbb
{L}^{\infty}}^2+1\bigr) \|\partial^2_x\reps{\Psi}(s)\|_{\mathbb{L}^2}^2\\
&&\hspace*{151pt}{}+
\|\reps{\Psi}(s)\|_{\mathbb{L}^{\infty}}^4\|\partial_{x}\reps
{\Psi}(s)\|_{\mathbb{L}^2}^{6}\bigr)\,ds.
\end{eqnarray*}
By Sobolev embeddings, properties of the cutoff function and again the
Gronwall lemma, we conclude
\[
\|\partial^2_x\reps{\Psi}(t)\|_{\mathbb{L}^2}^2 \leq\|
\partial^2_x\psi_0\|_{\mathbb{L}^2}^2C(R,T).
\]
A bound on $\|\partial^3_x\reps{X}\|_{\mathbb{L}^2}^2$ may be
obtained similarly using the previous estimates and the Gronwall lemma.
\end{pf}
%
%re4.1 #&#
\begin{rmq}
To prove the convergence result, we need initial data in $\mathbb
{H}^3(\mathbb{R})$. We will explain later where exactly we
need this extra regularity, but this is mainly due to the fact that we
prove tightness in $C([0,T], \mathbb{H}^1)$.
\end{rmq}

%re4.2 #&#
\begin{rmq}
Note that we first prove convergence in law for the couple of random
variables $( \reps{X}, \|\reps{X}(\cdot)\|_{\mathbb{H}^1}^2)$.
This is due to the fact that the cutoff is not continuous for the weak
topology in $\mathbb{H}^1$ or for the strong topology in $\mathbb
{H}^1_{\mathrm{loc}}$. These arguments have already been used in \cite
{bouardschema}.
\end{rmq}

%
%s4.2 #&#
\subsection{The perturbed test function method}

Note that the process $\reps{X}$ is not Markov due to the presence of $
\eps{\nu}$. However, $(\reps{X},\eps{\nu})$ is Markov, by construction
of $\nu$. We denote by $\reps{\mathscr{L}}$ its infinitesimal
generator. Let us compute $\reps{\mathscr{L}}f$ for $f$ sufficiently
smooth such that $f$ maps $\mathbb{H}^{-1} \times\mathbb{S}^3$ into
$\mathbb{R}$ and is of class $C^2_b$.
% such that
% \[
% f: \mathbb{H}^{-1} \times\mathbb{S}^3 \to\mathbb{R}.
% \]
Let $\langle\cdot,\cdot\rangle$ be the duality product between $\mathbb
{H}^{1}$ and
$\mathbb{H}^{-1}$. Then, for $\epsilon>0$ and for $\reps{X}$, the
solution of the Manakov PMD equation (\ref{manakovPMDtronquee}),
\begin{eqnarray*}
&&
f(\reps{X}(t), \eps{\nu}(t) ) - f( v,y)\\
&&\qquad=f(\reps{X}(t), \eps{\nu}(t) ) - f(v, \eps{\nu}(t)
) + f(v, \eps{\nu}(t)) - f( v,y)\\
&&\qquad=\langle D_vf(v, \eps{\nu}(t) ),\reps{X}(t)-v
\rangle+
R
(\reps{X}(t),v )\\
&&\qquad\quad{} +f(v, \eps{\nu}(t)) -
f(
v,y),
\end{eqnarray*}
where
\[
R(\reps{X}(t),v )=\int_0^1(1-\theta)\bigl\langle
D^2_{v}f\bigl(v +\theta\bigl(\reps{X}(t)-v\bigr)\bigr) \bigl(
\reps{X}(t)-v \bigr),\reps{X}(t)-v \bigr\rangle \,d\theta
\]
and $D^2_{v}f(v) \in\mathcal{L}( \mathbb
{H}^{-1},\mathbb{H}^{1} )$. Thus,
\begin{eqnarray*}
&&\frac{1}{t}\mathbb{E}\bigl(f(\reps{X}(t),\eps
{\nu
}(t))-f(v,y)|( X(0),\nu(0)) =(
v,y)
\bigr)\\
&&\qquad=\mathbb{E}\biggl(\biggl\langle D_vf(v,\eps{\nu}(t)
),\frac{\reps{X}(t)-v}{t} \biggr\rangle\Big|( X(0),\nu(0)
) =(
v,y
)\biggr)\\
&&\qquad\quad{}+
\mathbb{E}\biggl( \frac{R(\reps{X}(t),v )}{t}
\Big|X(0)=v \biggr) +
\mathbb{E}\biggl( \frac{f(v, \eps{\nu}(t)) -
f(
v,y)}{t}\Big| \nu(0)=y\biggr).
\end{eqnarray*}
We know by Theorem~\ref{theorem1} that if $v \in\mathbb{H}^3$, then
$\reps{X} \in C^1([0,T], \mathbb{H}^{1} ) $. Thus, by the
mean value theorem, equation (\ref{manakovPMDtronquee}), the almost
sure boundedness of $\nu$, Lemma~\ref{boundinH2} and the conservation
of the $\mathbb{L}^2$ norm,
\begin{eqnarray*}
&&
\frac{1}{t}\|\reps{X}(t)-v\|_{\mathbb{L}^{2}} \\
&&\qquad\leq
\sup
_{s\in[0,t]} \|\partial_s\reps{X}(s)\|_{\mathbb{L}^{2}} \\
&&\qquad\leq \sup_{s\in[0,t]}\biggl( \biggl\|\frac{b'}{\epsilon}\sigb
(\eps{\nu}(t)) \,\partial_x\reps{X}(s)\biggr\|_{\mathbb
{L}^2}+ \biggl\|\frac{d_0}{2}\,\partial^2_x\reps{X}(s)\biggr\|_{\mathbb{L}^2}\\
&&\qquad\quad\hspace*{57pt}{} + \|
\Theta_R(\|\reps{X}(s)\|_{\mathbb{H}^1}^2) F_{\eps
{\nu}(s)}{( \reps{X}(s)) } \|_{\mathbb{L}^2}\biggr) \\
&&\qquad\leq \biggl( \frac{b'}{\epsilon} +\frac{d_0}{2}\biggr)
C(R,T)+2RC\|v\|_{\mathbb{L}^2}.
\end{eqnarray*}
Thus, by the boundedness of $D_v^2f$, the continuity of $t \mapsto
\reps
{X}(t)$ in $\mathbb{L}^2$ and the previous bounds, we conclude that
\[
\frac{R(\reps{X}(t),v )}{t} \leq C(R,T,
\epsilon
) \sup_{w \in\mathbb{H}^1}\|D^2_{v}f(w)\|_{\mathcal{L}
( \mathbb{H}^{-1},\mathbb{H}^{1} )}(1+\|v\|_{\mathbb
{L}^2}) \|\eps{X}(t)-v\|_{\mathbb{L}^2}
\]
and the right-hand side above tends to zero as $t$ goes to zero.
Now, we perform the change of variables $t'=t/\epsilon^2$ to get
\[
\frac{1}{t}\mathbb{E}\bigl(f(v,\eps{\nu}(t)
)-f(v,y)|\nu(0)=y\bigr)= \frac{1}{\epsilon^2t'} \mathbb
{E}
\bigl( f(v,\nu(t'))-f(v,y)|\nu(0)=y\bigr).
\]
Thus, using the Markov property of the process $\nu$, and using
(\ref{manakovPMDtronquee}) again, we get an expression of the
infinitesimal generator $\reps{\mathscr{L}}$ of the Markov process
$(\reps{X}, \eps{\nu})$:
%
%e4.3 #&#
\begin{eqnarray}\label{generateur}
\reps{\mathscr{L}}f(v,y)&=&\lim _{t \to0}\frac{1}{t}
\bigl(\mathbb
{E}\bigl(f(\reps{X}(t),\eps{\nu}(t)
)-f(v,y)|
( X(0),\nu(0)) =( v,y) \bigr)\bigr)\nonumber\\
&=& \langle D_vf(v,y),\partial_t\reps{X}(t)|_{t=0}
\rangle+\frac
{1}{\epsilon^2}\mathscr{L}_{\nu}f(v,y)\nonumber\\[-8pt]\\[-8pt]
&=&\biggl\langle D_vf(v,y),\frac{id_0}{2}\,\frac{\partial^{2} v}{\partial
x^{2}}+i\Theta_R(\|
v\|_{\mathbb{H}^1}^2)F_y(v) \biggr\rangle
\nonumber\\
&&{}-\frac{1}{\epsilon}\biggl\langle D_vf(v,y),b'\sigb(y)\,\frac{\partial
v}{\partial x}
\biggr\rangle+\frac
{1}{\epsilon^2}\mathscr{L}_{\nu}f(v,y),\nonumber
\end{eqnarray}
where $\mathscr{L}_{\nu}$ is the infinitesimal generator of $\nu$ and
$\mathscr{D}_{\nu}$ its domain. %We define the set $\reps{
%$f(\cdot,y)$ is of class $C^3$ on $\mathbb{H}^{-1}$ such that $f$ and its
%first three derivatives are bounded on bounded set of $
% \begin{rmq}
% The definition of the set $\reps{\mathscr{D}}$ ensures that the
%brackets, given in the right hand side of \eqref{generateur}, are well
%defined.
% \end{rmq}
% \begin{rmq}
% The hypothesis of linear growth mean that for $v \in\mathbb{H}^1$
% \[
% &&\norm{d_vf(v)}{\mathcal{L}( \mathbb{H}^{2}, \mathbb{R})}
% &&\norm{d_{vv}f(v)}{\mathcal{L}( \mathbb{H}^{2} \times
% &&\norm{D_{vvv}f(v)}{\mathcal{L}( \mathbb{H}^{2} \times
% \]
The perturbed test function method gives (by identifying its
infinitesimal generator) an idea of the limit law of the sequence
$
( \reps{X})_{\epsilon>0} $. It provides in addition convergences
that are useful to prove the weak convergence of the sequence of
measures $(\mathcal{L}( \reps{X})
)_{\epsilon>0}$.
%
%pr4.1 #&#
\begin{prop}[(Perturbed test function method)]\label{fonctiontestperturbee}
There exists a limiting infinitesimal generator $( \mathscr{L}^R
, \mathscr{D}^R)$ such that for all sufficiently smooth and
real-valued functions $f \in\mathscr{D}^R$ and for all positive $
\epsilon$, there exists a test function $\eps{f}$ and positive
constants $C_1(K)$ and $C_2(K)$ satisfying
%
%e4.4 #&#
%e4.5 #&#
\begin{eqnarray}
\label{cvfonctions}
\mathop{\sup_{v \in\mathcal{B}(K)}}_
{y \in\mathbb{S}^3}
|\eps{f}(v,y)-f(v)| &\leq&\epsilon C_1(K),
\\
\label{cvgenerateurs}
\mathop{\sup_{v \in\mathcal{B}(K)}}_{y \in\mathbb{S}^3}
|\reps{\mathscr{L}}\eps{f}(v,y)-\mathscr
{L}^Rf(v)| &\leq&\epsilon C_2(K),
\end{eqnarray}
%
%for every compact set $C$ of $C( [0,T]; \mathbb{H}^1_{
where $\mathcal{B}(K)$ denotes the closed ball of $\mathbb{H}^3
(\mathbb{R})$ with radius $K$.
\end{prop}
\begin{pf}
The idea is to prove that for all suitable test functions $f$, one can
find a function $\eps{f}$ of the form
%
%e4.6 #&#
\begin{equation}\label{ftp}
\eps{f}(v,y)=f(v)+\epsilon f^1(v,y) +\epsilon^2f^2(v,y),
\end{equation}
such that Proposition~\ref{fonctiontestperturbee} holds. We plug\vspace*{1pt} this
expression of $\eps{f}$ into (\ref{generateur}) and formally compute
the expression of $\reps{\mathscr{L}}\eps{f}$:
%
%e4.7 #&#
\begin{eqnarray}\label{generateurmanakov}
\reps{\mathscr{L}}\eps{f}(v,y)&=&\biggl\langle D_vf(v),\frac{id_0}{2}\,\frac
{\partial^{2} v}{\partial x^{2}}+i\Theta_R(\|v\|_{\mathbb
{H}^1}^2)F_y
(v) \biggr\rangle\nonumber\\
&&{}
-\biggl\langle D_vf^1(v,y),b'\sigb(y)\,\frac{\partial v}{\partial x}
\biggr\rangle\nonumber\\
&&{}+\mathscr
{L}_{\nu
}f^2(v,y)+\frac{1}{\epsilon}\mathscr{L}_{\nu}f^1(v,y)-\frac
{1}{\epsilon
}\biggl\langle D_vf(v),b'\sigb(y)\,\frac{\partial v}{\partial x}
\biggr\rangle\nonumber\\[-8pt]\\[-8pt]
&&{}
+\epsilon\biggl\langle D_vf^1(v,y),\frac{id_0}{2}\,\frac{\partial^{2}
v}{\partial x^{2}}+i\Theta
_R(\|v\| _{\mathbb{H}^1}^2)F_y(v) \biggr\rangle\nonumber\\
&&{}
-\epsilon\biggl\langle D_vf^2(v,y),b'\sigb(y)\,\frac{\partial
v}{\partial x} \biggr\rangle
\nonumber\\
&&{}
+\epsilon^2\biggl\langle D_vf^2(v,y),\frac{id_0}{2}\,\frac{\partial^{2}
v}{\partial x^{2}}+i\Theta
_R(\|v\|_{\mathbb{H}^1}^2)F_y(v) \biggr\rangle
\nonumber,
\end{eqnarray}
and we notice that $\mathscr{L}_{\nu}f(v)$ is identically zero because
$f$ does not depend on $\nu=({\nu}_1,{\nu}_2)$. The aim is to wisely
choose the functions $f^1$ and $f^2$ and the regularity of $f$ so that
$\reps{\mathscr{L}}\eps{f}$ is well defined and that $\eps{f}$ and
$\reps{\mathscr{L}}\eps{f}$ converge in the sense of Proposition
\ref
{fonctiontestperturbee}. In particular, we need to cancel the terms
with a factor $1/\epsilon$ and we need the terms with factors
$\epsilon
$ or $\epsilon^2$ to be $\mathcal{O} (\epsilon)$ on bounded sets. In
order to cancel the $1/\epsilon$ terms, we look for a function $f^1$
solution of the Poisson equation
%
%e4.8 #&#
\begin{equation}\label{poisson1}
\mathscr{L}_{\nu}f^1(v,y)= \biggl\langle D_vf(v),b'\sigb(y)\,\frac{\partial
v}{\partial x}
\biggr\rangle.
\end{equation}
By Corollary~\ref{prop3}, we know that
\[
\mathbb{E}_{\Lambda}(g_j(\nu) ) =0\qquad
\forall j=1,2,3.
\]
We deduce that $\langle D_vf(v),b'\sigb(y)\,\frac{\partial
v}{\partial x} \rangle$,
which is a
linear combination of $m_j=g_j(y)$ [see (\ref{matrixsigb})], is of null
mass with respect to the invariant measure $\Lambda$. Hence, $\langle
D_vf(v),b'\sigb(y)\,\frac{\partial v}{\partial x} \rangle$ is a
function of $y \in
\mathbb
{S}^3$, which satisfies the assumptions of Proposition~\ref{prop2},
provided that $f$ is sufficiently smooth, that is, $f \in C^1
(\mathbb{H}^{-1})$ and $v \in\mathbb{L}^2$. It follows that the
solution $f^1$ of the Poisson equation (\ref{poisson1}) can be written as
%
%e4.9 #&#
\begin{eqnarray}\label{solutionpoisson1}
f^1(v,y)&=&\mathscr{L}_{\nu}^{-1}\biggl(\biggl\langle D_vf(v),b'\sigb
(\cdot)\,\frac{\partial v}{\partial x} \biggr\rangle
\biggr)(y)\nonumber\\[-8pt]\\[-8pt]
&=&-\biggl\langle D_vf(v),b'\widetilde{\sigb}(y)\,\frac{\partial
v}{\partial x} \biggr\rangle,\nonumber
% &=&-\int_0^{+\infty}\mathbb{E}(\dual{D_vf(v)}{b'\sigb(
\end{eqnarray}
where
%
%e4.10 #&#
\begin{equation}\label{notation}
\widetilde{\sigb}(y)= \int_0^{+\infty}\mathbb{E}\bigl(
\cdot\sigb
(\nu(t))| \nu(0)=y\bigr) \,dt.
\end{equation}
By Proposition~\ref{prop2}, there is a positive constant $M$ such that
%
%e4.11 #&#
\begin{equation}\label{csteC1}
\Vvert\widetilde{\sigb}(y) \Vvert_{\infty} \leq M\qquad
\forall y \in\mathbb{S}^3,
\end{equation}
and $f^1(v,y)$ is a continuous bounded function of $y$ for $v \in
\mathbb{L}^2$.
We now have to choose the function $f^2$, but we cannot choose
$\mathscr
{L}_{\nu}f^2$ cancelling the terms
\[
\biggl\langle D_vf(v),\frac{id_0}{2}\,\frac{\partial^{2} v}{\partial
x^{2}}+i\Theta_R(\|v\|
_{\mathbb{H}^1}^2)F_y(v) \biggr\rangle
-\biggl\langle D_vf^1(v,y),b'\sigb(y)\,\frac{\partial v}{\partial x}
\biggr\rangle,
\]
because they do not satisfy the null mass condition with respect to
$\Lambda$. Hence, we look for a solution $f^2$ of the Poisson equation
%
%e4.12 #&#
\begin{eqnarray}\label{poisson2}
\mathscr{L}_{\nu}f^2(v,y)&=&-\langle D_vf(v),i\Theta_R(\|v\|
_{\mathbb{H}^1}^2)F_y(v) \rangle\nonumber\\
&&{}+\langle
D_vf(v),i\Theta_R(\|v\|_{\mathbb{H}^1}^2)F
(v) \rangle\nonumber\\[-8pt]\\[-8pt]
&&{}+\biggl\langle D_vf^1(v,y),b'\sigb(y)\,\frac{\partial
v}{\partial x}
\biggr\rangle\nonumber\\
&&{}-\mathbb{E}_{\Lambda}
\biggl(\biggl\langle D_vf^1(v,y),b'\sigb(y)\,\frac{\partial v}{\partial x}
\biggr\rangle\biggr),\nonumber
\end{eqnarray}
where, due to (\ref{solutionpoisson1}),
%
%e4.13 #&#
\begin{eqnarray}\label{difff1}
&&
\biggl\langle D_vf^1(v,y),b'\sigb(y) \,\frac{\partial
v}{\partial x}
\biggr\rangle\nonumber\\
&&\qquad= -(b')^2\biggl\langle
D^2_vf(v)\widetilde{\sigb}(y)\,\frac{\partial v}{\partial
x},\sigb(y)\,\frac{\partial v}{\partial x}
\biggr\rangle\\
&&\qquad\quad{}-(b')^2\biggl\langle D_vf(v),\widetilde{\sigb}(y)\sigb
(y)\,\frac{\partial^{2} v}{\partial x^{2}} \biggr\rangle.\nonumber
\end{eqnarray}
Moreover, thanks to expression (\ref{difff1}), the Fubini theorem and
Corollary~\ref{prop3},
%
%e4.14 #&#
\begin{eqnarray}\label{gelimite2}
&&
-\mathbb{E}_{\Lambda}\biggl(\biggl\langle D_vf^1(v,y),b'\sigb
(y)\,\frac{\partial v}{\partial x} \biggr\rangle\biggr)\nonumber\\
&&\qquad=
(b')^2\sum_{j,k=1}^3\biggl\langle D^2_{v}f(v) \sigma_k\,\frac
{\partial v}{\partial x},\sigma_j\,\frac{\partial
v}{\partial x} \biggr\rangle\int_0^{+\infty}\mathbb
{E}_{\Lambda
}
( g_k(\nu(t)) g_j(\nu(0))
)\,dt\nonumber\\[-8pt]\\[-8pt]
&&\qquad\quad{}
+(b')^2\sum_{j,k=1}^3\biggl\langle D_vf(v),\sigma_k\sigma
_j\,\frac{\partial^{2} v}{\partial x^{2}} \biggr\rangle\int_0^{+\infty
}\mathbb{E}_{\Lambda}
(g_k(\nu
(t))g_j(\nu(0)))\,dt\nonumber\\
&&\qquad=\frac{\gamma}{2} \sum_{k=1}^3\biggl\langle D^2_{v}f(v) \sigma_k\,\frac
{\partial v}{\partial x},\sigma_k\,\frac{\partial
v}{\partial x} \biggr\rangle
+\frac{3\gamma}{2}\biggl\langle D_vf(v),\frac{\partial^{2} v}{\partial
x^{2}} \biggr\rangle,\nonumber
\end{eqnarray}
where $\gamma= ( b')^2/6\gamma_c $.
Provided that $f$ is of class $C^2(\mathbb{H}^{-1})$ and $v
\in\mathbb{H}^1$ and because $f^1(v,\cdot)$ is of class $C^2_b
(\mathbb
{S}^3)$ for any $v \in\mathbb{H}^1$, we can now define, by
Proposition~\ref{prop2}, a unique solution, up to a constant, to the
Poisson equation (\ref{poisson2}). This solution $f^2$ is expressed as
% \begin{eqnarray}\label{solutionpoisson2}
% f^2(v,y)&=&\int_0^{+\infty}\mathbb{E}(\dual{D_vf(v)}{i
%)-F(v)) }|\nu(0)=y)\,dt\nonumber\\&&-
%( \nu(t)) \pds{}{x}{v}}-\mathbb{E}_{\Lambda}(
% &=&\mathscr{L}_{\nu}^{-1}(\dual{D_vf(v)}{i\Theta_R(\norm{v}{
%) })\\
% &&-\mathscr{L}_{\nu}^{-1}(\dual{D_vf^1(v,y)}{b'\sigb( y
%) \pds{}{x}{v}}-\mathbb{E}_{\Lambda}(\dual{D_vf^1(v,y)}{b'
% \end{eqnarray}
%
%e4.15 #&#
\begin{eqnarray}\label{solutionpoisson2}
f^2(v,y)&=&\mathscr{L}_{\nu}^{-1}\bigl(\bigl\langle D_vf(v),i\Theta
_R(\|v\|_{\mathbb{H}^1}^2)\bigl( F_{y}(v
)-F(v)\bigr) \bigr\rangle\bigr)\nonumber\\
&&{}-\mathscr{L}_{\nu}^{-1}\biggl(\biggl\langle D_vf^1(v,y),b'\sigb(
y) \,\frac{\partial v}{\partial x} \biggr\rangle\nonumber\\
&&\hspace*{40.5pt}{}-\mathbb
{E}_{\Lambda}\biggl(\biggl\langle
D_vf^1(v,y),b'\sigb(y) \,\frac{\partial v}{\partial
x} \biggr\rangle
\biggr)\biggr)\nonumber\\[-8pt]\\[-8pt]
&=&\langle D_vf(v),i \Theta_R(\|v\|_{\mathbb{H}^1}^2)
\widetilde{F}(v,y) \rangle\nonumber\\
&&{}-(b')^2 \sum_{k,l=1}^3\biggl\langle D_v^2f(v)\sigma_k\,\frac
{\partial v}{\partial x},\sigma_l\,\frac{\partial
v}{\partial x} \biggr\rangle
\wwidetilde{g}_{k,l}(y)\nonumber\\
&&\hspace*{0pt}{}-\biggl\langle D_vf(v),( b')^2 \widetilde
{\widetilde{\sigb}}(y) \,\frac{\partial^{2} v}{\partial x^{2}}
\biggr\rangle,\nonumber
\end{eqnarray}
where
\[
\widetilde{F}(v,y) = \int_0^{+\infty}\mathbb{E}\bigl(
F_{\nu (t)}{( v) } -F_{}{( v) } | \nu(0) =y \bigr) \,dt
\]
and
\[
\wwidetilde{g}_{k,l}(y) = \int_0^{+\infty} \biggl( \int
_t^{+\infty}\mathbb{E}\bigl( g_k(\nu(s))
g_l(\nu
(t)) | \nu(0) =y \bigr) \,ds - \frac{\gamma}{2(b')^2}
\delta
_{kl}\biggr)\,dt
\]
and
\[
\widetilde{\widetilde{\sigb}}(y) = \int_0^{+\infty}\biggl( \int
_t^{+\infty}\mathbb{E}\bigl( \sigb(\nu(s))
\sigb
(\nu(t)) | \nu(0) =y \bigr) \,ds - \frac{3\gamma
}{2(b')^2}\biggr)\,dt.
\]
Replacing $\mathscr{L}_{\nu}f^1$ and $\mathscr{L}_{\nu}f^2$ in
(\ref{generateurmanakov}), respectively, by the right-hand side of
(\ref{poisson1}) and (\ref{poisson2}), and using expression (\ref
{gelimite2}), we get
%
%e4.16 #&#
\begin{eqnarray}\label{gelimite}
\reps{\mathscr{L}}\eps{f}(v,y)&=&\biggl\langle D_vf(v),\biggl(\frac
{id_0}{2}+\frac{3\gamma}{2}\biggr)\,\frac{\partial^{2} v}{\partial
x^{2}} +i\Theta_R
(\|v\|_{\mathbb{H}^1}^2)F_{}{( v) } \biggr\rangle\nonumber\\
&&{}+\frac{\gamma
}{2}\sum_{k=1}^3
\biggl\langle D^2_{v}f(v) \sigma_k\,\frac{\partial v}{\partial
x},\sigma_k\,\frac{\partial v}{\partial x}
\biggr\rangle\nonumber\\
&&{}+ \epsilon\biggl\langle D_vf^1(v,y),\frac{id_0}{2}\,\frac{\partial^{2}
v}{\partial x^{2}}+i
\Theta_R(\|v\|_{\mathbb{H}^1}^2) F_{y}{( v) } \biggr\rangle
\\
&&{}- \epsilon\biggl\langle D_vf^2(v,y),b'\sigb(y)\,\frac
{\partial v}{\partial x}
\biggr\rangle\nonumber\\
&&{}+ \epsilon^2\biggl\langle D_vf^2(v,y),\frac{id_0}{2}\,\frac{\partial^{2}
v}{\partial x^{2}}+i
\Theta_R(\|v\|_{\mathbb{H}^1}^2)F_{y}{( v) } \biggr\rangle
\nonumber,
\end{eqnarray}
and we define the limiting operator by
%
%e4.17 #&#
\begin{eqnarray}\label{generateurlimite}
\mathscr{L}^Rf(v)
&=&\biggl\langle D_vf(v),\biggl(\frac{id_0}{2}+\frac{3\gamma}{2}
\biggr)\,\frac{\partial^{2} v}{\partial x^{2}} +i\Theta_R(\|v\|
_{\mathbb{H}^1}^2)F_{}{( v) } \biggr\rangle\nonumber\\[-8pt]\\[-8pt]
&&{}+\frac{\gamma}{2}\sum_{k=1}^3 \biggl\langle D^2_{v}f(v) \sigma_k\,\frac
{\partial v}{\partial x},\sigma_k\,\frac{\partial
v}{\partial x} \biggr\rangle.\nonumber
\end{eqnarray}
Hence, if we define $\mathscr{D}^R$ as the space of functions which are
the restriction to $\mathbb{H}^3$ of functions $f$ from $\mathbb
{H}^{-1}$ into $\mathbb{R} $ of class $C^3(\mathbb
{H}^{-1})$
and such that $f$ and its first three derivatives are bounded on
bounded sets of $\mathbb{H}^{-1}$, then the functions $f^1$ and $f^2$
are well defined for $f \in\mathscr{D}^R$. Moreover, if $f \in
\mathscr
{D}^R$, then $\reps{\mathscr{L}}\eps{f}$ is well defined for $v \in
\mathbb{H}^3$.

We now write that
\[
\mathop{\sup_{v \in\mathcal{B}(K)}}_{y \in\mathbb{S}^3}
|\eps{f}(v,y)-f(v)|
\leq\epsilon
\mathop{\sup_{v \in\mathcal{B}(K)}}_{y \in\mathbb{S}^3}
|f^1(v,y)|+\epsilon^2
\mathop{\sup_{v \in\mathcal{B}(K)}}_{y \in\mathbb{S}^3}
|f^2(v,y)|
\]
and use the following result, which is proved in Section~\ref{A2}.
%
%le4.2 #&#
\begin{lemma}\label{firstcv}
Let $f \in\mathscr{D}^R$ and $f^1$ and $f^2$ be, respectively,
solutions of (\ref{poisson1}) and (\ref{poisson2}). Then
% \lim _{\epsilon\to0}( \epsilon\underset{y \in
%
\[
\mathop{\sup_{v \in\mathcal{B}(K)}}_{y \in\mathbb{S}^3}
|f^1(v,y)| \leq C_1(K) \quad\mbox{and}\quad
\mathop{\sup_{v \in\mathcal{B}(K)}}_{y \in\mathbb{S}^3}
|f^2(v,y)| \leq C_2(K).
\]
\end{lemma}

This proves the first convergence of Proposition~\ref{fonctiontestperturbee}.
With $\mathscr{L}^Rf(v)$ given by (\ref{generateurlimite}), the second
convergence (\ref{cvgenerateurs}) in Proposition \ref
{fonctiontestperturbee} follows from (\ref{gelimite}) and the next
lemma, which is proved in Section~\ref{A2}.
\end{pf}
%
%le4.3 #&#
\begin{lemma}\label{secondcv}
Let $f \in\mathscr{D}^R$ and $f^1$, $f^2$ be, respectively, solutions
of (\ref{poisson1}) and (\ref{poisson2}). Then
%&&\lim _{\epsilon\to0}\underset{y \in\mathbb{S}^2}{
%&&\lim _{\epsilon\to0} \underset{y \in\mathbb{S}^2}{
%
\begin{eqnarray*}
\mathop{\sup_{v \in\mathcal{B}(K)}}_{y \in\mathbb{S}^3}
\biggl|\biggl\langle D_vf^1(v,y),\frac{id_0}{2}\,\frac{\partial^{2} v}{\partial
x^{2}}+i \Theta_R(\|
v\|_{\mathbb{H}^1}^2) F_{y}{( v) } \biggr\rangle\biggr| &\leq& C_1(K),\\
\mathop{\sup_{v \in\mathcal{B}(K)}}_{y \in\mathbb{S}^3}
\biggl|\biggl\langle D_vf^2(v,y),b'\sigb(y)\,\frac{\partial
v}{\partial x} \biggr\rangle\biggr|
&\leq&
C_2(K),\\
\mathop{\sup_{v \in\mathcal{B}(K)}}_{y \in\mathbb{S}^3}
\biggl|\biggl\langle D_vf^2(v,y),\frac{id_0}{2}\,\frac{\partial^{2} v}{\partial
x^{2}}+i \Theta_R(\|
v\|_{\mathbb{H}^1}^2)F_{y}{( v) } \biggr\rangle\biggr| &\leq& C_3(K).
\end{eqnarray*}
\end{lemma}

%
%s4.3 #&#
\subsection{\texorpdfstring{Tightness of the family of probability measures $(\mathcal{L}(\reps{Z}))_{\epsilon>0}$}
{Tightness of the family of probability measures (L(Z R epsilon))epsilon>0}}

To prove tightness on $\mathcal{K}$ of the sequence of probability
measure $\mathcal{L}( \reps{Z}) = \mathbb{P}\circ
(\reps{Z})^{-1}$, we
need to obtain uniform bounds in $\epsilon$ on $\reps{Z}$ in the space
\[
\bigl( C([0,T], \mathbb{H}^{2})\cap
C^{\alpha
}([0,T], \mathbb{H}^{-1}) \bigr) \times
C^{\delta
}([0,T], \mathbb{R})
\]
for suitable $\alpha, \delta> 0$. Note that uniform bounds of $\reps
{X}$ in $C([0,T], \mathbb{H}^{2})$ are given by
Lemma~\ref{boundinH2}. The perturbed test function method will enable
us to get the uniform bound in $C^{\alpha}([0,T],
\mathbb{H}^{-1})$. Such bounds cannot be directly obtained using
(\ref{manakovPMDtronquee}) because of the $1/\epsilon$ term.
In order to obtain such bounds, we use again the perturbed test
function method for convenient test functions. Let $(\widetilde
{e}_j)_{j \in\mathbb{N}^*}$ be a complete orthonormal system in
$\mathbb{L}^{2}$. %such that $\widetilde{e_j}\in\mathbb{H}^1$ for any
%$j \in\mathbb{N}^*$.
Recall that $\langle\cdot,\cdot\rangle $ is the duality product between
$\mathbb
{H}^{1}$--$\mathbb{H}^{-1}$ and $(\cdot,\cdot)_{\mathbb{L}^2}$ the inner
product in
$\mathbb
{L}^2$. By definition of $\mathbb{H}^s, s \in\mathbb{R}$, we can define
a complete orthonormal system $(e_j)_{j \in\mathbb{N}^*}$ on $\mathbb
{H}^1$ from $(\widetilde{e}_j)_{j \in\mathbb{N}^*}$,
\begin{eqnarray*}
\|v\|_{\mathbb{H}^{-1}}^2
&=&\|( 1+\xi^2)^{-1/2}\widehat{v} \|_{\mathbb{L}^2}^2\\
&=&\sum_{j=1}^{+\infty}\bigl(( 1+\xi^2)^{-1/2}\widehat{v},
\widehattilde{e}_j \bigr)_{\mathbb{L}^2}^2\\
&=&\sum_{j=1}^{+\infty}\langle e_j,v \rangle^2,
\end{eqnarray*}
where $e_j=\mathcal{F}^{-1}( ( 1+\xi^2)^{-1/2}
\widehattilde{e}_j)$ %\in\mathbb{H}^2$
for any $j \in\mathbb{N}^*$. We denote by $(f_j)_{j \in
\mathbb{N}^*}$ the family of test functions in $\mathscr{D}^R$
defined by
\begin{eqnarray*}
f_j\dvtx \mathbb{H}^{-1} &\rightarrow& \mathbb{R},\\
v &\mapsto& f_j(v)=\langle e_j,v \rangle.
\end{eqnarray*}
For $v \in\mathbb{H}^3$, we also consider particular perturbed test
functions $f_{j,\epsilon}$ of the form
%
%e4.18 #&#
\begin{equation}\label{newphi}
f_{j,\epsilon}(v,y)=f_j(v)+\epsilon f_j^1(v,y),
\end{equation}
where, for all $j$ in $\mathbb{N}^*$, $ f_j^1(v,y)=\langle
e_j,\varphi^1(v,y) \rangle$ for a given function $\varphi^1$ with
values in $\mathbb
{H}^{2}$. %Applying the operator $\reps{\mathscr{L}}$ given by
%all $t$ in $[0, T]$ and $j \in\mathbb{N}^*$:
% \[
% \reps{\mathscr{L}}\compoeps{f}{j}(v, y ) &=& \dual{e_j}{
%)\nl{y}{v}}- \dual{e_j}{D_v\varphi^1(v,y).b'\sigb(y)
% && + \epsilon\dual{e_j}{D_v\varphi^1(v,y).(\frac{id_0}{2}
% \]
We now choose $\varphi^1$ as a solution of the Poisson equation in $y$:
%
%e4.19 #&#
\begin{equation}\label{poisson3}
\mathscr{L}_{\nu}\varphi^1(v, y )-b'\sigb
(y) \,\frac{\partial v}{\partial x}=0,
\end{equation}
whose explicit formulation is given by (see Proposition~\ref{prop2})
%
%e4.20 #&#
\begin{equation}\label{solpoisson3}
\varphi^1(v, y )=-b'\widetilde{\sigb}(y)\,\frac{\partial
v}{\partial x},
\end{equation}
where $\widetilde{\sigb}(y)$ is given by (\ref{notation}). We point out
that $\varphi^1$ behaves in its first variable like $\frac{\partial
}{\partial x}$ and
is linear in $v$. Consequently, for all $j$ in $\mathbb{N}^*$,
%
%e4.21 #&#
\begin{eqnarray}\label{generateurtension}
&&\reps{\mathscr{L}} f_{j,\epsilon}(\reps{X}(t),\eps{\nu
}(t)
)\nonumber\\
&&\qquad=\biggl\langle e_j,\frac{id_0}{2}\,\frac{\partial^{2} \reps
{X}(t)}{\partial x^{2}}+i\Theta
_R(\|\reps{X}(t)\|_{\mathbb{H}^1}^2 )F_{\eps{\nu
}(t)}{( \reps{X}(t)) } \biggr\rangle\nonumber\\
&&\qquad\quad{} +\biggl\langle e_j,
(b'
)^2 \widetilde{\sigb}(\eps{\nu}(t))\sigb(\eps
{\nu}(t))\,\frac{\partial^{2} \reps{X}(t)}{\partial x^{2}}
\biggr\rangle\\
&&\qquad\quad{}- \epsilon\biggl\langle e_j,b'\widetilde{\sigb}(\eps{\nu
}(t))\,\frac{\partial}{\partial x}\biggl(\frac{id_0}{2}
\,\frac{\partial^{2} \reps{X}(t)}{\partial x^{2}}\nonumber\\
&&\qquad\quad\hspace*{105pt}{} +i \Theta_R
(\|\reps{X}(t)\|_{\mathbb{H}^1}^2
)F_{\eps{\nu}(t)}{( \reps{X}(t)) }\biggr) \biggr\rangle\nonumber.
\end{eqnarray}
For all $t \in[0,T]$, we define the process $\reps{M}$ with values in
$\mathbb{H}^{-1}$ given for any $j$ in $\mathbb{N}^*$ by
\begin{eqnarray*}
&&\langle e_j,\reps{M}(t) \rangle\\
&&\qquad=f_{j,\epsilon}(\reps{X},\eps{\nu})(t)-f_{j,\epsilon
}(v,y)-\int_0^t\reps{\mathscr{L}} f_{j,\epsilon}(\reps
{X}(s),\eps{\nu}(s))\,ds\\
&&\qquad=\langle e_j,\reps{X}-v \rangle+\epsilon\langle e_j,\varphi
^1(\reps{X},\eps{\nu})-\varphi^1(v,y) \rangle\\
&&\qquad\quad{}-\int_0^t\reps
{\mathscr{L}}f_{j,\epsilon}
(\reps
{X}(s),\eps{\nu}(s))\,ds.
\end{eqnarray*}
Given the fact that $\reps{\mathscr{L}}$ is the infinitesimal generator
of the continuous Markov process $( \reps{X}, \eps{\nu}
) $
and $\reps{\mathscr{L}}f_{j,\epsilon}$ is well defined because $f_j
\in\mathscr{D}^R$, then $\langle e_j,\reps{M}(t) \rangle$ is a real-valued
continuous martingale. Moreover, it is a square integrable martingale,
as follows from the bounds on the $\mathbb{H}^3$ norm of $\reps{X}$
obtained in Lemma~\ref{boundinH2}. To prove tightness of the family of
probability measures $\mathcal{L}( \reps{Z})$ on
$\mathcal
{K}$, we need estimates of moments on the processes $\reps{X}$ and
$\|\reps{X}(\cdot)\|_{\mathbb{H}^1}^2$. Before proving these estimates we
introduce a process $\reps{Y}$ close in probability to $\reps{X}$ for
which it will be easier to get those estimates, using, in particular,
the Kolmogorov criterion. The idea is to use Lemma~\ref{procheproba}
below to get tightness of the family $\mathcal{L}( \reps
{Z}
)$ from convergence in law of a subsequence of $\reps{Y}$.
%
%le4.4 #&#
\begin{lemma}\label{procheproba}
Let us define the process $\reps{Y}$ as
%
%e4.22 #&#
\begin{equation}\label{processy}\quad
\reps{X}(t)- \reps{Y}(t) = \epsilon\bigl(\varphi^1(v,y
)-\varphi^1(\reps{X}(t), \eps{\nu}(t))\bigr)\qquad
\forall
t \in[0,T];
\end{equation}
then for all $\delta>0 $,
\[
\mathbb{P}\bigl( \|\reps{X}-\reps{Y}\|_{C([0,T
],\mathbb{H}^1)}>\delta\bigr)\leq\frac{\epsilon
}{\delta}C_1(T,R).
\]
\end{lemma}
\begin{pf}%[Proof of Lemma~\ref{procheproba}]
Using the Markov inequality and Lemma~\ref{boundinH2}, we get for all
$\delta>0$,
\begin{eqnarray*}
&&\mathbb{P}\Bigl(\sup_{t\in[0,T]}\|\reps{X}(t)-\reps
{Y}(t)\|_{\mathbb{H}^1}>\delta\Bigr)\\
&&\qquad\leq \frac{\epsilon}{\delta}\mathbb{E}\Bigl(\sup_{t\in
[0,T]}\|\varphi^1(\reps{X}(t),\eps{\nu}(t))-\varphi
^1(v,y)\|_{\mathbb{H}^1}\Bigr)\\
&&\qquad\leq \frac{\epsilon}{\delta}\mathbb{E}\biggl(\sup_{t\in
[0,T]}\biggl\|b'\widetilde{\sigb}(\eps{\nu}(t))\,\frac
{\partial\reps{X}(t)}{\partial x} -b'\widetilde{\sigb}
(y)\,\frac{\partial v}{\partial x}\biggr\|_{\mathbb{H}^1}
\biggr) \\
%&\leq& \frac{\epsilon}{\delta}\mathbb{E}(\sup_{t\in[0,T]}
%( b'\eps{\sigb}(u)\pds{}{x}{v} | \nu(0)=y )\,du }{
&&\qquad\leq\epsilon\frac{2M}{\delta}C(T,R),
\end{eqnarray*}
where $M$ is given by (\ref{csteC1}).
\end{pf}

Note that the process $\reps{Y}$ is also defined by the identity, for
all $j$ in $\mathbb{N}^*$,
%
%e4.23 #&#
\begin{eqnarray}\label{secondformulationforY}
\langle e_j,\reps{Y}(t) \rangle&=&\langle e_j,\reps{X}(t) \rangle-
\epsilon\langle e_j,\varphi^1(v,y)-\varphi^1
(\reps{X}(t), \eps{\nu}(t)) \rangle\nonumber\\
&=&\langle e_j,\reps{M}(t) \rangle+\langle e_j,v \rangle+\int
_0^t\reps{\mathscr
{L}}f_{j,\epsilon}(\reps{X}(s),\eps{\nu}(s))\,ds\\
&&\eqntext{\forall t \in[0,T].}
\end{eqnarray}

%le4.5 #&#
\begin{lemma}\label{lemmatension}
For all $1 \geq\epsilon>0$, there exist three positive constants
$C_1(T,\break R)$, $C_2(T,R)$ and $C_3(T,R)$ depending on final time $T$ and
on the cutoff radius $R$, but independent of $\epsilon$, such that
%
%e4.24 #&#
%e4.25 #&#
%e4.26 #&#
\begin{eqnarray}
\label{majoration1}
\mathbb{E}\bigl(\|\reps{Y}\|_{C([0,T],\mathbb{H}^2
)}^4\bigr) &\leq& C_1(T,R), \\
\label{majoration2}
\mathbb{E}\bigl(\|\reps{Y}\|_{C^{\alpha}([0,T
],\mathbb{H}^{-1})}\bigr)&\leq& C_2(T,R),\\
\label{majoration3}
\mathbb{E}\bigl(\|\|\reps{Y}\|_{\mathbb{H}^1}^2\|_{C^{\delta
}([0,T],\mathbb{R})}\bigr)&\leq& C_3(T,R),
\end{eqnarray}
where $0<\alpha<\frac{1}{2}$ and $\delta=\alpha/3 >0$.
\end{lemma}
\begin{pf}%[Proof of Lemma~\ref{lemmatension}]
Thanks to Lemma~\ref{boundinH2}, we know that the solution $\reps{X}$
of (\ref{manakovPMDtronquee}) is uniformly bounded, for all
$\epsilon$, in $\mathbb{H}^3$ by a constant $C$ depending on $R$ and~$T$.
We conclude, using the explicit formulation of $\varphi_1$ given
by (\ref{solpoisson3}) and (\ref{processy}), that (\ref
{majoration1}) holds.

To prove inequality (\ref{majoration2}), we first need an intermediate
estimate that will be proved in Section~\ref{A2}.
%
%le4.6 #&#
\begin{lemma} \label{lemmakolmogorov}
There exists a positive constant $C(R,T)$ such that for all $t,s \in[0,T]$
\[
\mathbb{E}\bigl(\|\reps{Y}(t)-\reps{Y}(s)\|_{\mathbb
{H}^{-1}}^4
\bigr)\leq C(R,T)(t-s)^2.
\]
\end{lemma}

Then we deduce from Lemma~\ref{lemmakolmogorov}
\[
\mathbb{E}\bigl(\|\reps{Y}\|_{\mathbb{W}^{\gamma,4}([0, T],
\mathbb{H}^{-1} ) }^4\bigr)\leq C(R,T)
\]
for any $\gamma< 1/2$. We use the Sobolev embedding $ \mathbb
{W}^{\gamma,4}([0, T], \mathbb{H}^{-1} ) \hookrightarrow
C^{\alpha}([0, T]$, $\mathbb{H}^{-1} )$ for $\gamma-\alpha>
1/4$ and $\gamma< 1/2$, which implies $\alpha<1/4$. Thus, we deduce the
second inequality (\ref{majoration2}).

It remains to prove the last bound (\ref{majoration3}). Note that for
$t,s \in[0,T]$
\begin{eqnarray*}
&&\bigl|\|\reps{Y}(t)\|_{\mathbb{H}^1}^2-\|\reps{Y}(s)\|_{\mathbb
{H}^1}^2\bigr|\\
&&\qquad\leq C\sup_{r \in[0,T]}\|\reps{Y}(r)\|_{\mathbb{H}^1}\|
\reps{Y}(t)-\reps{Y}(s)\|_{\mathbb{H}^1}\\
&&\qquad\leq C\sup_{r \in[0,T]}\|\reps{Y}(r)\|_{\mathbb{H}^1}\sup_{r
\in[0,T]}\|\reps{Y}(r)\|_{\mathbb{H}^2}^{2/3}\|\reps{Y}(t)-\reps
{Y}(s)\|_{\mathbb{H}^{-1}}^{1/3}.
\end{eqnarray*}
It follows that if $\delta= \alpha/3$,
\[
\|\|\reps{Y}(\cdot)\|_{\mathbb{H}^1}^2\|_{C^{\delta}([0,T],
\mathbb{R})}\leq C \sup_{r \in[0,T]}\|\reps{Y}(r)\|
_{\mathbb{H}^2}^{5/3}\|\reps{Y}\|_{C^{\alpha}([0,T],\mathbb
{H}^{-1})}^{1/3}.
\]
Inequality (\ref{majoration3}) is then implied by the H\"{o}lder
inequality, (\ref{majoration1}) and (\ref{majoration2}).
\end{pf}
%
%re4.3 #&#
\begin{rmq}
The extra $\mathbb{H}^3$ regularity is needed precisely in the first
step of the above proof\vspace*{1pt} in order to estimate the $\mathbb{H}^2$ norm of
$\reps{Y}$, which involves the gradient of $\reps{X}$.
\end{rmq}
%
%pr4.2 #&#
\begin{prop}\label{theoremtension}
The family of laws $(\mathcal{L}( \reps{Z})
)_{\epsilon>0}$ is tight on $\mathcal{K}$.
\end{prop}
\begin{pf}%[Proof of Proposition~\ref{theoremtension}]
We set $\reps{\widetilde{Z}} =(\reps{Y}, \|\reps{Y}(\cdot)\|
_{\mathbb{H}^1}^2)$. Denoting by $\mathcal{B}(K)$ the closed ball of $
(C([0,T];\mathbb{H}^2(\mathbb{R}) )\cap C^{\alpha
}([0,T];\mathbb{H}^{-1}(\mathbb{R}))) \times C^{\delta}([0,T];
\mathbb{R})$ with radius $K$, for $\alpha$ and $\delta$ as in Lemma
\ref{lemmatension}, we deduce using the Ascoli--Arzela and
Banach--Alaoglu theorems that $\mathcal{B}(K)$ is compact in
$\mathcal{K}$. Using the Markov inequality and Lem\-ma~\ref
{lemmatension}, we get
\begin{eqnarray*}
&&\mathbb{P}\bigl(\reps{\widetilde{Z}} \notin\mathcal{B}(K)\bigr)\\
&&\qquad\leq\frac{1}{K}\mathbb{E}\bigl(\max\bigl\{\|\reps{Y}\|
_{C([0,T];\mathbb{H}^2)}, \|\reps{Y}\|_{C^{\alpha
}([0,T];\mathbb{H}^{-1})}, \|\|\reps{Y}\|_{\mathbb
{H}^1}^2\|_{C^{\delta}([0,T])} \bigr\}\bigr)\\
&&\qquad\leq\frac{1}{K}\max(C_1^{1/4}(T,R),C_2(T,R),
C_3(T,R)).
\end{eqnarray*}
We conclude that the family of laws $(\mathcal{L}( \reps
{\widetilde{Z}}))_{\epsilon>0}$ is tight on $\mathcal{K}$
and by the Prokhorov theorem we obtain the relative compactness of the
sequence of laws $(\mathcal{L}( \reps{\widetilde
{Z}}
))_{\epsilon>0}$, that is, up to a subsequence, the sequence
$\mathcal{L}( \reps{\widetilde{Z}})$ weakly converges
to a
probability measure $\mathcal{L}( \widehat{Z}^R)$ where
$\widehat{Z}^R=(\widehat{X}^R, \gamma^R)$.
%we can find subsequences of $\mesure{(\reps{Y}, \norm{
%on $\mathcal{K}\times C([0,T], \mathbb{R})$ that we denote $
We may now use Lem\-ma~\ref{procheproba} to prove that the family of laws
$\mathcal{L}( \reps{Z})$ is tight. Indeed, it easily follows
from Lemma~\ref{procheproba} and the above convergence in law that for
all $g \in C_b(\mathcal{K})$
\[
\lim _{\epsilon\to0}\mathbb{E}(g(\reps{Z}
)
)=\mathbb{E}(g(\widehat{Z}^R)).
\]
\upqed\end{pf}
%
%s4.4 #&#
\subsection{\texorpdfstring{Convergence in law of the process $\eps{X}^R$}
{Convergence in law of the process X R epsilon}}
\label{subproblememartingale}

In order to get the convergence in law of the whole sequence $
(\reps{X})_{\epsilon>0}$, it remains to characterize the limit,
that is, to prove that $ \widehat{X}^R = X^R$, the solution\vspace*{1pt} of
(\ref{manakovlimitetronquee}), and that $\gamma^R(t)= \|X^R(t)\|
_{\mathbb{H}^1}^2$ for any $t \in[0,T]$. The tool here will be
the use of the martingale problem formulation introduced by
Stroock and Varadhan in~\cite{Stroock}.

%pr4.3 #&#
\begin{prop}\label{finalcv}
The whole sequence $\reps{X}$ converges in law to $X^R$ in $C
([0,T], \mathbb{H}^1)$.
\end{prop}
\begin{pf}%[Proof of the Proposition~\ref{finalcv}]
In order to prove that any subsequence of $\reps{X}$ converges to the
same limit $X^R$, the solution of (\ref{manakovlimitetronquee}), we
will prove the convergence of the martingale problem for suitable test
functions $f \in\mathscr{D}^R$. To this purpose let us define, for $a
\in\mathbb{H}^{1}$ with compact support, the particular test function
$f_a(\cdot) = \langle a,\cdot\rangle$, so that \mbox{$f_a \in\mathscr{D}^R$}. From
this particular choice, we construct a perturbed test function
$f_{a,\epsilon}$,
\[
f_{a,\epsilon}(v,y)=f_a(v)+\epsilon f^1_a(v,y)+\epsilon^2 f^2_a(v,y),
\]
obtained thanks to Proposition~\ref{fonctiontestperturbee}. The
correctors $f_a^1$ and $f_a^2$ are chosen to be the solution of the
Poisson equations (\ref{poisson1}) and (\ref{poisson2}) for $f_a$.
% \begin{eqnarray*}
% \mathscr{L}_{\nu}f_a^1(v,y)= \dual{a}{b'\sigb(y)\pds{}{x}{v}},
% \]
% i.e $f_a^1(v,y)= - \dual{a}{b'\widetilde{\sigb}(y)\pds{}{x}{v}}$ and
% \begin{eqnarray*}
% \mathscr{L}_{\nu}f_a^2(v,y)&=&-\dual{a}{i\Theta_R(\norm{v}{
% - \dual{a}{(b')^2(\widetilde{\sigb}(y)\sigb(y)- \frac{3
% \]
% Using \eqref{gelimite}, we get
% \begin{eqnarray*}
% \reps{\mathscr{L}}\compoeps{f}{a}(v,y)&=&\dual{a}{(
% &&+ \epsilon\dual{D_vf_a^1(v,y)}{\frac{id_0}{2}\pds{2}{x}{v}+i
% &&- \epsilon\dual{D_vf_a^2(v,y)}{b'\sigb(y)\pds{}{x}{v}}\\
% &&+ \epsilon^2\dual{D_vf_a^2(v,y)}{\frac{id_0}{2}\pds{2}{x}{v}+i
% \]
% because $D^2_vf_a(v)=0$. Thus $\reps{\mathscr{L}}\compoeps{f}{a}$ is
%well defined because $f_a \in\mathscr{D}^R$.
Let us denote by $\reps{Z}$ a subsequence converging to $\widehat{Z}^R$
and define the $\mathbb{H}^{-1}$ valued process $\reps{\mathbf
{N}}
(\reps{Z}(t))$, associated to (\ref{manakovPMDtronquee}),
\begin{eqnarray*}
\langle a,\reps{\mathbf{N}}(\reps{Z}(t)) \rangle
&=&f_{a,\epsilon
}(\reps{X}(t),\eps{\nu}(t))-f_{a,\epsilon}(v,y) \\
&&{}-\int
_0^t\reps
{\mathscr{L}} f_{a,\epsilon}(\reps{X}(s),\eps{\nu}(s))\,ds,
\end{eqnarray*}
where $\reps{\mathscr{L}}$ is given by (\ref{generateurmanakov}).
% Consequently we get
% \begin{eqnarray}\label{diffmartingale}
% \lefteqn{\dual{a}{\reps{\mathbf{N}}(\reps{Z}(t))-
% \epsilon(f_a^1(\reps{X}(t), \eps{\nu}(t))-f_a^1
%(v, y)) + \epsilon^2 (f_a^2(\reps{X}(t),
% -\epsilon\int_0^t\dual{D_vf_a^1(\reps{X}(s), \eps{\nu}(s)
%)}{\frac{id_0}{2}\pds{2}{x}{\reps{X}(s)}+i\Theta_R(\norm{
%)}\,ds\nonumber\\&&+\epsilon\int_0^t\dual{D_vf_a^2(
% -\epsilon^2\int_0^t\dual{D_vf_a^2(\reps{X}(s), \eps{\nu}(s)
%)}{\frac{id_0}{2}\pds{2}{x}{\reps{X}(s)}+i\Theta_R(\norm{
%)}\,ds\nonumber,
% \end{eqnarray}
% where
% \[
% \dual{a}{\mathbf{N}^R(\reps{Z}(t))}=\dual{a}{
%) \pds{2}{x}{\reps{X}(s)}+i\Theta_R(\norm{\reps{X}(s)}{
% \]
We also define the process $\mathbf{N}^R(\reps{Z}(t))$,
\[
\langle a,\mathbf{N}^R(\reps{Z}(t)) \rangle=f_a(
\reps{X}(t)
)-f_a(v) -\int_0^t \mathscr{L}^Rf_a(\reps
{X}(s))\,ds,
\]
where $\mathscr{L}^R$ is given by expression (\ref{generateurlimite}).
%&&+ \epsilon^2\dual{D_vf^2(\reps{X}(t),\eps{\nu}(t))}{
%and where
Moreover, we denote by $\mathscr{L}_{\gamma^R}^R$ the operator whose
expression is given by (\ref{generateurlimite}) replacing\vspace*{1pt}
$\|\widehat {X}^R(t)\|_{\mathbb{H}^1}$ by $\gamma^R(t)$ in the cutoff
function. Let us now define $\langle a,\mathbf{N}^R(\widehat {Z}^R(t))
\rangle$ by
%
%e4.27 #&#
\begin{equation}\label{martingalelimiteZ}
\langle a,\mathbf{N}^R(\widehat{Z}^R(t)) \rangle
=f_a(
\widehat
{X}^R(t) )-f_a(v) -\int_0^t \mathscr{L}_{\gamma
^R}^Rf_a( \widehat{X}^R(s))\,ds.
\end{equation}
%
%and prove the following convergence result: for all $t,s \in[0,T]$
%and $0\leq t_1\ldots<t_m\leq s\leq t$ and for all test
%functions $h_1,\ldots, h_m \in C_b(\mathbb{H}^1_{\mathrm{loc}}
%)}\prod_{j=1}^mh_j(\reps{Z}(t_j)))}\\
%&=&\mathbb{E}(\dual{a}{\mathbf{N}^R(\widehat{Z}^R(t))-
The process $\langle a,\reps{\mathbf{N}}(\reps{Z}(t))
\rangle$
is a
real continuous martingale because $(\reps{X}, \eps{\nu})$
is a Markov process and because $\reps{\mathscr{L}}f_{a,\epsilon}$ is
well defined since $\reps{X}(t) \in\mathbb{H}^3$. Moreover, it is a
square integrable martingale, as follows from the bounds on the
$\mathbb
{H}^3$ norm of $\reps{X}$ obtained in Lemma~\ref{boundinH2}. The above
martingale property implies that for all $t,s \in[0, T], t \geq s$,
\[
\mathbb{E}[ \langle a,\reps{\mathbf{N}}(\reps
{Z}(t))-\reps{\mathbf{N}}(\reps{Z}(s)) \rangle
|\sigma
(\reps
{Z}(u), \eps{\nu}(u)), u\leq s ]=0.
\]
%
% The Markov property of the couple $(\reps{X}, \eps{m})$
%allows us to reduce the conditioning to
% \[
% \mathbb{E}[( \eps{f}(\reps{X},\eps{m})(t)-
%)\,ds) |(\reps{X}, \eps{m})(s) ]=0
% \]
It follows, in particular, that for all test functions $h_1,\ldots, h_m
\in C_b(\mathbb{H}^1_{\mathrm{loc}}\times\mathbb{R})$ and
$0\leq t_1<\cdots<t_m\leq s\leq t$,
\[
\mathbb{E}\Biggl[\langle a,\reps{\mathbf{N}}(\reps
{Z}(t))-\reps{\mathbf{N}}(\reps{Z}(s)) \rangle
\prod_{j=1}^mh_j(\reps
{Z}(t_j))\Biggr]=0.
\]
Using Proposition~\ref{fonctiontestperturbee}, Lemma~\ref{boundinH2}
and the boundedness of the functions $h_j$, we get
\begin{eqnarray*}
& & \mathbb{E} \Biggl(\langle a,\reps{\mathbf{N}}(\reps
{Z}(t))-\mathbf{N}^R(\reps{Z}(t))-\reps{\mathbf
{N}}(\reps{Z}(s))+\mathbf{N}^R(\reps
{Z}(s)) \rangle\prod
_{j=1}^mh_j(\reps{Z}(t_j))\Biggr)\\
& &\qquad\leq\epsilon C(R,T).
\end{eqnarray*}
%
%%%%%%%%%%%%%%%%%%%%%%%%%%%
% Let us now consider the term
% \[\lefteqn{\mathbb{E}(\dual{a}{\mathbf{N}^R(\reps{Z}(t)
%)-\mathbf{N}^R(\reps{Z}(s))}\chi_{R_0}(\reps{Z}
%)\prod_{j=1}^mh_j(\reps{Z}(t_j)))}\\
% &-&\mathbb{E}(\dual{a}{\mathbf{N}^R(\widehat{Z}^R(t))-
%)\prod_{j=1}^mh_j(\widehat{Z}^R(t_j))),
% \]
Let us consider a cutoff function $\chi_{R_0} \in C_c^{\infty
}(\mathcal
{K})$ satisfying
\[
\chi_{R_0}(u) = \cases{
1, &\quad if $u \in\mathcal{B}_{\mathcal{K}}(R_0)$,\cr
0, &\quad if $u \notin\mathcal{B}_{\mathcal{K}}(2R_0)$,}
\]
where $\mathcal{B}_{\mathcal{K}}(R_0)$ denotes the closed
ball of radius $R_0$ of the space $\mathcal{K}$ and $R_0$ is chosen
such that $\reps{X} \in\mathcal{B}_{\mathcal{K}}(R_0)$ a.s.
(see Lemma~\ref{boundinH2}). Note that by continuity of the functions
$\chi_{R_0}$ and $\{h_j\}_{j\in\{1,\ldots,
m\}
}$, respectively, in $\mathcal{K}$ and $\mathbb{H}^{1}_{\mathrm
{loc}}\times
\mathbb{R}$, by continuity of $f_a(\cdot)$ for the weak topology in
$\mathbb
{H}^{1}$, by continuity and boundedness of $\Theta_R$ in $C([0,T];
\mathbb{R})$, by continuity of $F$ from $\mathbb{H}^1$ to
$\mathbb
{H}^{-1}$ and the bounds on $F(\reps{X}(t))$ obtained thanks
to Lemma~\ref{boundinH2}, the function
\[
\langle a,\mathbf{N}^R(\reps{Z}(t)) \rangle\chi
_{R_0}
(\reps
{Z})\prod_{j=1}^mh_j(\reps{Z}(t_j))
\]
is a bounded and continuous function of $\reps{Z}$ from $\mathcal{K}$
into $\mathbb{R}$. We deduce by convergence in law of $\reps{Z}$ to
$\widehat{Z}^R$ in $\mathcal{K}$, since the test function $a$ is
compactly supported, that for all $t,s \in[0,T], t \geq s$
%
%e4.28 #&#
\begin{equation}\label{martingalelimite}\quad
\mathbb{E}\Biggl(\langle a,\mathbf{N}^R(\widehat
{Z}^R(t))-\mathbf{N}^R(\widehat{Z}^R(s)) \rangle
\chi_{R_0}
(\widehat
{Z}^R)\prod_{j=1}^mh_j(\widehat{Z}^R(t_j))\Biggr)=0.
\end{equation}
Since, almost surely, $\reps{X}$ belongs to the closed ball $\mathcal
{B}_{\mathcal{K}}(R_0)$, we deduce that almost surely
$\widehat{X}^R \in\mathcal{B}_{\mathcal{K}}(R_0)$.
Thus, we
conclude from (\ref{martingalelimite}) that $\langle a,\break\mathbf
{N}^R(\widehat{Z}^R(\cdot)) \rangle$ is a continuous square
integrable martingale
with respect to the filtration $\mathcal{G}_t=\sigma(\widehat
{Z}^R(s), s \leq t)$ and this holds for any $a \in\mathbb
{H}^1$ with compact support.

In order to identify the equation satisfied by $\widehat{X}^R$, we
consider, for $a,b \in\mathbb{H}^1$ with compact support, the function
$g_{a,b}(v)=f_a(v)f_b(v) \in\mathscr{D}^R$ and the perturbed test
function $g_{a,b,\epsilon}$,
\[
g_{a,b,\epsilon}(v,y)=g_{a,b}(v) +\epsilon g_{a,b}^1(v,y) +\epsilon
^2g_{a,b}^2(v,y),
\]
obtained thanks to Proposition~\ref{fonctiontestperturbee}. Thus,
functions $g_{a,b}^1(v,y)$ and $g_{a,b}^2(v,y)$ are chosen to be
solutions of the Poisson equations (\ref{poisson1}) and (\ref
{poisson2}) for $g_{a,b}$.
% \[
% \mathscr{L}_{\nu}g_{a,b}^1(v,y)= f_b(v)\dual{a}{b'\sigb(y)
% \]
% i.e $g_{a,b}^1(v,y)= - f_b(v)\dual{a}{b'\widetilde{\sigb}(y)
% \[
% \mathscr{L}_{\nu}g_{a,b}^2(v,y)&=&-f_b(v)\dual{a}{i\Theta_R(
%(v)}\\&&-f_a(v)\dual{b}{i\Theta_R(\norm{v}{
% +\dual{D_vg^1_{a,b}(v,y) }{b'\sigb(y)\pds{}{x}{v}}+\frac{3
% +\gamma\sum_{k=1}^{3}\dual{a}{\sigma_k\pds{}{x}{v}}\dual{b}{\sigma_k
% \]
% Thus using \eqref{gelimite}
% \[
% \reps{\mathscr{L}}g_{a,b,\epsilon}(v,y)&=& f_b(v) \dual{a}{(
% +f_a(v) \dual{b}{( \frac{id_0}{2} +\frac{3\gamma}{2})
%) }\\&&
% +\gamma\sum_{k=1}^{3}\dual{a}{\sigma_k\pds{}{x}{v}}\dual{b}{\sigma_k
% &&+ \epsilon\dual{D_vg^1_{a,b}(v,y)}{\frac{id_0}{2}\pds{2}{x}{v}+i
% &&- \epsilon\dual{D_vg^2_{a,b}(v,y)}{b'\sigb(y)
% &&+ \epsilon^2\dual{D_vg^2_{a,b}(v,y)}{\frac{id_0}{2}\pds{2}{x}{v}+i
% \]
Let us now define the real-valued continuous martingale
\begin{eqnarray*}
\mathbf{H}_{a,b,\epsilon}^R(\reps{Z}(t)
)&=&{g}_{a,b,\epsilon
}(\reps{X}(t),\eps{\nu}(t))-g_{a,b,\epsilon}(v,y)\\
&&{} -\int_0^t\reps
{\mathscr{L}} g_{a,b,\epsilon}(\reps{X}(s),\eps{\nu
}(s))\,ds.
\end{eqnarray*}
Using the same arguments as before, we may prove that
\begin{eqnarray*}
&&\lim _{\epsilon\to0}\mathbb{E}\Biggl(\bigl(
\mathbf{H}_{a,b,\epsilon}^R(\reps{Z}(t))-\mathbf
{H}_{a,b,\epsilon}^R(\reps{Z}(s))\bigr)\chi
_{R_0}(\reps
{Z}
) \prod_{j=1}^mh_j(\reps{Z}(t_j))\Biggr)\\
&&\qquad=\mathbb{E}\Biggl(\bigl( \mathbf{H}_{a,b}^R(\widehat
{Z}^R(t)
)-\mathbf{H}_{a,b}^R(\widehat{Z}^R(s))\bigr)\chi
_{R_0}
(\widehat{Z}^R) \prod_{j=1}^mh_j(\widehat{Z}^R(t_j))\Biggr),
\end{eqnarray*}
where
\begin{eqnarray*}
\mathbf{H}_{a,b}^R(\widehat{Z}(t))&=&{g}_{a,b}(\widehat
{X}^R(t)) - {g}_{a,b}(v) \\
&&{}-\int_0^t\mathscr{L}_{\gamma^R}^R
{g}_{a,b}(\widehat{X}^R(s))\,ds.
\end{eqnarray*}
From the above convergence and the martingale property of $\mathbf
{H}_{a,b,\epsilon}^R(\reps{Z}(t))$, we deduce that
$\mathbf
{H}_{a,b}^R(\widehat{Z}^R(\cdot))$ is a continuous real-valued
martingale. A classical computation then shows that the quadratic
variation of the martingale $\mathbf{N}^R(\widehat
{Z}^R(t))$
defined in (\ref{martingalelimiteZ}) is given by
\begin{eqnarray*}
&& \langle b,\llangle\mathbf{N}^R(\widehat{Z}^R(t)) \ggangle a
\rangle\\
&&\qquad = \int_0^t\mathscr{L}_{\gamma^R}^R(f_a(\widehat
{Z}^R(s)) f_b(\widehat{Z}^R(s)) ) -f_a
(\widehat{Z}^R(s))\mathscr{L}_{\gamma^R}^Rf_b(\widehat
{Z}^R(s))\\
&&\qquad\quad{} -f_b(\widehat{Z}^R(s))\mathscr
{L}_{\gamma
^R}^Rf_a(\widehat{Z}^R(s))\,ds.
\end{eqnarray*}
Applying the operator $\mathscr{L}_{\gamma^R}^R$, respectively, to the
test functions $f_a$ and $g_{a,b}$, we obtain that
\[
\mathscr{L}_{\gamma^R}^Rf_a(\widehat{Z}^R(t))
=\biggl\langle a,\biggl(\frac{id_0}{2}+\frac{3\gamma}{2}\biggr)\,\frac
{\partial^{2} \widehat{X}^R}{\partial x^{2}} +i\Theta_R(\|
\gamma^R(t)\|_{\mathbb
{H}^1}^2)F_{}{( \widehat{X}^R) } \biggr\rangle
\]
and
\begin{eqnarray*}
&& \mathscr{L}_{\gamma^R}^Rg_{a,b}( \widehat{Z}^R(t))
\\
&&\qquad= f_b( \widehat{X}^R(t))\biggl\langle a,\biggl( \frac
{id_0}{2} +\frac{3\gamma}{2}\biggr) \,\frac{\partial^{2} \widehat
{X}^R(t)}{\partial x^{2}}+i\Theta_R(\gamma^R(t))F
(\widehat
{X}^R(t)) \biggr\rangle\\
&&\qquad\quad{}
+f_a( \widehat{X}^R(t)) \biggl\langle b,\biggl( \frac{id_0}{2}
+\frac{3\gamma}{2}\biggr) \,\frac{\partial^{2} \widehat
{X}^R(t)}{\partial x^{2}}+i\Theta
_R(\gamma^R(t))F(\widehat{X}^R(t))
\biggr\rangle\\
&&\qquad\quad{}
+\gamma\sum_{k=1}^{3}\biggl\langle a,\sigma_k\,\frac{\partial\widehat
{X}^R(t)}{\partial x} \biggr\rangle\biggl\langle b,\sigma_k\,\frac{\partial
\widehat{X}^R(t)}{\partial x}
\biggr\rangle.
\end{eqnarray*}
We deduce that the quadratic variation is given by formula (\ref
{quadraticvariationlimit}) with $\widetilde{X}$ replaced by
$\widehat
{X}^R$. Thus, using the martingale representation theorem, we can write
the $\mathcal{G}_t$-martingale $\mathbf{N}^R (\widehat
{Z}^R(t))$ as the stochastic integral
\[
\langle a,\mathbf{N}^R (\widehat{Z}^R(t)) \rangle=\sqrt
{\gamma
}\int
_0^t \sum_{k=1}^3 \biggl\langle a,\sigma_k\,\frac{\partial\widehat
{X}^R(s)}{\partial x}
\biggr\rangle \,dW_k(s),
\]
where $W=(W_1,W_2,W_3)$ is a real-valued Brownian motion on
a possibly enlarged space $(\Omega,\mathcal{G},\mathcal{G}_t,
\mathbb
{P})$. We deduce that $( \widehat{X}^R, W)$ is a weak
solution in $C([0,T]; \mathbb{H}_{\mathrm
{loc}}^1)\cap
C_w([0,T]; \mathbb{H}^1) \cap L_w^{\infty
}
(0,T; \mathbb{H}^2)$ of the equation
%
%e4.29 #&#
\begin{equation}\label{identifyinglimit}
\cases{
\displaystyle id\widehat{X}^R(t)+\biggl( \frac{d_0}{2}\partial_x^2{\widehat{X}^R(t)}
+\Theta_R(\gamma^R(t))F_{}{( \widehat{X}^R(t)) }
\biggr)\,dt \vspace*{2pt}\cr
\displaystyle\qquad{} +i\sqrt{\gamma}{ \sum_{k=1}^3} \sigma
_k\,\partial_x{\widehat{X}^R(t)}\circ dW_k(t)=0,\vspace*{2pt}\cr
X_0=v \in\mathbb{H}^3.}
\end{equation}
The next step consists in proving that almost surely $\gamma^R(t)=
\|\widehat{X}^R(t)\|_{\mathbb{H}^1}^2$. Using the Skorokhod representation
theorem, we can construct new random variables (that we still denote
$\reps{Z}$, $\widehat{Z}^R$) on a new common probability space $
(\Omega, \mathcal{F}, \mathcal{F}_t,\mathbb{P})$ with,
respectively, $\mathcal{L}( \reps{Z}) $ and $\mathcal
{L}
( \widehat{Z}^R)$ as probability measures and with values in
$\mathcal{K}$ such that
\[
\lim _{\epsilon\to0}\reps{Z}= \widehat{Z}^R , \qquad\mathbb
{P}
\mbox{-a.s. in } \mathcal{K}.
\]
Since $\widehat{X}^R \in L^{\infty}(0,T; \mathbb{H}^2
)$, we
deduce using (\ref{identifyinglimit}) that $\widehat{X}^R
\in
C([0,T]; \mathbb{L}^2)$. Hence, applying the It\^{o}
formula, it is easy to see, since $\Theta_R$ is a real-valued function,
that almost surely
\[
\|\widehat{X}^R(t)\|_{\mathbb{L}^2}=\|v\|_{\mathbb{L}^2}=
\|\reps{X}(t)\|_{\mathbb{L}^2}\qquad  \forall t \in[0,T],  \forall
\epsilon>0.
\]
Thus,\vspace*{2pt} we deduce the strong convergence of $\reps{X}(t)$ to $\widehat
{X}^R(t)$ in $\mathbb{L}^2$, a.s. for each $t \in[0,T]$. Since $\reps
{X}$ converges to $\widehat{X}^R$ in $L_w^{\infty}(0,T; \mathbb
{H}^2)$, we get using Lemma~\ref{boundinH2} that
\[
\|\widehat{X}^R\|_{L^{\infty}(0,T; \mathbb{H}^2
)}\leq
\liminf_{\epsilon\to0} \|\reps{X}\|_{L^{\infty}(0,T;
\mathbb{H}^2)} \leq C(R,T),  \qquad\mathbb{P}\mbox{-a.s.}
\]
Interpolating $\mathbb{H}^1$ between $\mathbb{L}^2$ and $\mathbb{H}^2$,
we conclude that
%
%e4.30 #&#
\begin{equation}\label{cvpointwise}
\lim _{\epsilon\to0}\|\reps{X}(t)-\widehat{X}^R(t)\|
_{\mathbb{H}^1} = 0\qquad  \forall t \in[0,T], \qquad \mathbb{P}\mbox{-a.s.},
\end{equation}
and $\widehat{X}^R \in C([0,T]; \mathbb{H}^1)$;
it follows that, almost surely for all $t$ in $[0,T]$, $\gamma
^R(t)=\|\widehat{X}^R(t)\|_{\mathbb{H}^1}^2$ and $\widehat{X}^R$
is a
solution of (\ref{manakovlimitetronquee}). Thus, the limit in law of
$\reps{X}$ is unique and is given by the solution $X^R$ of
(\ref{manakovlimitetronquee}).

The final step consists in recovering the convergence in law in $C
([0,T], \mathbb{H}^1)$.
Since $\reps{Y}$ is uniformly bounded in $\epsilon$ in $C^{\alpha
}
([0,T], \mathbb{H}^{-1})\cap C([0,T];
\mathbb
{H}^2) $ with $0\leq\alpha<1/2$, we deduce that it is a.s.
uniformly bounded
in $\epsilon$ in $C^{\beta}([0,T], \mathbb
{H}^1)$
with $\beta=\alpha/3$. Moreover, using pointwise convergence (\ref
{cvpointwise}), expression (\ref{processy}) and uniform bounds (\ref
{boundinH2}), we get pointwise convergence in $\mathbb{H}^1$ of $\reps
{Y}$ to $X^R$. We conclude that $\reps{Y}$ converges in law to $X^R$ in
$C([0,T], \mathbb{H}^1(\mathbb{R}))$ and by Lemma
\ref{procheproba}, the convergence in law of $\reps{X}$ to $X^R$ in
$C([0,T], \mathbb{H}^1(\mathbb{R}))$ follows.
\end{pf}
%
%re4.4 #&#
\begin{rmq}
Using the Arzela--Ascoli and Banach--Alaoglu theorems, Lemma \ref
{lemmatension} and the Tychonov theorem, we deduce that $( \mathcal
{L}( \reps{X}) )_{R \in\mathbb{N}} $ is tight
on~$\mathcal{K}^{\mathbb{N}}$. Thus, the same arguments as above lead
to the convergence in law of $( \reps{X} )_{R \in\mathbb{N}}$ to $
(X^R)_{R \in\mathbb{N}}$ (see~\cite{bouardschema}).
\end{rmq}

%s4.5 #&#
\subsection{\texorpdfstring{Convergence of $(\eps{X})_{\epsilon>0}$ to $X$}
{Convergence of (X epsilon)epsilon>0 to X}}

Using the Skorokhod theorem, we can construct new random variables $
\reps{\widetilde{X}}$, $\widetilde{X}^R$ on a common probability space
$(\widetilde{\Omega}, \widetilde{\mathcal{F}}, \widetilde
{\mathcal
{F}_t},\widetilde{\mathbb{P}})$ and with values in $C
(
[0,T], \mathbb{H}^1)$ such that for any $R>0$,
\[
\cases{
\displaystyle \reps{\widetilde{\mu}}=\reps{\mu},\vspace*{2pt}\cr
\widetilde{\mu}^R=\mu^R,}
\quad\mbox{and}\quad \reps{\widetilde{X}}
\mathop{\rightarrow}_{\epsilon\to0}
\widetilde{X}^R ,\qquad \widetilde{\mathbb{P}} \mbox{-a.s. in }
C([0,T], \mathbb{H}^1).
\]
We define the escape times $\widetilde{\tau}^R$ and $\reps
{\widetilde
{\tau}}$ associated to the cutoff:
\[
\widetilde{\tau}^R=\inf\{t \in[0,T], \|\widetilde{X}^R(t)\|
_{\mathbb{H}^1}>R \}
\]
and
\[
\reps{\widetilde{\tau}}=\inf\{t \in[0,T], \|\eps
{\widetilde{X}^R}(t)\|_{\mathbb{H}^1}>R \}.
\]
Let $\eps{\widetilde{X}}$ and $\widetilde{X}$ be the processes, with
values in $\mathcal{E}(\mathbb{H}^1)$, defined,
respectively, by $\eps{\widetilde{X}}(t)=\reps{\widetilde{X}}(t)$ for
$t<\reps{\widetilde{\tau}}$ and $\widetilde{X}(t)=\widetilde{X}^R(t)$
for $t< \widetilde{\tau}^R$, $\widetilde{X}(t)=\Delta$ for $t
\geq
\tau^*= \lim_{R \to+\infty} \widetilde{\tau}^R$. Then if $\tau<
\tau
^*$ a.s. is a stopping time, the process $\eps{\widetilde{X}}$ converges
to $\widetilde{X}$ a.s. in $C([0,\tau], \mathbb{H}^1
(\mathbb
{R}))$. Hence, the convergence in law in $\mathcal
{E}
(\mathbb{H}^1)$ follows.
\section{\texorpdfstring{Study of the driving process $\nu$}
{Study of the driving process nu}} \label{A1}

We recall in this appendix some results obtained in~\cite{marty,2}
about the driving process $\nu$.
%
%pr5.1 #&#
\begin{prop}\label{prop2}
The process $\nu=(\nu_1, \nu_2)^t$ is a Feller process that evolves on
the unit sphere $\mathbb{S}^3$ of $\mathbb{C}^2 \sim\mathbb{R}^4$.
Furthermore, it admits a unique invariant measure\vadjust{\goodbreak} $\Lambda$, which is
the uniform measure on $\mathbb{S}^3$, under which it is ergodic. For
all $f \in C^2_b(\mathbb{S}^3)$ satisfying the Fredholm
alternative (or null mass condition) $\mathbb{E}_{\Lambda}
(f(\nu
))=\int_{\mathbb{S}^3}f(y)\Lambda(dy)=0$, the Poisson equation
$\mathscr{L}_{\nu}u(y)+f(y)=0$
admits a unique solution of class $C^2_b(\mathbb{S}^3)$, up to a
constant, which can be written as $u(y)=\int_0^{+\infty}\mathbb
{E}
[f(\nu(t))|\nu_0=y]\,dt$.
\end{prop}

Let us recall that $\sigb(\nu(t)) =\sigma_1m_1+\sigma
_2m_2+\sigma_3m_3$ where $m_j(t)= g_j(\nu(t)) $. We now
state a result related to the effect of the random PMD on the pulse evolution.
%
%co5.1 #&#
\begin{cor}\label{prop3}
\textup{(1)} The process $m=(m_1,m_2,m_3) \in\mathbb{S}^3$ is a Feller
process with a unique invariant measure $\Lambda\circ g^{-1}$ under
which it is ergodic.\vspace*{-4pt}
%We suppose fullfiled the null mass condition i.e. $
%)=0$, then the Poisson equation $\mathscr{L}_mu(y)+f(y)=0$
%admits a unique solution, up to a constant, of class $C_b^2(
%[f( m(t)) | m(0)=y ]$.
%
\begin{longlist}[(2)]
\item[(2)] For $j=1,2,3$,
$\mathbb{E}_{\Lambda}(g_j(\nu))=\mathbb{E}_{\Lambda
\circ
g_j^{-1}}(m)=0$ and $
\mathbb{E}_{\Lambda}(g_j(\nu(t))g_k(\nu(t)))= \delta
_{jk} / 3$.
As a consequence,
\begin{eqnarray*}
\mathbb{E}_{\Lambda}(N_{1, \nu}(X ))&=&\tfrac
{2}{3}(2|X_2|^2-|X_1|^2
)X_1,\\
\mathbb{E}_{\Lambda}(N_{2, \nu}(X ))&=&\tfrac
{2}{3}(2|X_1|^2-|X_2|^2)X_2.
\end{eqnarray*}

\item[(3)] For $j,k=1,2,3$,
\[
\int_0^{+\infty}\mathbb{E}_{\Lambda}[g_j( \nu(0))
g_k( \nu(t)) ]\,dt=
\cases{
\displaystyle \frac{1}{12\gamma_c}, &\quad if $j=k$,\vspace*{2pt}\cr
0, &\quad if $j\neq k$,}
\]
\end{longlist}
where $\gamma_c$ is the constant appearing in (\ref{nu1}).
\end{cor}

% Then in the asymptotic dynamic of the electric field, the nonlinear
%PMD effect is averaged to zero.% and has no impact on the pulse
%propagation for the scaling we considered.
% \begin{pf}[Sketch of the proof of Corollary~\ref{prop3}]
% The proof of $1.$ follows from Proposition~\ref{prop2} and the change
%of variable $m=g(\nu)$.
%
% To prove $2.$, write the process $\nu$ in spherical coordinates in $
%has a density function given by $\frac{1}{2\pi^2}\sin^2(\theta_1)\sin(
%moments is just a calculation using trigonometric identities. The
%proof of Corollary~\ref{prop4} follows with the same arguments.
% \end{proof}

%s6 #&#
\section{Proof of technical lemmas}\label{A2}

\mbox{}

\begin{pf*}{Proof of Lemma~\ref{firstcv}} \label{prooffirstcv}
Let $v$ be in $\mathbb{H}^3$. Using the explicit representation (\ref
{solutionpoisson1}) of $f^1$, we obtain, since $D_vf(v) \in\mathbb
{H}^1(\mathbb{R})$, that
\begin{eqnarray*}
|f^1(v,y)|&=&\biggl|\biggl\langle D_vf(v),b'\widetilde{\sigb}(y
)\,\frac{\partial v}{\partial x} \biggr\rangle\biggr|\\
&\leq& b' \|
D_vf(v)\|_{\mathbb
{H}^{1}}\biggl\|\frac{\partial v}{\partial x}\biggr\|_{\mathbb{H}^{-1}}
\sum_{j=1}^3 \biggl|\int_0^{+\infty}
\mathbb{E}\bigl( g_j(\nu(t)) | \nu
(0)=y\bigr) \,dt\biggr|.
\end{eqnarray*}
Moreover,\vspace*{1pt} by Proposition~\ref{prop2} the integral
$\int_0^{+\infty} \mathbb{E}( g_j(\nu(t)) | \nu (0)=y) \,dt$ converges
because $g_j$ is a bounded function of $\nu\in\mathbb {S}^3$. Since $v
\mapsto D_vf(v)$ is a continuous function which is bounded on bounded
sets of $\mathbb{H}^{-1}$, we deduce that
\[
\mathop{\sup_{v \in\mathcal{B}(K)}}_{y \in\mathbb{S}^3}
|\eps{f}^1(v,y)|\leq b'C(K).
\]
The function $f^2$ given by (\ref{solutionpoisson2}) may be bounded
using the same arguments. Indeed,
\[
\langle D_vf(v),i\Theta_R(\|v\|_{\mathbb{H}^1}^2
)\widetilde{F}(v,y) \rangle
\leq\|D_vf(v)\|_{\mathbb{H}^{1}}\|\widetilde{F}
(v,y)\|_{\mathbb{H}^{-1}}.\vadjust{\goodbreak}
\]
Since for all $v \in\mathbb{H}^3$, $y \mapsto F_{y}(v
)-F(v)$ is a function of class $C^2_b$ on $\mathbb{S}^3$,
with values in $\mathbb{H}^{-1}$, satisfying the null mass condition of
Proposition~\ref{prop2}, the term $\widetilde{F}(v,y)$ is
bounded. Moreover, $v \mapsto F_{y}(v)-F(v)$ is
bounded in $\mathbb{H}^{-1}$ on bounded sets of $\mathbb{H}^1$ by the
continuous embeddings $\mathbb{H}^{1}(\mathbb{R})
\hookrightarrow\mathbb{L}^4(\mathbb{R})$ and $ \mathbb
{L}^{4/3}(\mathbb{R})\hookrightarrow\mathbb
{H}^{-1}
(\mathbb{R})$. In addition,
\begin{eqnarray*}
&&
\Biggl|(b')^2 \sum_{k,l=1}^3\biggl\langle D_v^2f(v)\sigma
_k\,\frac{\partial v}{\partial x},\sigma_l\,\frac{\partial
v}{\partial x} \biggr\rangle
\wwidetilde{g}_{k,l}(y)+\biggl\langle D_vf(v),( b')^2 \widetilde
{\widetilde{\sigb}}(y) \,\frac{\partial^{2} v}{\partial x^{2}}
\biggr\rangle\Biggr|\\
&&\qquad\leq C\sum
_{k,l=1}^3
\biggl(\biggl|\biggl\langle D^2_{v}f(v)\sigma_k\,\frac{\partial v}{\partial
x},\sigma_l\,\frac{\partial v}{\partial x}
\biggr\rangle\biggr|+\biggl|\biggl\langle D_vf(v),\sigma_k\sigma_l\,\frac{\partial^{2}
v}{\partial x^{2}} \biggr\rangle
\biggr|\biggr)\\
&&\qquad\leq C \biggl( \|D^2_{v}f(v)\|_{\mathscr{L}( \mathbb
{H}^{-1}, \mathbb{H}^1)}\biggl\|\frac{\partial v}{\partial
x}\biggr\|_{\mathbb{H}^{-1}}^2 +
\|D_vf(v)\|_{\mathbb{H}^{1}} \biggl\|\frac{\partial^{2} v}{\partial
x^{2}}\biggr\|_{\mathbb
{H}^{-1}}\biggr).
\end{eqnarray*}
Since $v \mapsto D_{v}f(v)$ and $v \mapsto D^2_{v}f(v)$ are bounded on
bounded sets of $\mathbb{H}^{-1}(\mathbb{R})$, we conclude
the proof of the lemma.
\end{pf*}
\begin{pf*}{Proof of Lemma~\ref{secondcv}}\label{proofsecondcv}
Replacing $f^1$ by its expression (\ref{solutionpoisson1}), we get
\begin{eqnarray*}
&&\biggl\langle  D_vf^1(v,y),\frac{id_0}{2}\,\frac{\partial^{2} v}{\partial
x^{2}}+i \Theta_R
(\|v\|_{\mathbb{H}^1}^2)F_{y}{( v) } \biggr\rangle \\
&&\qquad= -\biggl\langle D^2_{v}f(v)b'\widetilde{\sigb}(y) \,\frac
{\partial v}{\partial x},\frac{id_0}{2}\,\frac{\partial^{2}
v}{\partial x^{2}}+i \Theta_R(\|v\|
_{\mathbb{H}^1}^2) F_{y}{( v) } \biggr\rangle\\
&&\qquad\quad{}-\biggl\langle D_vf(v),b'\widetilde{\sigb}(y)\frac
{id_0}{2}\,\frac{\partial^{3} v}{\partial x^{3}}+ib'\widetilde{\sigb
}(y) \Theta
_R(\|v\| _{\mathbb{H}^1}^2) \,\partial_xF_{y}{( v) }
\biggr\rangle.
\end{eqnarray*}
By the assumptions on $f$, $v \mapsto D_vf(v)$ and $v \mapsto
D^2_{v}f(v)$ are continuous bounded functions on bounded sets of
$\mathbb{H}^{-1}(\mathbb{R}) $. Moreover, $D^2_{v}f(v)
\in
\mathcal{L}( \mathbb{H}^{-1},\break \mathbb{H}^{1} ) $, $D_{v}f(v)
\in\mathbb{H}^{1}$ and $\frac{\partial^{3} v}{\partial x^{3}} \in
\mathbb{L}^{2}$. Using the
bound (\ref{csteC1}), we deduce that
\[
\mathop{\sup_{v \in\mathcal{B}(K)}}_{y \in\mathbb{S}^3}
\biggl|\biggl\langle D_vf^1(v,y),\frac{id_0}{2}\,\frac{\partial^{2} v}{\partial
x^{2}}+i \Theta_R(\|
v\|_{\mathbb{H}^1}^2)F_{y}{( v) } \biggr\rangle\biggr|\leq C(K
).
\]
Let us now compute the first derivative of $f^2$ using expression
(\ref{solutionpoisson2}); for all $h$ in $\mathbb{H}^{1}$ and $v$ in
$\mathbb{H}^3$,
\begin{eqnarray*}
&&\langle D_vf^2(v,y),h \rangle\\
&&\qquad= \langle D^2_{v}f(v)h,i \Theta_R(\|v\|_{\mathbb
{H}^1}^2)\widetilde{F}( v,y ) \rangle\\
&&\qquad\quad{}+\langle
D_vf(v),2i\Theta_R'(\|v\|_{\mathbb{H}^1}^2)
(v,h)_{\mathbb{H}^1}\widetilde{F}( v,y )
\rangle\\
&&\qquad\quad{}+ \langle D_vf(v),i \Theta_R(\|v\|_{\mathbb{H}^1}^2)
D_v\widetilde{F}( v,y ).h \rangle\\
&&\qquad\quad{}-(b')^2 \sum_{k,l=1}^3D_v^3f(v).\biggl( \sigma_k\,\frac
{\partial v}{\partial x},\sigma_l\,\frac{\partial
v}{\partial x}, h\biggr) \wwidetilde{g}_{k,l}(y)\\
&&\qquad\quad{}-2(b')^2 \sum_{k,l=1}^3\biggl\langle D_v^2f(v)\sigma_k\,\frac
{\partial h}{\partial x},\sigma_l\,\frac{\partial
v}{\partial x} \biggr\rangle\wwidetilde{g}_{k,l}(y)\\
&&\qquad\quad{}-\biggl\langle D^2_vf(v)h,(b')^2 \widetilde{\widetilde
{\sigb}}(y)\,\frac{\partial^{2} v}{\partial x^{2}} \biggr\rangle\\
&&\qquad\quad{}-\biggl\langle D_vf(v),(b')^2 \widetilde{\widetilde{\sigb
}}(y) \,\frac{\partial^{2} h}{\partial x^{2}} \biggr\rangle.
\end{eqnarray*}
%
% \[
% && \dual{D_vf^2(v,y)}{h}\\
% &=&\int_0^{+\infty}\mathbb{E}\{
%(\nl{\nu(t)}{v}-F( v ) ) }| \nu(0)=y
% &+&\int_0^{+\infty}\mathbb{E} \{ \dual{D_vf(v)}{
%)_{\mathbb{H}^1} \Theta_R(\norm{v}{\mathbb{H}^1}^2)
%(\nl{\nu(t)}{v}-F( v )) } | \nu(0)=y
% &+&\int_0^{+\infty}\mathbb{E} \{ \dual{D_vf(v)}{i
%( v )).h } | \nu(0)=y\} dt\\
% &-&\int_0^{+\infty}\mathbb{E}\{ \dual{D^2_{v}f^1(v,
%(\dual{D^2_{v}f^1(v,\nu(t))h }{b'\sigb( \nu(t))
% &-&\int_0^{+\infty}\mathbb{E} \{ \dual{D_vf^1(v,
%)\pds{}{x}{h}}) | \nu(0)=y\} dt.
% \]
Taking, respectively, $h=\frac{id_0}{2}\,\frac{\partial^{2} v}{\partial
x^{2}}+i \Theta_R
(\|v\|_{\mathbb{H}^1}^2) F_{y}{( v) }$ and $h=b'\sigb(y)\,\frac
{\partial v}{\partial x}$,
we conclude
\[
\mathop{\sup_{v \in\mathcal{B}(K)}}_{y \in\mathbb{S}^3}
\biggl|\biggl\langle D_vf^2(v,y),\frac{id_0}{2}\,\frac{\partial^{2} v}{\partial
x^{2}}+i \Theta_R(\|
v\|_{\mathbb{H}^1}^2) F_{y}{( v) } \biggr\rangle\biggr|
\leq C(K)
\]
and
\[
\mathop{\sup_{v \in\mathcal{B}(K)}}_{y \in\mathbb{S}^3}
\biggl|\biggl\langle D_vf^2(v,y),b'\sigb(y)\,\frac{\partial v}{\partial x}
\biggr\rangle\biggr|
\leq C(K),
\]
since $v \mapsto D^3_{v}f(v)$ is bounded on the bounded set of $\mathbb
{H}^{-1}(\mathbb{R}) $ with values in $\mathcal
{L}_3
(\mathbb{H}^{-1}, \mathbb{R}) $ and $\frac{\partial^{4}
v}{\partial x^{4}} \in
\mathbb{H}^{-1}$.
\end{pf*}
\begin{pf*}{Proof of Lemma~\ref{lemmakolmogorov}}\label{prooflemmakolmogorov}
Let us recall that the family $\{e_i\}_{i \in\mathbb{N}^*}$
denotes a complete orthonormal system of $\mathbb{H}^1$ % with values
%in $\mathbb{H}^2$,
constructed from a complete orthonormal system $\{\widetilde
{e}_i\}_{i \in\mathbb{N}^*}$ in $\mathbb{L}^2$ and $\langle
\cdot,\cdot\rangle$
is the duality product between $\mathbb{H}^{1}$--$\mathbb{H}^{-1}$. Then
\[
\|\reps{Y}(t)-\reps{Y}(s)\|_{\mathbb{H}^{-1}}^4= \Biggl\{
\sum
_{i=1}^{+\infty}\langle e_i,\reps{Y}(t)-\reps{Y}(s) \rangle^2
\Biggr\}^2.
\]
Using twice the Young inequality and the expression of $\reps{Y}$ given
by (\ref{secondformulationforY}) and (\ref{generateurtension}), we obtain
\begin{eqnarray*}
& & \|\reps{Y}(t)-\reps{Y}(s)\|_{\mathbb{H}^{-1}}^4 \\
&&\qquad\leq C\biggl\|\frac{d_0}{2}\int_s^t \partial_x^2{\reps{X}(t')}
\,dt'\biggr\|_{\mathbb{H}^{-1}}^4\\
&&\qquad\quad{}+C\biggl\|\int_s^t \Theta_R(\|\reps
{X}(t')\|_{\mathbb{H}^1}^2) F_{\eps{\nu}(t')}{( \reps
{X}(t')) }\,dt'\biggr\|_{\mathbb{H}^{-1}}^4
\\
&&\qquad\quad{}
+C\biggl\|\int_s^t (b')^2\widetilde{\sigb} (\eps
{\nu}(t'))\sigb(\eps{\nu}(t'))\,\partial
_x^2{\reps{X}(t')}\,dt'\biggr\|_{\mathbb{H}^{-1}}^4\\
&&\qquad\quad{}+C\Biggl[ \sum
_{i=1}^{+\infty}\langle e_i,\reps{M}(t)-\reps{M}(s) \rangle^2
\Biggr]^2 \\
&&\qquad\quad{}
+C\epsilon^4 \biggl\|\int_s^tb'\widetilde{\sigb} (\eps{\nu
}(t'))\,\partial_x\biggl(\frac{d_0}{2}\,\partial_x^2{\reps
{X}(t')}\\
&&\qquad\quad\hspace*{116.5pt}{} + \Theta_R(\|\eps{X}(t')\|_{\mathbb{H}^1}^2)
F_{\eps{\nu}(t')}{( \reps{X}(t')) }\biggr) \,dt'\biggr\|_{\mathbb{H}^{-1}}^4.
\end{eqnarray*}
We bound each term separately. Using Lemma~\ref{boundinH2},
\[
\biggl\|\int_s^t\frac{d_0}{2}\,\frac{\partial^{2} \reps{X}(t')}{\partial
x^{2}}\,dt'\biggr\|_{\mathbb
{H}^{-1}}^4\leq C(R,T)(t-s)^4.
\]

Using that $F$ is cubic and Lemma~\ref{boundinH2},% and the following
%sobolev embedding $\mathbb{H}^1 \hookrightarrow\mathbb{L}^4$ and $
%
\[
\biggl\|\int_s^t\Theta_R(\|\eps{X}(t')\|_{\mathbb{H}^1}^2
)F_{\eps{\nu}(t')}{( \reps{X}(t')) }\,dt'\biggr\|_{\mathbb{H}^{-1}}^4
\leq C(R,T)(t-s)^4
\]
and using Lemma~\ref{boundinH2} and the bound (\ref{csteC1}),
\[
\biggl\|\int_s^t(b')^2\widetilde{\sigb} (\eps{\nu
}(t'))\sigb(\eps{\nu}(t'))\,\frac{\partial^{2}
\reps{X}(t')}{\partial x^{2}}\,dt'\biggr\|_{\mathbb{H}^{-1}}^4
\leq C(R,T)(t-s)^4.
\]
Finally, we bound the $\epsilon^4$ term that is well defined because
$\reps{X}$ has values in $\mathbb{H}^3$. Using the Cauchy--Schwarz
inequality, Lemma~\ref{boundinH2} and (\ref{csteC1}), we get for all
$\epsilon<1$,
\begin{eqnarray*}
&&
\epsilon^4 \biggl\|\int_s^tb'\widetilde{\sigb} (\eps
{\nu}(t'))\,\frac{\partial}{\partial x}\biggl( \frac
{d_0}{2}\,\frac{\partial^{2} \reps{X}(t')}{\partial x^{2}} + \Theta
_R(\|\eps{X}\|_{\mathbb
{H}^1}^2) F_{\eps{\nu}(t')}{( \reps{X}(t')) }\biggr) \,dt'\biggr\|
_{\mathbb{H}^{-1}}^4\\
&&\qquad\leq \epsilon^4 (b')^4M^4 \biggl\|\int_s^t \frac{d_0}{2}\,\frac
{\partial^{3} \reps{X}(t')}{\partial x^{3}} + \frac{\partial
}{\partial x}\bigl( \Theta_R(\|\eps
{X}\|_{\mathbb{H}^1}^2) F_{\eps{\nu}(t')}{( \reps {X}(t'))
}\bigr) \,dt'\biggr\|_{\mathbb{H}^{-1}}^4\\
&&\qquad\leq C(R,T)(t-s)^4.
\end{eqnarray*}
Taking the expectation and adding the previous estimates, we deduce that
\begin{eqnarray*}
& & \mathbb{E}\bigl(\|\reps{Y}(t)-\reps{Y}(s)\|_{\mathbb
{H}^{-1}}^4\bigr)\\
&&\qquad\leq C(R,T)(t-s)^4+C\mathbb{E}\biggl(\biggl[ \sum
_{j\in
\mathbb{N}^*}\langle e_j,\reps{M}(t)-\reps{M}(s) \rangle^2
\biggr]^2
\biggr).
\end{eqnarray*}
In order to prove a uniform bound, with respect to $\epsilon$, of the
second term, we will use the Burkholder--Davis--Gundy\vspace*{1pt} inequality and,
consequently, we have to compute the quadratic variation $\llangle\reps
{M}(t) \ggangle$ of $\reps{M}(t)$ defined, for all $j\in\mathbb{N}^*$, by
\[
\langle e_j,\reps{M}(t) \rangle=f_{j,\epsilon}(\reps{X}(t),
\eps
{\nu
}(t))-f_{j,\epsilon}(v,y)-\int_0^t\reps
{\mathscr
{L}}f_{j,\epsilon}(\reps{X}(s),\eps{\nu}(s))\,ds,
\]
where\vspace*{1pt} $\reps{\mathscr{L}}f_{j,\epsilon}(\reps{X}(s),\eps
{\nu
}(s))$ is given by (\ref{generateurtension}). The next lemma
states that the process $\llangle\reps{M}(t) \ggangle$ can be
expressed only in terms of the infinitesimal generator $\mathscr
{L}_{\nu
}$ of the Markov process $\nu$.
%&=&\int_0^t\mathscr{L}_{m} (\dual{\varphi^1}{e_i}\dual{
%
%le6.1 #&#
\begin{lemma}\label{generateurXnu}
For all $j$ in $\mathbb{N}^*$
\begin{eqnarray*}
&&
\langle e_j,\llangle\reps{M}(t) \ggangle e_j \rangle\\
&&\qquad=
( b')^2 \int_0^t \mathscr{L}_{\nu}\biggl(\biggl\langle
e_j,\widetilde{\sigb} (\eps{\nu}(s)) \,\frac{\partial
\reps{X}}{\partial x}(s) \biggr\rangle^2\biggr)\,ds\\
&&\qquad\quad{}-2( b')^2 \int_0^t\biggl\langle e_j,\widetilde{\sigb}
(\eps{\nu}(s)) \,\frac{\partial\reps{X}}{\partial
x}(s) \biggr\rangle\biggl\langle
e_j,\mathscr{L}_{\nu}\widetilde{\sigb} (\eps{\nu}(s)
) \,\frac{\partial\reps{X}}{\partial x}(s) \biggr\rangle \,ds.
\end{eqnarray*}
\end{lemma}

Thus, using the Burkholder--Davis--Gundy inequality,
\[
\mathbb{E}\Biggl(\Biggl( \sum_{j=1}^{+\infty}\langle e_j,\reps
{M}(t)-\reps{M}(s) \rangle^2\Biggr)^2 \Biggr)\leq C(R,T)
|t-s|^2,
\]
thanks to Lemma~\ref{boundinH2} and Proposition~\ref{prop2}. Adding the
previous estimates,
\[
\mathbb{E}\bigl(\|\reps{Y}(t)-\reps{Y}(s)\|_{\mathbb
{H}^{-1}}^4
\bigr)\leq C(R,T)|t-s|^2,
\]
and Lemma~\ref{lemmakolmogorov} is proved.
\end{pf*}
\begin{pf*}{Proof of Lemma~\ref{generateurXnu}}
A classical computation shows that for all \mbox{$j \in\mathbb{N}^*$},
\begin{eqnarray*}
\langle e_j,\llangle\reps{M}(t) \ggangle e_j \rangle
&=&
\int_0^t\reps{\mathscr{L}}( f_{j,\epsilon}(\reps
{X}(s), \eps
{\nu}(s))^2)\,ds\\
&&{} -2 \int_0^tf_{j,\epsilon}
(\reps
{X}(s), \eps{\nu}(s)) \reps{\mathscr{L}} f_{j,\epsilon
}
(\reps{X}(s), \eps{\nu}(s))\,ds.
\end{eqnarray*}
Now, for all $j$ in $\mathbb{N}^*$,
\begin{eqnarray*}
( f_{j,\epsilon}(\reps{X}(s), \eps{\nu}(s)
))^2
&=&\langle e_j,\reps{X}(s) \rangle^2-2b'\epsilon\langle e_j,\reps
{X}(s) \rangle\biggl\langle e_j,\widetilde{\sigb} (\eps{\nu
}(s)) \,\frac{\partial\reps{X}}{\partial x}(s) \biggr\rangle\\
&&{}+(b'
)^2 \epsilon^2\biggl\langle e_j,\widetilde{\sigb} (\eps{\nu
}(s)) \,\frac{\partial\reps{X}}{\partial x}(s) \biggr\rangle^2.
\end{eqnarray*}
Thus, we get
\begin{eqnarray*}
\hspace*{-4pt}& & \reps{\mathscr{L}} ( f_{j,\epsilon}
(\reps
{X}(s), \eps{\nu}(s)) )^2 \\[-2pt]
\hspace*{-4pt}&&\qquad= 2\langle e_j,\reps
{X}(s) \rangle\biggl\langle e_j,\frac{id_0}{2}\,\frac{\partial^{2} \reps
{X}(s)}{\partial x^{2}}+i\Theta_R(\|\reps{X}(s)\|_{\mathbb
{H}^1}^2
)F_{\eps{\nu}(s)}(\reps{X}(s)) \biggr\rangle\\[-2pt]
\hspace*{-4pt}&&\qquad\quad{}-2b'\epsilon\biggl\langle e_j,\frac{id_0}{2}\,\frac{\partial^{2} \reps
{X}(s)}{\partial x^{2}}+i\Theta_R(\|\reps{X}(s)\|_{\mathbb
{H}^1}^2
)F_{\eps{\nu}(s)}(\reps{X}(s)) \biggr\rangle\\[-2pt]
\hspace*{-4pt}&&\qquad\quad\hspace*{9.5pt}{}\times\biggl\langle
e_j,\widetilde{\sigb} (\eps{\nu}(s)) \,\frac{\partial
\reps{X}}{\partial x}(s) \biggr\rangle\\[-2pt]
\hspace*{-4pt}&&\qquad\quad{}-2b'\epsilon\langle e_j,\reps{X}(s) \rangle\\[-2pt]
\hspace*{-4pt}&&\qquad\quad\hspace*{9.5pt}{}\times\biggl\langle e_j,\widetilde
{\sigb} (\eps{\nu}(s)) \,\frac{\partial}{\partial
x}\biggl( \frac
{id_0}{2}\,\frac{\partial^{2} \reps{X}(s)}{\partial x^{2}}+i\Theta
_R(\|\reps{X}(s)\|
_{\mathbb{H}^1}^2)F_{\eps{\nu}(s)}(\reps{X}(s)
)\biggr) \biggr\rangle\\[-2pt]
\hspace*{-4pt}&&\qquad\quad{} -2(b')^2 \biggl\langle e_j,\mathscr{L}_{\nu}\widetilde
{\sigb} (\eps{\nu}(s)) \,\frac{\partial\reps
{X}}{\partial x}(s) \biggr\rangle
\biggl\langle e_j,\widetilde{\sigb} (\eps{\nu}(s)) \,\frac
{\partial\reps{X}}{\partial x}(s) \biggr\rangle\\[-2pt]
\hspace*{-4pt}&&\qquad\quad{} +(b')^2\mathscr{L}_{\nu}\biggl(\biggl\langle
e_j,\widetilde{\sigb} (\eps{\nu}(s)) \,\frac{\partial
\reps{X}}{\partial x}(s) \biggr\rangle^2\biggr)\\[-2pt]
\hspace*{-4pt}&&\qquad\quad{}
+2b'\biggl\langle e_j,\widetilde{\sigb} (\eps{\nu}(s)) \,\frac
{\partial}{\partial x}\biggl( b'\sigb(\eps{\nu
}(s))\,\frac{\partial\reps{X}(s)}{\partial x}\biggr)
\biggr\rangle\langle e_j,\reps{X}(s) \rangle\\[-2pt]
\hspace*{-4pt}&&\qquad\quad{}
+ 2(b')^2 \epsilon^2\biggl\langle e_j,\widetilde{\sigb}
(\eps{\nu}(s)) \,\frac{\partial\reps{X}}{\partial
x}(s) \biggr\rangle\biggl\langle
e_j,\widetilde{\sigb} (\eps{\nu}(s))\frac
{id_0}{2}\,\frac{\partial^{3} \reps{X}(s)}{\partial x^{3}} \biggr\rangle\\[-2pt]
\hspace*{-4pt}&&\qquad\quad{}
+ 2(b')^2\epsilon^2\biggl\langle e_j,\widetilde{\sigb}
(\eps{\nu}(s)) \,\frac{\partial\reps{X}}{\partial
x}(s) \biggr\rangle\\[-2pt]
\hspace*{-4pt}&&\qquad\quad\hspace*{9.5pt}{}\times\biggl\langle
e_j,i\widetilde{\sigb} (\eps{\nu}(s)) \Theta_R
(\|\reps{X}(s)\| _{\mathbb{H}^1}^2) \,\frac{\partial
}{\partial x}F_{\eps{\nu
}(s)}(\reps{X}(s)) \biggr\rangle\\[-2pt]
\hspace*{-4pt}&&\qquad\quad{}
-2(b')^2\epsilon\biggl\langle e_j,\widetilde{\sigb}
(\eps{\nu}(s)) \,\frac{\partial\reps{X}}{\partial
x}(s) \biggr\rangle\biggl\langle
e_j,b'\widetilde{\sigb} (\eps{\nu}(s))\sigb
(\eps{\nu}(s))\,\frac{\partial^{2} \reps{X}(s)}{\partial
x^{2}} \biggr\rangle.
\end{eqnarray*}
The same kinds of computations for the term $2f_{j,\epsilon} \reps
{\mathscr{L}}f_{j,\epsilon}$ lead to the result.
\end{pf*}

\section*{Acknowledgments}

The authors wish to thank the referees for fruitful comments.

%suskaldyti doi

% imsref loaded by lrinkeviciute, 2012-05-04 13:52:37
% imsref loaded by lrinkeviciute, 2012-05-04 14:32:20

\printaddresses


\begin{thebibliography}{33}
% BibTex style file: ims.bst, 2011-05-30
% Default style options (sort=0,type=number).
% Used options (sort=1,type=number).

%b1 ###
\bibitem{Agrawalappli}
\begin{bbook}[auto:STB|2012/04/30|08:06:40]
\bauthor{\bsnm{Agrawal},~\bfnm{G.~P.}\binits{G.~P.}}
(\byear{2001}).
\btitle{Applications of Nonlinear Fiber Optics}.
\bpublisher{Academic Press}, \baddress{San Diego}.
\bptok{imsref}%
\end{bbook}
\endbibitem

%b2 ###
\bibitem{Agrawal}
\begin{bbook}[auto:STB|2012/04/30|08:06:40]
\bauthor{\bsnm{Agrawal},~\bfnm{G.~P.}\binits{G.~P.}}
(\byear{2001}).
\btitle{Nonlinear Fiber Optics}, \bedition{3rd} ed.
\bpublisher{Academic Press}, \baddress{San Diego}.
\bptok{imsref}%
\end{bbook}
\endbibitem

%b3 ###
\bibitem{azencott}
\begin{bincollection}[mr]
\bauthor{\bsnm{Azencott},~\bfnm{R.}\binits{R.}}
(\byear{1980}).
\btitle{Grandes d\'eviations et applications}.
In \bbooktitle{Eighth {S}aint {F}lour {P}robability {S}ummer {S}chool---1978
  ({S}aint {F}lour, 1978)}.
\bseries{Lecture Notes in Math.}
\bvolume{774}
\bpages{1--176}.
\bpublisher{Springer}, \baddress{Berlin}.
\bid{mr={0590626}}
\bptok{imsref}%
\end{bincollection}
\endbibitem

%b4 ###
\bibitem{billingsley}
\begin{bbook}[mr]
\bauthor{\bsnm{Billingsley},~\bfnm{Patrick}\binits{P.}}
(\byear{1968}).
\btitle{Convergence of Probability Measures}.
\bpublisher{Wiley}, \baddress{New York}.
\bid{mr={0233396}}
\bptok{imsref}%
\end{bbook}
\endbibitem

%b5 ###
\bibitem{papa}
\begin{barticle}[mr]
\bauthor{\bsnm{Blankenship},~\bfnm{G.}\binits{G.}} \AND
  \bauthor{\bsnm{Papanicolaou},~\bfnm{G.~C.}\binits{G.~C.}}
(\byear{1978}).
\btitle{Stability and control of stochastic systems with wide-band noise
  disturbances. {I}}.
\bjournal{SIAM J. Appl. Math.}
\bvolume{34}
\bpages{437--476}.
\bid{issn={0036-1399}, mr={0476129}}
\bptok{imsref}%
\end{barticle}
\endbibitem

%b6 ###
\bibitem{brzezniak}
\begin{barticle}[mr]
\bauthor{\bsnm{Brze{\'z}niak},~\bfnm{Zdzis{\l}aw}\binits{Z.}}
(\byear{1995}).
\btitle{Stochastic partial differential equations in {M}-type {$2$} {B}anach
  spaces}.
\bjournal{Potential Anal.}
\bvolume{4}
\bpages{1--45}.
\bid{doi={10.1007/BF01048965}, issn={0926-2601}, mr={1313905}}
\bptok{imsref}%
\end{barticle}
\endbibitem

%b7 ###
\bibitem{caz}
\begin{bbook}[mr]
\bauthor{\bsnm{Cazenave},~\bfnm{Thierry}\binits{T.}}
(\byear{2003}).
\btitle{Semilinear {S}chr\"odinger Equations}.
\bseries{Courant Lecture Notes in Mathematics}
\bvolume{10}.
\bpublisher{Amer. Math. Soc.},
  \baddress{Providence, RI}.
\bid{mr={2002047}}
\bptok{imsref}%
\end{bbook}
\endbibitem

%b8 ###
\bibitem{DPZ}
\begin{bbook}[mr]
\bauthor{\bsnm{Da~Prato},~\bfnm{Giuseppe}\binits{G.}} \AND
  \bauthor{\bsnm{Zabczyk},~\bfnm{Jerzy}\binits{J.}}
(\byear{1992}).
\btitle{Stochastic Equations in Infinite Dimensions}.
\bseries{Encyclopedia of Mathematics and Its Applications}
\bvolume{44}.
\bpublisher{Cambridge Univ. Press}, \baddress{Cambridge}.
\bid{doi={10.1017/CBO9780511666223}, mr={1207136}}
\bptok{imsref}%
\end{bbook}
\endbibitem

%b9 ###
\bibitem{bouard}
\begin{barticle}[mr]
\bauthor{\bparticle{de} \bsnm{Bouard},~\bfnm{A.}\binits{A.}} \AND
  \bauthor{\bsnm{Debussche},~\bfnm{A.}\binits{A.}}
(\byear{1999}).
\btitle{A stochastic nonlinear {S}chr\"odinger equation with multiplicative
  noise}.
\bjournal{Comm. Math. Phys.}
\bvolume{205}
\bpages{161--181}.
\bid{doi={10.1007/s002200050672}, issn={0010-3616}, mr={1706888}}
\bptok{imsref}%
\end{barticle}
\endbibitem

%b10 ###
\bibitem{bouardDebussche}
\begin{barticle}[mr]
\bauthor{\bparticle{de} \bsnm{Bouard},~\bfnm{A.}\binits{A.}} \AND
  \bauthor{\bsnm{Debussche},~\bfnm{A.}\binits{A.}}
(\byear{2003}).
\btitle{The stochastic nonlinear {S}chr\"odinger equation in {$H\sp 1$}}.
\bjournal{Stoch. Anal. Appl.}
\bvolume{21}
\bpages{97--126}.
\bid{doi={10.1081/SAP-120017534}, issn={0736-2994}, mr={1954077}}
\bptok{imsref}%
\end{barticle}
\endbibitem

%b11 ###
\bibitem{bouardschema}
\begin{barticle}[mr]
\bauthor{\bsnm{De~Bouard},~\bfnm{A.}\binits{A.}} \AND
  \bauthor{\bsnm{Debussche},~\bfnm{A.}\binits{A.}}
(\byear{2004}).
\btitle{A semi-discrete scheme for the stochastic nonlinear {S}chr\"odinger
  equation}.
\bjournal{Numer. Math.}
\bvolume{96}
\bpages{733--770}.
\bid{doi={10.1007/s00211-003-0494-5}, issn={0029-599X}, mr={2036364}}
\bptok{imsref}%
\end{barticle}
\endbibitem

%b12 ###
\bibitem{bD}
\begin{barticle}[mr]
\bauthor{\bparticle{de} \bsnm{Bouard},~\bfnm{Anne}\binits{A.}} \AND
  \bauthor{\bsnm{Debussche},~\bfnm{Arnaud}\binits{A.}}
(\byear{2010}).
\btitle{The nonlinear {S}chr\"odinger equation with white noise dispersion}.
\bjournal{J. Funct. Anal.}
\bvolume{259}
\bpages{1300--1321}.
\bid{doi={10.1016/j.jfa.2010.04.002}, issn={0022-1236}, mr={2652190}}
\bptok{imsref}%
\end{barticle}
\endbibitem

%b13 ###
\bibitem{vovelle}
\begin{bmisc}[auto:STB|2012/04/30|08:06:40]
\bauthor{\bsnm{Debussche},~\bfnm{A.}\binits{A.}} \AND
  \bauthor{\bsnm{Vovelle},~\bfnm{J.}\binits{J.}}
(\byear{2011}).
\bhowpublished{Diffusion limit for a stochastic kinetic problem. Preprint.}
\bptok{imsref}%
\end{bmisc}
\endbibitem

%b14 ###
\bibitem{Doss}
\begin{barticle}[mr]
\bauthor{\bsnm{Doss},~\bfnm{Halim}\binits{H.}}
(\byear{1977}).
\btitle{Liens entre \'equations diff\'erentielles stochastiques et ordinaires}.
\bjournal{Ann. Inst. H. Poincar\'e Sect. B (N.S.)}
\bvolume{13}
\bpages{99--125}.
\bid{mr={0451404}}
\bptok{imsref}%
\end{barticle}
\endbibitem

%b15 ###
\bibitem{ethier}
\begin{bbook}[mr]
\bauthor{\bsnm{Ethier},~\bfnm{Stewart~N.}\binits{S.~N.}} \AND
  \bauthor{\bsnm{Kurtz},~\bfnm{Thomas~G.}\binits{T.~G.}}
(\byear{1986}).
\btitle{Markov Processes: Characterization and Convergence}.
\bpublisher{Wiley}, \baddress{New York}.
\bid{doi={10.1002/9780470316658}, mr={0838085}}
\bptok{imsref}%
\end{bbook}
\endbibitem

%b16 ###
\bibitem{flandoli}
\begin{barticle}[mr]
\bauthor{\bsnm{Flandoli},~\bfnm{Franco}\binits{F.}} \AND
  \bauthor{\bsnm{G{\c{a}}tarek},~\bfnm{Dariusz}\binits{D.}}
(\byear{1995}).
\btitle{Martingale and stationary solutions for stochastic {N}avier--{S}tokes
  equations}.
\bjournal{Probab. Theory Related Fields}
\bvolume{102}
\bpages{367--391}.
\bid{doi={10.1007/BF01192467}, issn={0178-8051}, mr={1339739}}
\bptok{imsref}%
\end{barticle}
\endbibitem

%b17 ###
\bibitem{garnier}
\begin{bbook}[mr]
\bauthor{\bsnm{Fouque},~\bfnm{Jean-Pierre}\binits{J.-P.}},
  \bauthor{\bsnm{Garnier},~\bfnm{Josselin}\binits{J.}},
  \bauthor{\bsnm{Papanicolaou},~\bfnm{George}\binits{G.}} \AND
  \bauthor{\bsnm{S{\o}lna},~\bfnm{Knut}\binits{K.}}
(\byear{2007}).
\btitle{Wave Propagation and Time Reversal in Randomly Layered Media}.
\bseries{Stochastic Modelling and Applied Probability}
\bvolume{56}.
\bpublisher{Springer}, \baddress{New York}.
\bid{mr={2327824}}
\bptok{imsref}%
\end{bbook}
\endbibitem

%b18 ###
\bibitem{marty}
\begin{barticle}[mr]
\bauthor{\bsnm{Garnier},~\bfnm{Josselin}\binits{J.}} \AND
  \bauthor{\bsnm{Marty},~\bfnm{Renaud}\binits{R.}}
(\byear{2006}).
\btitle{Effective pulse dynamics in optical fibers with polarization mode
  dispersion}.
\bjournal{Wave Motion}
\bvolume{43}
\bpages{544--560}.
\bid{doi={10.1016/j.wavemoti.2006.05.001}, issn={0165-2125}, mr={2252753}}
\bptok{imsref}%
\end{barticle}
\endbibitem

%b19 ###
\bibitem{ginibre}
\begin{bmisc}[auto:STB|2012/04/30|08:06:40]
\bauthor{\bsnm{Ginibre},~\bfnm{J.}\binits{J.}}
(\byear{1994/95}).
\bhowpublished{Introduction aux \'equations de Schr\"odinger non
lin\'eaires. Cours de DEA.}
\bptok{imsref}%
\end{bmisc}
\endbibitem

%b20 ###
\bibitem{kushner}
\begin{bbook}[mr]
\bauthor{\bsnm{Kushner},~\bfnm{Harold~J.}\binits{H.~J.}}
(\byear{1984}).
\btitle{Approximation and Weak Convergence Methods for Random Processes, with
  Applications to Stochastic Systems Theory}.
\bseries{MIT Press Series in Signal Processing, Optimization, and Control}
\bvolume{6}.
\bpublisher{MIT Press}, \baddress{Cambridge, MA}.
\bid{mr={0741469}}
\bptok{imsref}%
\end{bbook}
\endbibitem

%b21 ###
\bibitem{menyukappli}
\begin{barticle}[auto:STB|2012/04/30|08:06:40]
\bauthor{\bsnm{Marcuse},~\bfnm{D.}\binits{D.}},
  \bauthor{\bsnm{Wai},~\bfnm{P.~K.~A.}\binits{P.~K.~A.}} \AND
  \bauthor{\bsnm{Menyuk},~\bfnm{C.~R.}\binits{C.~R.}}
(\byear{1997}).
\btitle{Application of the Manakov-PMD equation to studies of signal
  propagation in optical fibers with randomly varying birefringence}.
\bjournal{J. Lightwave Technology}
\bvolume{15}
\bpages{1735--1746}.
\bptok{imsref}%
\end{barticle}
\endbibitem

%b22 ###
\bibitem{2}
\begin{bmisc}[auto:STB|2012/04/30|08:06:40]
\bauthor{\bsnm{Marty},~\bfnm{Renaud}\binits{R.}}
(\byear{2005}).
\bhowpublished{Probl\`emes d'\'evolution en milieux al\'eatoires: Th\'eor\`emes
  limites, sch\'emas num\'eriques et applications en optique. Ph.D. thesis,
  Univ. Paul Sabatier, Toulouse III.}
\bptok{imsref}%
\end{bmisc}
\endbibitem

%b23 ###
\bibitem{5}
\begin{barticle}[auto:STB|2012/04/30|08:06:40]
\bauthor{\bsnm{Menyuk},~\bfnm{C.~R.}\binits{C.~R.}}
(\byear{1989}).
\btitle{Pulse propagation in an elliptically birefringent Kerr medium}.
\bjournal{IEEE J. Quantum Electronics}
\bvolume{25}
\bpages{2674--2682}.
\bptok{imsref}%
\end{barticle}
\endbibitem

%b24 ###
\bibitem{papanicolaou}
\begin{bincollection}[mr]
\bauthor{\bsnm{Papanicolaou},~\bfnm{G.~C.}\binits{G.~C.}},
  \bauthor{\bsnm{Stroock},~\bfnm{D.}\binits{D.}} \AND
  \bauthor{\bsnm{Varadhan},~\bfnm{S.~R.~S.}\binits{S.~R.~S.}}
(\byear{1977}).
\btitle{Martingale approach to some limit theorems}.
In \bbooktitle{Papers from the {D}uke {T}urbulence {C}onference ({D}uke
  {U}niv., {D}urham, {N}.{C}., 1976), {P}aper {N}o. 6}.
\bseries{Duke Univ. Math. Ser.}
\bvolume{III}
\bpages{ii+120 pp.}
\bpublisher{Duke Univ. Press}, \baddress{Durham, NC}.
\bid{mr={0461684}}
\bptok{imsref}%
\end{bincollection}
\endbibitem

%b25 ###
\bibitem{piatnitski}
\begin{barticle}[mr]
\bauthor{\bsnm{Pardoux},~\bfnm{E.}\binits{E.}} \AND
  \bauthor{\bsnm{Piatnitski},~\bfnm{A.~L.}\binits{A.~L.}}
(\byear{2003}).
\btitle{Homogenization of a nonlinear random parabolic partial differential
  equation}.
\bjournal{Stochastic Process. Appl.}
\bvolume{104}
\bpages{1--27}.
\bid{doi={10.1016/S0304-4149(02)00221-1}, issn={0304-4149}, mr={1956470}}
\bptok{imsref}%
\end{barticle}
\endbibitem

%b26 ###
\bibitem{Stroock}
\begin{bbook}[mr]
\bauthor{\bsnm{Stroock},~\bfnm{Daniel~W.}\binits{D.~W.}} \AND
  \bauthor{\bsnm{Varadhan},~\bfnm{S.~R.~Srinivasa}\binits{S.~R.~S.}}
(\byear{1979}).
\btitle{Multidimensional Diffusion Processes}.
\bseries{Grundlehren der Mathematischen Wissenschaften [Fundamental Principles
  of Mathematical Sciences]}
\bvolume{233}.
\bpublisher{Springer}, \baddress{Berlin}.
\bid{mr={0532498}}
\bptnote{check year}%
\bptok{imsref}%
\end{bbook}
\endbibitem

%b27 ###
\bibitem{Sussman}
\begin{barticle}[mr]
\bauthor{\bsnm{Sussmann},~\bfnm{H{\'e}ctor~J.}\binits{H.~J.}}
(\byear{1978}).
\btitle{On the gap between deterministic and stochastic ordinary differential
  equations}.
\bjournal{Ann. Probab.}
\bvolume{6}
\bpages{19--41}.
\bid{mr={0461664}}
\bptok{imsref}%
\end{barticle}
\endbibitem

%b28 ###
\bibitem{31}
\begin{barticle}[auto:STB|2012/04/30|08:06:40]
\bauthor{\bsnm{Wai},~\bfnm{P.~K.~A.}\binits{P.~K.~A.}},
  \bauthor{\bsnm{Kath},~\bfnm{W.~L.}\binits{W.~L.}} \AND
  \bauthor{\bsnm{Menyuk},~\bfnm{C.~R.}\binits{C.~R.}}
(\byear{1997}).
\btitle{Nonlinear polarization mode dispersion in optical fibers with randomly
  varying birefringence}.
\bjournal{J. Opt. Soc. Amer. A}
\bvolume{14}
\bpages{2697--2979}.
\bptok{imsref}%
\end{barticle}
\endbibitem

%b29 ###
\bibitem{4}
\begin{barticle}[auto:STB|2012/04/30|08:06:40]
\bauthor{\bsnm{Wai},~\bfnm{P.~K.~A.}\binits{P.~K.~A.}} \AND
  \bauthor{\bsnm{Menyuk},~\bfnm{C.~R.}\binits{C.~R.}}
(\byear{1994}).
\btitle{Polarization decorrelation in optical fibers with randomly varying
  birefringence}.
\bjournal{Optics Letters}
\bvolume{19}
\bpages{1517--1519}.
\bptok{imsref}%
\end{barticle}
\endbibitem

%b30 ###
\bibitem{1}
\begin{barticle}[auto:STB|2012/04/30|08:06:40]
\bauthor{\bsnm{Wai},~\bfnm{P.~K.~A.}\binits{P.~K.~A.}} \AND
  \bauthor{\bsnm{Menyuk},~\bfnm{C.~R.}\binits{C.~R.}}
(\byear{1994}).
\btitle{Polarization evolution and dispersion in fibers with spatially varying
  birefringence}.
\bjournal{J. Opt. Soc.}
\bvolume{11}
\bpages{1288--1296}.
\bptok{imsref}%
\end{barticle}
\endbibitem

%b31 ###
\bibitem{3}
\begin{barticle}[auto:STB|2012/04/30|08:06:40]
\bauthor{\bsnm{Wai},~\bfnm{P.~K.~A.}\binits{P.~K.~A.}} \AND
  \bauthor{\bsnm{Menyuk},~\bfnm{C.~R.}\binits{C.~R.}}
(\byear{1996}).
\btitle{Polarization mode dispersion, decorrelation, and diffusion in optical
  fibers with randomly varying birefringence}.
\bjournal{J.~Lightwave Technology}
\bvolume{14}
\bpages{148--157}.
\bptok{imsref}%
\end{barticle}
\endbibitem

%b32 ###
\bibitem{yamato}
\begin{barticle}[mr]
\bauthor{\bsnm{Yamato},~\bfnm{Yuiti}\binits{Y.}}
(\byear{1979}).
\btitle{Stochastic differential equations and nilpotent {L}ie algebras}.
\bjournal{Z.~Wahrsch. Verw. Gebiete}
\bvolume{47}
\bpages{213--229}.
\bid{doi={10.1007/BF00535284}, issn={0044-3719}, mr={0523171}}
\bptok{imsref}%
\end{barticle}
\endbibitem

\end{thebibliography}
\end{document}